\documentclass{amsart}
\usepackage{graphicx}
\usepackage{amssymb}
\usepackage{epstopdf}
\usepackage{nicefrac}
\usepackage{enumerate}
\usepackage{float}
\usepackage{color}

\usepackage{tabularx}
\usepackage{hyperref}

\newtheorem{thm}{THEOREM}[section]
\newtheorem{conj}[thm]{CONJECTURE}
\newtheorem{cor}[thm]{COROLLARY}
\newtheorem{defn}[thm]{DEFINITION}

\newtheorem{lemma}[thm]{LEMMA}

\newtheorem{prop}[thm]{PROPOSITION}

\newtheorem{remark}[thm]{REMARK}

\newcommand{\ds}{\displaystyle}

\newcommand{\ceil}[1]{{\left\lceil #1 \right\rceil}}

\newcommand{\F}{{\mathcal F}} 
\newcommand{\G}{\Gamma}

\newcommand{\e}{{\epsilon}} 
\newcommand{\eU}{{\epsilon_{\cU}}} 
\newcommand{\eFU}{{\epsilon^{\F}_{\cU}}} 
\newcommand{\eTU}{{\epsilon^{\cT}_{\cU}}} 

\newcommand{\ve}{\varepsilon}
\newcommand{\eM}{{\epsilon_{\fM}}} 

\newcommand{\dF}{d_{\F}} 
\newcommand{\dFU}{\delta^{\F}_{\cU}} 
\newcommand{\dTU}{\delta^{\cT}_{\cU}} 
\newcommand{\dT}{\delta_{\cT}} 
\newcommand{\dM}{d_{\fM}} 

\newcommand{\lF}{{\lambda_{\mathcal F}}} 
\newcommand{\rp}{\rho_{\pi}} 
\newcommand{\rt}{\rho_{\tau}} 

\newcommand{\cGF}{\cG_{\F}} 
\newcommand{\cRF}{\cR_{\F}} 
\newcommand{\GF}{\Gamma_{\F}} 

\newcommand{\Dom}{{\rm Dom}} %

\newcommand{\diamX}{{\rm diam}_{\fX}} %
\newcommand{\diamM}{{\rm diam}_{\fM}} %
\newcommand{\dX}{d_{\fX}} %

\newcommand{\kwW}{{\kappa_w^W}} 

\newcommand{\whvarp}{{\widehat \varphi}}
\newcommand{\whfU}{\widehat{\mathfrak{U}}}

\newcommand{\wtgamma}{{\widetilde{\gamma}}}

\newcommand{\whf}{{\widehat f}}

\newcommand{\whK}{{\widehat K}}

\newcommand{\whU}{{\widehat U}}

\newcommand{\whbeta}{\widehat{\beta}}

\newcommand{\whdelta}{\widehat{\delta}}

\newcommand{\oU}{{\overline{U}}}

\newcommand{\bT}{{\bf T}}

\newcommand{\mG}{{\mathbb G}}

\newcommand{\mR}{{\mathbb R}}
\newcommand{\mS}{{\mathbb S}}

\newcommand{\mZ}{{\mathbb Z}}

\newcommand{\cA}{{\mathcal A}}

\newcommand{\cC}{{\mathcal C}}

\newcommand{\cG}{{\mathcal G}}
\newcommand{\cH}{{\mathcal H}}
\newcommand{\cI}{{\mathcal I}}
\newcommand{\cJ}{{\mathcal J}}

\newcommand{\cN}{{\mathcal N}}
\newcommand{\cO}{{\mathcal O}}
\newcommand{\cP}{{\mathcal P}}
\newcommand{\cQ}{{\mathcal Q}}
\newcommand{\cR}{{\mathcal R}}
\newcommand{\cS}{{\mathcal S}}
\newcommand{\cT}{{\mathcal T}}
\newcommand{\cU}{{\mathcal U}}
\newcommand{\cV}{{\mathcal V}}
\newcommand{\cW}{{\mathcal W}}
\newcommand{\cX}{{\mathcal X}}

\newcommand{\fB}{{\mathfrak{B}}}

\newcommand{\fM}{{\mathfrak{M}}}
\newcommand{\fN}{{\mathfrak{N}}}

\newcommand{\fT}{{\mathfrak{T}}}
\newcommand{\fU}{{\mathfrak{U}}}

\newcommand{\fX}{{\mathfrak{X}}}

\newcommand{\fZ}{{\mathfrak{Z}}}
\newcommand{\whfN}{\widehat{{\mathfrak{N}}}}

\input xy
\xyoption {all}

\newcommand{\vp}{{\varphi}}

\begin{document}

 \begin{abstract}
In this work, we develop   shape expansions  of minimal matchbox manifolds without holonomy, in terms of branched manifolds formed from their leaves.
Our approach is based on the method of coding the holonomy groups for the foliated spaces,  to define leafwise regions which are transversely stable and are adapted to the foliation dynamics. Approximations are obtained by collapsing appropriately chosen neighborhoods onto these regions along a ``transverse Cantor foliation''. The existence of the ``transverse Cantor foliation'' allows us to generalize standard techniques known for Euclidean and fibered cases to arbitrary matchbox manifolds with Riemannian leaf geometry and without holonomy. The transverse Cantor foliations used here are constructed by purely intrinsic and topological  means, as we do not assume that our matchbox manifolds are embedded into a smooth foliated manifold, or a smooth manifold. 
\end{abstract}

\title{Shape of matchbox manifolds}


\author{Alex Clark}
\thanks{AC and OL supported in part by  EPSRC grant EP/G006377/1}
\address{Alex Clark, Department of Mathematics, University of Leicester, University Road, Leicester LE1 7RH, United Kingdom}
\email{adc20@le.ac.uk}

\author{Steven Hurder}
\address{Steven Hurder, Department of Mathematics, University of Illinois at Chicago, 322 SEO (m/c 249), 851 S. Morgan Street, Chicago, IL 60607-7045}
\email{hurder@uic.edu}

\author{Olga Lukina}
\address{Olga Lukina, Department of Mathematics, University of Illinois at Chicago, 322 SEO (m/c 249), 851 S. Morgan Street, Chicago, IL 60607-7045}
\email{lukina@uic.edu}

\thanks{Version date: August 15, 2013; revised February 5, 2014}

\date{}



\maketitle

\section{Introduction} \label{sec-intro}

In this work, we consider    topological spaces $\fM$ which are continua; that is,   compact, connected metric spaces.  We assume that $\fM$ has the additional structure of a codimension-zero foliated space, so are \emph{matchbox manifolds}. 
The path-connected components of $\fM$ form the leaves of a foliation $\F$   of dimension $n \geq 1$. The precise definitions are given in Section~\ref{sec-concepts} below. Matchbox manifolds arise naturally in the study of exceptional minimal sets for foliations of compact manifolds, and as the tiling spaces associated to  repetitive, aperiodic tilings of Euclidean space $\mR^n$ which have      finite local complexity. They also arise  in some aspects of group representation theory and index theory for leafwise elliptic operators for foliations, as discussed in the books \cite{CandelConlon2000,MS2006}.

The class of Williams solenoids and the results about these spaces provide one motivation for this work. Recall that an expanding attractor $\Lambda$  for an  Axiom A diffeomorphism $f \colon M \to M$ of a  compact manifold is a     continuum; that is, a compact, connected metric space. 
  Williams  developed a structure theory for these spaces in his seminal works  \cite{Williams1967,Williams1970,Williams1974}.  
  The hyperbolic splitting of the tangent bundle to $TM$ along $\Lambda$ yields a foliation of the space $\Lambda$ by leaves of the expanding foliation for $f$, and the contracting foliation for $f$ gives a transverse foliation on an open neighborhood   $\Lambda \subset U \subset M$.  Williams used this additional structure on a neighborhood of $\Lambda$ to obtain a ``presentation'' of $\Lambda$  as an inverse limit of ``branched $n$-manifolds'', $\whf \colon M_0 \to M_0$, where $n$ is the dimension of the expanding bundle for $f$, and the map $f$ induces the map $\whf$ between the approximations. The notion of a $1$-dimensional branched manifold is easiest to define, as the branches are required only to meet each other at disjoint vertices. In higher dimensions, the definition of branched manifolds becomes more subtle, and especially for the ``transversality condition'' imposed on the cell attachment maps.  The  spaces $\Lambda$ with this structure  are called \emph{Williams solenoids} in the literature.  
The topological properties of the approximating map $\whf \colon M_0 \to M_0$ is used to study the   the dynamical system defined by $f$, for example as discussed in the works \cite{FJ1981,Jones1983,ShubSullivan1975,SullivanWilliams1976} and others.

    The  Riemann surface laminations introduced by Sullivan   \cite{Sullivan1988}  are compact topological spaces locally homeomorphic to a complex disk times a Cantor set, and a similar notion is used by Lyubich and Minsky in \cite{LM1997}, and Ghys in \cite{Ghys1999}. These are also well-known examples of matchbox manifolds.

Associated to the foliation $\F$ of a matchbox manifold $\fM$ is a compactly generated pseudogroup $\cGF$    acting on a ``transverse'' totally disconnected  space $\fX$, which determines the transverse dynamical properties of $\F$. The first two authors showed in  \cite{ClarkHurder2013} that if the action of $\cGF$ on $\fX$ is equicontinuous, then $\fM$ is homeomorphic to a weak solenoid (in the sense of \cite{McCord1965,Schori1966,FO2002}). Equivalently,  the shape of $\fM$ is defined by a tower of proper covering maps between compact manifolds of dimension $n$.

The purpose of this work is to  study the shape properties of   an arbitrary  matchbox manifold  $\fM$, without the assumption that the dynamics of the associated action of $\cGF$ is equicontinuous, though with the assumption  that $\fM$ is minimal;  that is, that every leaf of $\F$ is dense. We also assume that $\F$ is without holonomy. Our main result shows that all such spaces have an analogous structure as that of a William solenoid.  An important difference between the case of Williams solenoids, and the general case we consider, is that the tower of approximations is not in general defined by a single map, but uses a sequence of maps between compact branched manifolds.

  A  \emph{presentation} of a space $\Omega$   is a collection of continuous maps $\cP = \{ p_{\ell} \colon M_{\ell} \to M_{\ell -1} \mid \ell \geq 1\}$, where each $M_{\ell}$ is a connected compact branched $n$-manifold, and each   $p_{\ell} \colon M_{\ell} \to M_{\ell -1}$  is a proper surjective map of branched manifolds, as defined in Section~\ref{sec-Cantorprojections}.
It is assumed that there is given a homeomorphism $h$ between  $\Omega$ and   the inverse limit space  defined by 
\begin{equation}\label{eq-presentationinvlim}
  \cS_{\cP} \equiv \lim_{\longleftarrow} ~ \{ p_{\ell} \colon M_{\ell} \to M_{\ell -1}\} ~ .
\end{equation}
 A \emph{Vietoris solenoid}  \cite{Vietoris1927}  is   a $1$-dimensional solenoid, where each $M_{\ell}$ is a circle, and each $p_{\ell} \colon \mS^1 \to \mS^1$ is a covering map of degree greater than $1$. More generally, if  each $M_{\ell}$ is a compact manifold and each $p_{\ell}$ is a proper covering map, 
 then we say that      $\cS_{\cP}$ is a \emph{weak solenoid},  as discussed in \cite{McCord1965,Schori1966,FO2002}.   

Here is our main result.
\begin{thm}\label{thm-main1}
Let  $\fM$ be a minimal matchbox manifold   of dimension $n$. Assume that the foliation $\F$ of $\fM$ is without holonomy. Then  there exists a presentation  $\ds \cP = \{ p_{\ell} \colon M_{\ell} \to M_{\ell -1}  \mid \ell \geq \ell_0\}$ where each $M_{\ell}$ is a triangulated branched $n$-manifold and  each $p_{\ell}$ is a    proper  simplicial map, so  that $\fM$   is   homeomorphic to the inverse limit of the   system of maps they define, 
  \begin{equation} \label{eq-mainthmbonding}
 h \colon  \fM ~  \cong  ~  \cS_{\cP} = \varprojlim ~ \{ p_{\ell} \colon M_{\ell} \to M_{\ell -1} \mid \ell \geq \ell_0\}. 
  \end{equation}
\end{thm}
The   hypothesis that $\F$ is without holonomy may possibly be removed, using a more general definition of the bonding maps allowed. The work by 
 Benedetti and  Gambaudo in \cite{BG2003} gives a model for such a result, in the case of tiling spaces associated with an action of a connected Lie group $G$. Also, the results of the first two authors in \cite{ClarkHurder2013} applies to equicontinuous matchbox manifolds which may have non-trivial leafwise holonomy. On the other hand, a generalization of Theorem~\ref{thm-main1} to arbitrary minimal matchbox manifolds would yield solutions to several classical problems in the study of exceptional minimal sets for codimension-one foliations, 
 for example as listed in \cite{Hurder2002,Hurder2006}, so an extension of Theorem~\ref{thm-main1} to this generality would be important, but is expected to involve additional subtleties.

 In this work,   no assumptions are made on the geometry and topology of the leaves in $\fM$, beyond that they are complete Riemannian manifolds. In order to construct the branched manifold approximations to $\fM$,  
  we use   a technique based on   ``dynamical codings'' for the orbits of the pseudogroup associated to the foliation, a method which generalizes that used by Gambaudo and Martens  \cite{GambaudoMartens2006} for flows.   This coding method can also be seen as extending   a technique used by  Thomas in his study of equicontinuous flows in \cite{EThomas1973}. The first two  authors used this coding approach in their study of equicontinuous actions of groupoids on Cantor sets in \cite{ClarkHurder2013}, and it also appears implicitly in the work by Forrest \cite{Forrest2000} in his study of minimal actions of $\mZ^n$ on Cantor sets. 
 The coding method we use is  more formal, so applies in complete  generality. The other ingredient required is the existence of a transverse Cantor foliation in the generality we consider. We construct these foliations in this work, extending the   results of the authors' paper \cite{CHL2013a}. The existence of a compatible transverse Cantor foliation $\cH$ on $\fM$  is fundamental for defining the branched manifold approximations.

As mentioned above,      tiling spaces provide an important class of examples for which Theorem~\ref{thm-main1} applies.  Let $\bT$ be  a tiling of $\mR^n$ which is     repetitive, aperiodic, and has    finite local complexity. The tiling space $\Omega_{\bT}$   is defined as the   closure of the  set of tilings obtained via the translation action of $\mR^n$,
in a suitable metric topology \cite{AP1998,BBG2006,FHK2002,Sadun2008}. The assumptions   imply  that $\Omega_{\bT}$ is locally homeomorphic to a disk in $\mathbb{R}^n$ times a Cantor set, so that $\Omega_{\bT}$ is a matchbox manifold. (See \cite{FS2009} for a discussion of variants of this result.)
The seminal  result of Anderson and Putnam in \cite{AP1998} showed that under these assumptions,  and under the assumption that the tiling $\bT$ is defined by an expanding substitution, then   $\Omega_{\bT}$  is homeomorphic to the inverse limit of a tower of branched flat $n$-manifolds. The approach of Anderson and Putnam, and the alternative approaches  to generalizing their method and extending the conclusion to all   tilings  by G\"{a}hler and Sadun \cite{Sadun2003},
 and  in   works \cite{ALM2011,BBG2006,BG2003,GambaudoMartens2006,LR2013},
 used variations on a technique   called ``collaring'' or  a method called ``inflation'' of a Voronoi tessellation of $\mR^n$,
 to obtain the branched manifold approximations for $\Omega_{\bT}$. The more recent work \cite{BDHS2010} gives a version of this method closest to the approach we take in this paper.  The work here is an important step in developing the general theory of point patterns and tilings for general, non-Euclidean spaces. As described in Senechal~\cite{Sen1995}, the theory of point patterns and tilings is directly connected to the theory of quasicrystals, and so our work can be seen as a step towards understanding possible quasicrystalline structures in non-Euclidean spaces.

 We discuss next the contents of the rest of this work, which culminates in  the proof of Theorem~\ref{thm-main1}. 
Section~\ref{sec-concepts} introduces the basic definitions of a matchbox manifold, and Section~\ref{sec-holonomy} introduces the definitions and properties of the holonomy pseudogroups for    matchbox manifolds, as developed in \cite{ClarkHurder2013} and \cite{CHL2013a}. Section~\ref{sec-voronoi}  discusses the Voronoi tessellations of the leaves associated to transversal clopen sets, then Section~\ref{sec-induced} develops some fundamental properties of restricted pseudogroups.

A coding of the orbits of the   action of the holonomy pseudogroup of the foliation $\F$ on $\fM$ is  developed in Section~\ref{sec-partialcoding}, and  used to construct ``dynamically defined''  nested sequences of clopen coverings $\cV_{\ell}$ for $\ell \geq \ell_0$ of a Cantor transversal to $\F$, which are ``centered'' on the transverse orbit defined by a fixed leaf $L_0$. The study of the dynamics of the induced action on a descending chain of clopen neighborhoods  is a standard technique in almost all  approaches to the study of the dynamics of minimal actions on a Cantor set. It corresponds to the initial steps in the construction of a Kakutani tower for the measurable theory of group actions.

In Section~\ref{sec-reebslabs}, we associate to each  clopen covering of the transversal space $\fX$  a collection of open sets, called \emph{Reeb neighborhoods}  
which form a covering of $\fM$.  The results of \cite{CHL2013a}  show that   each such Reeb neighborhood  projects  along a transverse Cantor foliation    to its compact base.  
Section~\ref{sec-compatiblecantor} uses  the results in  \cite{CHL2013a}   to show that   these bi-foliated structures on the Reeb neighborhoods can be chosen to be  compatible on overlaps.  
In Section~\ref{sec-Cantorprojections}, the proof of Theorem~\ref{thm-main1} is   completed, by showing that 
the  Reeb neighborhoods collapsed along the Cantor transversal foliation, can be identified (or ``glued'')  each to the other on their overlaps, to obtain a 
  branched  $n$-manifold. By iterating this process for a   chain of clopen covers for $\fX$,    we obtain a tower of maps whose inverse limit is homeomorphic to $\fM$.

Note that the works \cite{ALM2011,LR2013} follow  a similar approach to the above, but using an intuitive extension of the    zooming technique for Voronoi cells   as in \cite{BBG2006,BG2003}. The geometry and topology of a generic leaf $L_0$   can be quite  complicated, and not at all as intuitively straightforward as for the case of Euclidean spaces. Hence the geometry and topology  of the Voronoi cells which decompose $L_0$  may be similarly non-intuitive. Some of the subtleties for    Voronoi partitions at large scale in non-Euclidean manifolds  are discussed, for example, in the work \cite{LeibonLetscher2000}. The more formal approach using codings that we use in this work  avoids the need to control the geometry of the Voronoi cells introduced.   The reader will also find the care taken with the domains of holonomy maps makes the proofs of some results quite  tedious at times, but a careful consideration of examples shows that these concerns are required.

This work is part of a program to generalize the results of the thesis of Fokkink \cite{Fokkink1991}, started during a visit by the authors to the University of Delft in August 2009. The papers \cite{ClarkHurder2013,CHL2013a} are the initial results of this study, and this work is the preparation for the forthcoming paper \cite{CHL2013c} which completes the program. The authors would like to thank Robbert Fokkink   for the invitation to meet in Delft, and the University of Delft for its generous support for the visit. The authors' stay in Delft  was also supported by a   travel grant No.~040.11.132 from the \emph{Nederlandse Wetenschappelijke Organisatie}.

\section{Foliated spaces and matchbox manifolds} \label{sec-concepts}

 In this section, we give the precise definitions and results for  matchbox manifolds, as required for this work.  More detailed  discussions with examples can be found in the books of Candel and Conlon \cite[Chapter 11]{CandelConlon2000} and Moore and Schochet \cite[Chapter 2]{MS2006}. In the following, details of proofs are omitted whenever possible, as they are presented in detail in these two texts and  the  paper \cite{ClarkHurder2013}. 

\begin{defn} \label{def-fs}
A \emph{foliated space of dimension $n$} is a   continuum $\fM$, such that  there exists a compact   metric space $\fX$, and
for each $x \in \fM$ there is a compact subset $\fT_x \subset \fX$, an open subset $U_x \subset \fM$, and a homeomorphism $\vp_x \colon \oU_x \to [-1,1]^n \times \fT_x$ defined on the closure of $U_x$ in $\fM$,  such that   $\vp_x(x) = (0, w_x)$ where $w_x \in int(\fT_x)$. 
Moreover, it is assumed that each $\vp_x$  admits an extension to a foliated homeomorphism
$\whvarp_x \colon \whU_x \to (-2,2)^n \times \fT_x$ where $\oU_x \subset \whU_x$ is an open neighborhood.
\end{defn}
The subspace  $\fT_x$ of $\fX$ is called   the \emph{local transverse model} at $x$.

Let $\pi_x \colon \oU_x \to \fT_x$ denote the composition of $\vp_x$ with projection onto the second factor.

For $w \in \fT_x$ the set $\cP_x(w) = \pi_x^{-1}(w) \subset \oU_x$ is called a \emph{plaque} for the coordinate chart $\vp_x$. We adopt the notation, for $z \in \oU_x$, that $\cP_x(z) = \cP_x(\pi_x(z))$, so that $z \in \cP_x(z)$. Note that each plaque $\cP_x(w)$ for $w \in \fT_x$ is given the topology so that the restriction $\vp_x \colon \cP_x(w) \to [-1,1]^n \times \{w\}$ is a homeomorphism. Then $int (\cP_x(w)) = \vp_x^{-1}((-1,1)^n \times \{w\})$.

Let $U_x = int (\oU_x) = \vp_x^{-1}((-1,1)^n \times int(\fT_x))$.
Note that if $z \in U_x \cap U_y$, then $int(\cP_x(z)) \cap int( \cP_y(z))$ is an open subset of both
$\cP_x(z) $ and $\cP_y(z)$.
The collection of sets
$$\cV = \{ \vp_x^{-1}(V \times \{w\}) \mid x \in \fM ~, ~ w \in \fT_x ~, ~ V \subset (-1,1)^n ~ {\rm open}\}$$
forms the basis for the \emph{fine topology} of $\fM$. The connected components of the fine topology are called \emph{leaves}, and define the foliation $\F$ of $\fM$.
For $x \in \fM$, let $L_x \subset \fM$ denote the leaf of $\F$ containing $x$.

 Definition~\ref{def-fs} does not impose any smoothness conditions on the leaves of $\F$, so they may just be topological manifolds. The next definition imposes a \emph{uniform} smoothness condition on the leaves. 
\begin{defn} \label{def-sfs}
A \emph{smooth foliated space} is a foliated space $\fM$ as above, such that there exists a choice of local charts $\vp_x \colon \oU_x \to [-1,1]^n \times \fT_x$ such that for all $x,y \in \fM$ with $z \in U_x \cap U_y$, there exists an open set $z \in V_z \subset U_x \cap U_y$ such that $\cP_x(z) \cap V_z$ and $\cP_y(z) \cap V_z$ are connected open sets, and the composition
$$\psi_{x,y;z} \equiv \vp_y \circ \vp_x ^{-1}\colon \vp_x(\cP_x (z) \cap V_z) \to \vp_y(\cP_y (z) \cap V_z)$$
is a smooth map, where $\vp_x(\cP_x (z) \cap V_z) \subset \mR^n \times \{w\} \cong \mR^n$ and $\vp_y(\cP_y (z) \cap V_z) \subset \mR^n \times \{w'\} \cong \mR^n$. The leafwise transition maps $\psi_{x,y;z}$ are assumed to depend continuously on $z$ in the $C^{\infty}$-topology on maps between subsets of $\mR^n$.
\end{defn}

A map $f \colon \fM \to \mR$ is said to be \emph{smooth} if for each flow box
$\vp_x \colon \oU_x \to [-1,1]^n \times \fT_x$ and $w \in \fT_x$ the composition
$y \mapsto f \circ \vp_x^{-1}(y, w)$ is a smooth function of $y \in (-1,1)^n$, and depends continuously on $w$ in the $C^{\infty}$-topology on maps of the plaque coordinates $y$. As noted in \cite{MS2006} and \cite[Chapter 11]{CandelConlon2000}, this allows one to define leafwise smooth partitions of unity, vector bundles, and tensors for smooth foliated spaces. In particular, one can define leafwise smooth Riemannian metrics. We recall a standard result, whose proof for foliated spaces can be found in \cite[Theorem~11.4.3]{CandelConlon2000}.
\begin{thm}\label{thm-riemannian}
Let $\fM$ be a smooth foliated space. Then there exists a leafwise Riemannian metric for $\F$, such that for each $x \in \fM$, the leaf $L_x$ inherits the structure of a complete Riemannian manifold with bounded geometry, and the Riemannian geometry of $L_x$ depends continuously on $x$. 
\end{thm}

Bounded geometry implies, for example, that for each $x \in \fM$, there is a leafwise exponential map
$\exp^{\F}_x \colon T_x\F \to L_x$ which is a surjection, and the composition $\exp^{\F}_x \colon T_x\F \to L_x \subset \fM$ depends continuously on $x$ in the compact-open topology on maps.

 \begin{defn} \label{def-mm}
A \emph{matchbox manifold} is a continuum with the structure of a
smooth foliated space $\fM$, such that  the transverse model space $\fX$ is totally disconnected, and for each $x \in \fM$, $\fT_x \subset \fX$ is a clopen (closed and open) subset.  A matchbox manifold $\fM$ is \emph{minimal} if every leaf of $\F$ is dense.
\end{defn}
 
    Intuitively, a $1$-dimensional matchbox manifold   has local coordinate charts
  which are homeomorphic to a ``box of matches'',  which gave   rise to this  terminology in the works   \cite{AHO1991,AO1995,AM1988}.

The  maximal path-connected components of $\fM$ define the leaves of a foliation $\F$ of $\fM$. 
All matchbox manifolds  are assumed to be smooth, with a   leafwise Riemannian metric, and  metric $\dM$ on $\fM$.

\subsection{Metric properties and regular covers}\label{subsec-metricregular}
We   introduce some local metric considerations for a matchbox manifold, and give  the definition of a \emph{regular covering} of $\fM$.  
 
For $x \in \fM$ and $\e > 0$, let $D_{\fM}(x, \e) = \{ y \in \fM \mid \dM(x, y) \leq \e\}$ be the closed $\e$-ball about $x$ in $\fM$, and $B_{\fM}(x, \e) = \{ y \in \fM \mid \dM(x, y) < \e\}$ the open $\e$-ball about $x$.

Similarly, for $w \in \fX$ and $\e > 0$, let $D_{\fX}(w, \e) = \{ w' \in \fX \mid d_{\fX}(w, w') \leq \e\}$ be the closed $\e$-ball about $w$ in $\fX$, and $B_{\fX}(w, \e) = \{ w' \in \fX \mid d_{\fX}(w, w') < \e\}$ the open $\e$-ball about $w$.

Each leaf $L \subset \fM$ has a complete path-length metric, induced from the leafwise Riemannian metric:
$$\dF(x,y) = \inf \left\{\| \gamma\| \mid \gamma \colon [0,1] \to L ~{\rm is ~ piecewise ~~ C^1}~, ~ \gamma(0) = x ~, ~ \gamma(1) = y ~, ~ \gamma(t) \in L \quad \forall ~ 0 \leq t \leq 1\right\}$$
  where $\| \gamma \|$ denotes the path-length of the piecewise $C^1$-curve $\gamma(t)$. If $x,y \in \fM$   are not on the same leaf, then set $\dF(x,y) = \infty$. 
  
  For each $x \in \fM$ and $r > 0$, let $D_{\F}(x, r) = \{y \in L_x \mid \dF(x,y) \leq r\}$ be the closed leafwise ball.

For each $x \in \fM$, the  {Gauss Lemma} implies that there exists $\lambda_x > 0$ such that $D_{\F}(x, \lambda_x)$ is a \emph{strongly convex} subset for the metric $\dF$. That is, for any pair of points $y,y' \in D_{\F}(x, \lambda_x)$ there is a unique shortest geodesic segment in $L_x$ joining $y$ and $y'$ and  contained in $D_{\F}(x, \lambda_x)$ (cf. \cite[Chapter 3, Proposition 4.2]{doCarmo1992}, or \cite[Theorem 9.9]{Helgason1978}). Then for all $0 < \lambda < \lambda_x$ the disk $D_{\F}(x, \lambda)$ is also strongly convex. As $\fM$ is compact and the leafwise metrics have uniformly bounded geometry, we obtain:
\begin{lemma}\label{lem-stronglyconvex}
There exists $\lF > 0$ such that for all $x \in \fM$, $D_{\F}(x, \lF)$ is strongly convex.
\end{lemma}

The following proposition summarizes results in \cite[Sections 2.1 - 2.2]{ClarkHurder2013}.
\begin{prop}\label{prop-regular} 
For a smooth foliated space $\fM$, given $\eM > 0$, there exist  constants  $\lF>0$ and $0< \dFU < \lF/5$, and  a covering of $\fM$ by foliation  charts 
$\ds \left\{\vp_i \colon \oU_i \to [-1,1]^n \times \fT_i \mid 1 \leq i \leq \nu \right\}$
 with the following properties: For each  $1 \leq i \leq \nu$, let   $\pi_i = \pi_{x_i} \colon \oU_i \to \fT_i$ be the projection, then
\begin{enumerate}
\item  Interior: $U_i \equiv int(\oU_i) = \vp_i^{-1}\left( (-1,1)^n \times \fT_i\right)$
\item Locality: for   $x_i \equiv \vp_i^{-1}(w_i, 0) \in \fM$,    $\oU_i \subset B_{\fM}(x_i, \eM)$ where  $w_i = \pi_i(x_i)$.
\end{enumerate}
For $z \in \oU_i$, the \emph{plaque} of the chart $\vp_i$ through $z$ is denoted by $\cP_i(z) = \cP_i(\pi_i(z)) \subset \oU_i$. 
\begin{enumerate}\setcounter{enumi}{2}
\item Convexity:  the plaques of $\vp_i$ are   strongly convex subsets for the leafwise metric. 
\item  \label{item-uniform}  Uniformity:  for   $w \in \fT_i$ let $x_{w} = \vp_{x_i}^{-1}(0 , w)$, then 
$$
D_{\F}(x_{w} , \dFU/2) ~ \subset ~ \cP_i(w) ~ \subset ~ D_{\F}(x_{w} , \dFU) 
$$
\item \label{item-clopen} The projection $\pi_i(U_i \cap U_j) = \fT_{i,j} \subset \fT_i$ is a clopen subset for all $1 \leq i, j \leq \nu$. 
\end{enumerate}
A \emph{regular foliated covering} of $\fM$ is one that satisfies   these conditions. 
\end{prop}

We assume in the following that a   regular foliated covering  of $\fM$ as in Proposition~\ref{prop-regular}  has been chosen.
 Let $\cU = \{U_{1}, \ldots , U_{\nu}\}$ denote the corresponding open covering of $\fM$.
We can   assume without loss of generality,  that the spaces $\fT_i$ form a \emph{disjoint  clopen covering} of $\fX$, so that 
  $\ds \fX = \fT_1 \ \dot{\cup} \cdots \dot{\cup} \ \fT_{\nu}$.

Let $\eU > 0$ be a Lebesgue number for $\cU$. That is, given any $z \in \fM$ there exists some index $1 \leq i_z \leq \nu$ such that the open metric ball $B_{\fM}(z, \eU) \subset U_{i_z}$. Also, introduce a form of ``leafwise Lebesgue number'', defined by
\begin{equation}\label{eq-leafdiam}
\eFU = \min \left\{ \eFU(y) \mid ~ \forall ~ y \in \fM \right\} ~ , ~ \eFU(y) = \max \left\{ \e \mid  ~ D_{\F}(y, \e) \subset D_{\fM}(y, \eU/4)\right\}.
\end{equation}
Thus, for all $y \in \fM$, $D_{\F}(y, \eFU) \subset D_{\fM}(y, \eU/4)$.
Note that for all $r > 0$ and $z' \in D_{\F}(z, \eFU)$, the triangle inequality implies that
$D_{\fM}(z', r) \subset D_{\fM}(z, r + \eU/4)$.
Note that  for $y \in \fM$, we have 
$$ D_{\F}(y, \eFU) \subset D_{\fM}(y, \eU/4) \subset D_{\fM}(y, \eU) \subset U_{i_y}, ~\cP_{i_y}(\pi_{i_y}(y)) \subset D_{\F}(y , \dFU)$$
which   imply that $\eFU < \dFU$, and hence $\dFU/\eFU > 1$. This last inequality  will be used in Section~\ref{sec-induced}.

For $1 \leq i \leq \nu$, let  $ \lambda_i \colon \oU_i \to [-1,1]^n$ be the projection, so that for each $z \in U_i$   the restriction $\lambda_i \colon \cP_i(z) \to [-1,1]^n$ is    a smooth coordinate system on the plaque $\cP_i(z)$.

For each $1 \leq i \leq \nu$, the set $\cT_i =  \vp_i^{-1}(0 , \fT_i)$ is a compact transversal to $\F$. Without loss of generality, we can assume   that the transversals 
$\ds \{ \cT_{1} , \ldots , \cT_{\nu} \}$ are pairwise disjoint in $\fM$, so there exists a constant $0 < \dT <  \dFU$ such that 
\begin{equation}\label{eq-dT}
\dF(x,y) \geq \dT \quad {\rm for} ~ x \ne y ~,   x \in \cT_i ~ , ~ y \in \cT_j ~ , ~ 1 \leq i, j \leq \nu ~.
\end{equation}
In particular, this implies that 
 the centers of disjoint plaques on the same leaf are separated by distance at least $\dT$.
Also, define sections
\begin{equation}\label{eq-taui}
\tau_i \colon \fT_i \to \oU_i ~ , ~ {\rm defined ~ by} ~ \tau_i(\xi) = \vp_i^{-1}(0 , \xi) ~ , ~ {\rm so ~ that} ~ \pi_i(\tau_i(\xi)) = \xi.
\end{equation}
Then $\cT_i = \cT_{x_i}$ is the image of $\tau_i$ and we let $\cT = \cT_1 \cup \cdots \cup \cT_{\nu} \subset \fM$ denote their disjoint union, and $\tau \colon \fX \to \cT$ is defined by the union of the maps $\tau_i$.

 Define the metric  $d_{\fX}$ on $\fX$ via the restriction of $\dM$ to $\cT$, and use the map  $\tau$ to pull it back to $\fX$.

\subsection{Local estimates}\label{subsec-locestimates}

The local projections $\pi_i \colon \oU_i \to \fT_i$ and sections
$\tau_i \colon \fT_i \to \oU_i$ are continuous maps of compact spaces, so admit uniform metric estimates as shown in \cite{ClarkHurder2013}. 
\begin{lemma} \label{lem-modpi}
There exists a    \emph{modulus of continuity} function $\rp$   which is continuous and increasing,     such that:
\begin{equation}\label{eq-modpi}
\forall ~ 1 \leq i \leq \nu ~ \text{and}~ x,y \in \oU_i \quad , \quad \dM(x,y) <\rp(\e) ~ \Longrightarrow ~ d_{\fX}(\pi_i(x), \pi_i(y)) < \e ~.
\end{equation}
\end{lemma}
\proof
Set $\ds\rp(\e) = \min \left\{\e, \min \left\{ \dM(x,y) \mid 1 \leq i \leq \nu ~ , ~ x,y \in \oU_i ~ , ~ d_{\fX}(\pi_i(x), \pi_i(y)) \geq \e\right\}\right\}$.
\endproof

\begin{lemma}  \label{lem-modtau}
There exists a \emph{modulus of continuity} function    $\rt$ which is continuous and increasing,     such that:
\begin{equation}\label{eq-modtau}
\forall ~ 1 \leq i \leq \nu ~ \text{and}~ w, w' \in \fT_i \quad , \quad d_{\fX}(w,w') <\rt(\e) ~ \Longrightarrow ~ \dM(\tau_i(w), \tau_i(w')) < \e ~ .
\end{equation}
\end{lemma}
\proof
Set $\ds\rt(\e) = \min \left\{\e, \min \left\{ d_{\fX}(w,w') \mid 1 \leq i \leq \nu ~ , ~ w,w' \in \fT_i ~ , ~ \dM(\tau_i(w), \tau_i(w')) \geq \e\right\}\right\}$.
\endproof

Next, for each $1 \leq i \leq \nu$, consider the projection map 
$\pi_i \colon \oU_i \to \fT_i$, and define 
$$ \epsilon^{\cT}_{i} = \max \left\{\e \mid  \forall ~ x \in \oU_i ~{\rm such ~ that}  ~ D_{\fM}(x, \eU/2) \subset \oU_i ~ , {\rm then} ~ D_{\fX}(\pi_i(x),\e) \subset \pi_i\left( D_{\fM}(x, \eU/2)\right)\right\}.
$$ 
 The   assumption  that $D_{\fM}(x, \eU/2) \subset \oU_i$   implies that $x$  has distance at least $\eU/2$ from the exterior of  $\oU_i$. Then the projection of the closed ball $D_{\fM}(x, \eU/2)$ to the transversal $\fT_i$ contains an open neighborhood of $\pi_i(x)$ by the continuity of projections, and $\epsilon^{\cT}_{i}$ is the distance from this center $\pi_i(x)$ to the exterior of the projected ball.
Then introduce $\eTU \geq \rt(\eU/2 )$ given by 
 \begin{equation}\label{eq-transdiam}
\eTU = \min \left\{\epsilon^{\cT}_{i} \mid \forall ~ 1 \leq i \leq \nu \right\}.
\end{equation}
  
Finally, we give a plaque-wise estimate on the Riemannian metrics.  
The assumption that the leafwise Riemannian metric on $\fM$ is continuous, means that for each coordinate chart   
  $\vp_i \colon \oU_i \to [-1,1]^n \times \fT_i$, the push-forwards of the Riemannian metric to the slices $[-1,1]^n \times \{w\}$ vary continuously with $w \in \fT_i$.
Let $\| \cdot \|_w$ denote the norm defined on the tangent bundle $T(-1,1)^n \times \{w\}$, and denote a tangent vector by $(\vec{v}, \xi, w)$, for $\vec{v} \in  \mR^n$ and $\xi \in (-1,1)^n$. Let $\| \cdot \|$ denote the Euclidean norm on $\mR^n$. Then by the compactness of $\oU_i$ and the continuity of the metric, for 
$\e> 0$,    there exists $\delta_{\cU}(\e) > 0$ such that the following holds, for each $1 \leq i \leq \nu$:
\begin{equation}\label{eq-metricbounds}
\sup \left \{ \left| \max \left\{ \frac{\|\vec{v}\|_w}{\|\vec{v}\|_{w'}} , \frac{\|\vec{v}\|_{w'}}{\|\vec{v}\|_{w}}\right\} - 1 \right|    \mid d_{\fX}(w,w') \leq  \delta_{\cU}(\e) ~, ~ \vec{v} \in \mR^n ~, ~ \|\vec{v}\| = 1 \right\}  ~ < ~ \e ~ .
\end{equation}

\subsection{Foliated maps}

A map $f \colon \fM \to \fM'$ between foliated spaces is said to be a \emph{foliated map} if the image of each leaf of $\F$ is contained in a leaf of $\F'$. If $\fM'$ is a matchbox manifold, then each leaf of $\F$ is path connected, so its image is path connected, hence must be contained in a leaf of $\F'$.  

A \emph{leafwise path}  is a continuous map $\gamma \colon [0,1] \to \fM$ such that there is a leaf $L$ of $\F$ for which $\gamma(t) \in L$ for all $0 \leq t \leq 1$. 
If $\fM$ is a matchbox manifold and $\gamma \colon [0,1] \to \fM$ is continuous, then   $\gamma$ is a leafwise path. This yields:
\begin{lemma} \label{lem-foliated1}
Let $\fM$ and $\fM'$ be matchbox manifolds, and $h \colon \fM \to \fM'$ a continuous map. Then $h$ maps the leaves of $\F$ to leaves of $\F'$. In  particular, any homeomorphism $h \colon \fM \to \fM$ of a matchbox manifold is a foliated map. \hfill $\Box$
\end{lemma}

\section{Holonomy of foliated spaces} \label{sec-holonomy}

The holonomy pseudogroup of a smooth foliated manifold $(M, \F)$ generalizes the induced discrete dynamical system     associated to a section of a flow. 
 A standard construction in foliation theory (see  \cite{CN1985}, \cite[Chapter 2]{CandelConlon2000} or  \cite{Haefliger1984}) associates to a leafwise path $\gamma$ a holonomy map $h_{\gamma}$. 
The collection of all such maps with initial and ending points on a transversal to $\F$ defines the holonomy pseudogroup $\cGF$ for a matchbox manifold $(\fM, \F)$.  We recall below the ideas and notations,  as  required in the proofs of our main theorems, especially the  delicate issues of domains which must be considered. 
See \cite{ClarkHurder2013,CHL2013a} for a complete discussion of the ideas of this section  and related technical results.

\subsection{Holonomy pseudogroup} 
Let $\cU = \{U_{1}, \ldots , U_{\nu}\}$   be a    regular foliated covering  of $\fM$ as in Proposition~\ref{prop-regular}. 
A pair of indices $(i,j)$, for $1 \leq i,j \leq \nu$, is said to be \emph{admissible} if  $U_i \cap U_j \ne \emptyset$.

For $(i,j)$ admissible, set $\fT_{i,j} = \pi_i(U_i \cap U_j) \subset \fT_i$.  The regular foliated covering assumption  implies that plaques in admissible charts are either disjoint, or have connected intersection. This implies that there is a well-defined transverse change of coordinates homeomorphism $h_{i,j} \colon \fT_{i,j} \to \fT_{j,i}$ with domain $\Dom(h_{i,j}) = \fT_{i,j}$ and range $R(h_{i,j}) = \Dom(h_{j,i}) = \fT_{j,i}$.   By definition they satisfy $h_{i,i} = Id$, $h_{i,j}^{-1} = h_{j,i}$, and if $U_i \cap U_j\cap U_k \ne \emptyset$ then $h_{k,j} \circ h_{j,i} = h_{k,i}$ on their common domain of definition. 
Note that the domain and range of $h_{i,j}$ are clopen subsets of $\fX$ by Proposition~\ref{prop-regular}.(\ref{item-clopen}).

 Recall that for each $1 \leq i \leq \nu$, $\tau_i \colon \fT_i \to \cT_i$ denotes the transverse section for the   chart $U_i$, and $\cT = \cT_1 \cup \cdots \cup \cT_{\nu} \subset \fM$ denotes their disjoint union. Then $\pi \colon \cT \to \fX$ is the coordinate projection restricted to $\cT$, which is a homeomorphism, and $\tau \colon \fX \to \cT$ denotes  its inverse.

The \emph{holonomy pseudogroup} $\cGF$ of $\F$ is the topological pseudogroup modeled on $\fX$ generated by   the elements of   $\cGF^{(1)} = \{h_{j,i} \mid (i,j) ~{\rm admissible}\}$. 
A sequence $\cI = (i_0, i_1, \ldots , i_{\alpha})$ is \emph{admissible} if each pair $(i_{\ell -1}, i_{\ell})$ is admissible  for $1 \leq \ell \leq \alpha$, and the composition
\begin{equation}\label{eq-defholo}
 h_{\cI} = h_{i_{\alpha}, i_{\alpha-1}} \circ \cdots \circ h_{i_1, i_0}
\end{equation}
 has non-empty domain $\Dom(h_{\cI})$, which is defined to be    the  maximal \emph{clopen} subset of $\fT_{i_0}$ for which the compositions are defined.
Given any open subset $U \subset \Dom(h_{\cI})$,  we obtain a new element $h_{\cI} | U \in \cGF$ by restriction.  Introduce
\begin{equation}\label{eq-restrictedgroupoid}
\cGF^* = \left\{ h_{\cI} |  U \mid   \cI ~ {\rm admissible~ and} ~ U \subset \Dom(h_{\cI}) \right\} \subset \cGF ~ .
\end{equation}
 The range of $g = h_{\cI} |  U$ is the open set $R(g) = h_{\cI}(U) \subset \fT_{i_{\alpha}} \subset \fX$. Note that each map $g \in \cGF^*$ admits a
continuous extension $\overline{g} \colon \overline{\Dom(g)} = \overline{U} \to \fT_{i_{\alpha}}$ as $\Dom( h_{\cI})$ is a clopen set for each $\cI$.

\subsection{Plaque-chains} \label{subsec-plaquechains}

Let $\cI = (i_0, i_1, \ldots , i_{\alpha})$ be an  admissible sequence. 
For each $1 \leq \ell \leq \alpha$, set 
$\cI_{\ell} = (i_0, i_1, \ldots, i_{\ell})$, and let $h_{\cI_{\ell}}$ denote the corresponding holonomy map. For $\ell = 0$, let $\cI_0 = (i_0 , i_0)$.
Note that $h_{\cI_{\alpha}} = h_{\cI}$ and $h_{\cI_{0}} = Id \colon \fT_0 \to \fT_0$.

Given $w \in \Dom(h_{\cI})$,  let $x = \tau_{i_0}(w) \in L_{w}$. For each 
$0 \leq \ell \leq \alpha$, set $w_{\ell} = h_{\cI_{\ell}}(w)$ and
$x_{\ell}= \tau_{i_{\ell}}(w_{\ell})$. 
Recall that $\cP_{i_{\ell}}(x_{\ell}) = \cP_{i_{\ell}}(w_{\ell})$, where 
each $\cP_{i_{\ell}}(w_{\ell})$ is a strongly convex subset of the   leaf $L_w$ in the leafwise metric $d_{\F}$. 
 Introduce the   \emph{plaque-chain}
 \begin{equation}\label{eq-plaquechain}
\cP_{\cI}(w) = \{\cP_{i_0}(w_0), \cP_{i_1}(w_1), \ldots , \cP_{i_{\alpha}}(w_{\alpha}) \} ~ .
\end{equation}
Adopt the notation $\cP_{\cI}(x) \equiv \cP_{\cI}(w)$.
  Intuitively, a plaque-chain $\cP_{\cI}(x)$ is a sequence of successively overlapping convex ``tiles'' in $L_{w}$ starting at $x = \tau_{i_0}(w)$, ending at
$y = x_{\alpha} = \tau_{i_{\alpha}}(w_{\alpha})$, where each $\cP_{i_{\ell}}(x_{\ell})$ is ``centered'' on the point $x_{\ell} = \tau_{i_{\ell}}(w_{\ell})$.

 \subsection{Leafwise paths to plaque-chains}  \label{subsec-pathstochains}
 
    Let $\gamma \colon  [0,1] \to \fM$ be a   path. Set $x_0 = \gamma(0) \in U_{i_0}$,   $w = \pi(x_0)$ and $x = \tau(w) \in \cT_{i_0}$.  
Let $\cI$ be an admissible sequence with  $w \in \Dom(h_{\cI})$. We say that $(\cI , w)$ \emph{covers} $\gamma$,  
if there is    a partition $0 = s_0 < s_1 < \cdots < s_{\alpha} = 1$ such that   $\cP_{\cI}(w)$   satisfies
\begin{equation}\label{eq-cover}
\gamma([s_{\ell} , s_{\ell + 1}]) \subset    \cP_{i_{\ell}}(\xi_{\ell})  ~ , ~ 0 \leq \ell < \alpha, ~ {\rm and} ~  \gamma(1) \in   \cP_{i_{\alpha}}(\xi_{\alpha}) .
\end{equation}
  For  a   path $\gamma$, we   construct an admissible sequence
$\cI = (i_0, i_1, \ldots, i_{\alpha})$ with $w \in \Dom(h_{\cI})$ so that $(\cI , w)$ covers $\gamma$, and has ``uniform domains''.
Inductively choose a partition of the interval $[0,1]$, say $0 = s_0 < s_1 < \cdots < s_{\alpha} = 1$,  such that for each $0 \leq \ell \leq \alpha$,
$$\gamma([s_{\ell}, s_{\ell + 1}]) \subset D_{\F}(x_{\ell}, \eFU) \quad , \quad x_{\ell} = \gamma(s_{\ell}).$$
As a notational convenience, we have let
$s_{\alpha+1} = s_{\alpha}$, so that $\gamma([s_{\alpha}, s_{\alpha + 1}]) = x_{\alpha}$.
Choose $s_{\ell + 1}$ to be the largest value of $s_{\ell} < s \leq 1$ such that $\dF(\gamma(s_{\ell}), \gamma(t)) \leq \eFU$ for all  $s_{\ell} \leq t \leq s$,  then    $\alpha \leq   \| \gamma \|/\eFU$. 

For each $0 \leq \ell \leq \alpha$, choose an index $1 \leq i_{\ell} \leq \nu$ so that $ B_{\fM}(x_{\ell}, \eU) \subset U_{i_{\ell}}$.
Note that, for all $s_{\ell} \leq t \leq s_{\ell +1}$, $B_{\fM}(\gamma(t), \eU/2) \subset U_{i_{\ell}}$, so that
$x_{\ell+1} \in U_{i_{\ell}} \cap U_{i_{\ell +1}}$. It follows that $\cI_{\gamma} = (i_0, i_1, \ldots, i_{\alpha})$ is an admissible sequence.
Set $h_{\gamma} = h_{\cI_{\gamma}}$ and note that  $h_{\gamma}(w) = w'$.

The construction of the admissible sequence $\cI_{\gamma}$ above has the important   property,  that 
  $h_{\cI_{\gamma}}$ is the composition of generators of $\cGF^*$ which each have a uniform lower bound estimate $ \eTU$ on the radii of the metric balls centered at the orbit and which are contained in their domains.
To see this, let  $0 \leq \ell < \alpha$, and note that $x_{\ell+1} \in D_{\F}(x_{\ell +1}, \eFU)$ implies that for some $s_{\ell} < s_{\ell + 1}' < s_{\ell + 1}$, we have that
$\gamma([s_{\ell + 1}', s_{\ell + 1}]) \subset D_{\F}(x_{\ell +1}, \eFU)$. Hence,
\begin{equation} \label{eq-unifest}
B_{\fM}(\gamma(t), \eU/2) \subset U_{i_{\ell}} \cap U_{i_{\ell +1}} ~ , ~ {\rm for ~ all} ~ s_{\ell + 1}' \leq t \leq s_{\ell + 1} ~.
\end{equation}
Then for all $s_{\ell + 1}' \leq t \leq s_{\ell + 1}$, the uniform estimate defining $\eTU > 0$ in \eqref{eq-transdiam} implies that
\begin{equation} \label{eq-domains}
B_{\fX}(\pi_{i_{\ell}}(\gamma(t)), \eTU ) \subset \fT_{i_{\ell} , i_{\ell +1}} ~ {\rm and} ~
B_{\fX}(\pi_{i_{\ell +1}}(\gamma(t)), \eTU ) \subset \fT_{i_{\ell +1} , i_{\ell}} ~ .
\end{equation}
For the admissible sequence $\cI_{\gamma} = (i_0, i_1, \ldots, i_{\alpha})$,
recall that $x_{\ell} = \gamma(s_{\ell})$ and $w_{\ell} = \pi_{i_{\ell}}(x_{\ell})$.
By   definition \eqref{eq-defholo} of $\ds h_{\cI_{\gamma}}$, the condition \eqref{eq-domains} implies that
$D_{\fX}(w_{\ell} , \eTU) \subset \Dom(h_{\ell})$ as was claimed. 

\begin{defn}\label{def-goodPC}
Let $\gamma \colon  [0,1] \to \fM$ be a   path starting at $x = \tau(w)$ and ending at $y = \tau(w')$. Then a \emph{good plaque-chain covering} of $\gamma$ is the plaque-chain, starting at $x$,    associated to an admissible sequence $\cI_{\gamma} = (i_0, i_1, \ldots, i_{\alpha})$   as constructed above with $\alpha \leq   \| \gamma \|/\eFU$, and which satisfies \eqref{eq-domains}.
\end{defn}

\subsection{Equivalence of leafwise holonomy} \label{subsec-homotopyindepence}

We give several criteria for when two holonomy maps must agree. Consider first the case where we are given a path $\gamma$ with $x = \gamma(0)$ and $y = \gamma(1)$. Let  $\cI = (i_0, i_1, \ldots, i_{\alpha})$ and
$\cJ = (j_0, j_1, \ldots, j_{\beta})$ be admissible sequences such that both $(\cI, \xi_0)$ and $(\cJ, \xi_0')$ cover the   path $\gamma$. 
Then  
$$x \in  \cP_{i_0}(\xi_0) \cap \,   \cP_{j_0}(\xi'_0) ~ , \quad y \in \cP_{i_{\alpha}}(\xi_{\alpha})  \cap \,  \cP_{j_{\beta}}(\xi'_{\beta}) .$$ 

Thus both $(i_0 , j_0)$ and $(i_{\alpha} , j_{\beta})$ are admissible, and
$\xi'_0 = h_{j_{0} , i_{0}}(\xi_0)$, $\xi_{\alpha} = h_{i_{\alpha} , j_{\beta}}(\xi'_{\beta})$.
The proof of the following is intuitively clear, with details  in \cite{ClarkHurder2013}.
\begin{lemma}\label{lem-copc}
The point $\xi_0$ is contained in the domains of both   $h_{\cI}$ and
$\ds h_{i_{\alpha} , j_{\beta}} \circ h_{\cJ} \circ h_{j_{0} , i_{0}}$, and the intersection of their domains is a non-empty clopen subset on which the maps agree. \hfill $\Box$ 
\end{lemma}

Next, consider paths $\gamma, \gamma' \colon [0,1] \to \fM$   with $x = \gamma(0) = \gamma'(0)$ and $y = \gamma(1) = \gamma'(1)$. Suppose that $\gamma$ and $\gamma'$ are homotopic relative endpoints. That is, assume  there exists a continuous map $H \colon [0,1] \times [0,1] \to \fM$ with 
$$H(0,t) = \gamma(t) ~, ~ H(1,t) = \gamma'(t) ~ , ~ H(s,0) = x ~ {\rm and} ~ H(s,1) = y \quad {\rm for ~ all} ~ 0 \leq s \leq 1$$
Then there exists partitions $0 = s_0 < s_1 < \cdots < s_{\beta} = 1$ and $0 = t_0 < t_1 < \cdots < t_{\alpha} = 1$ such that for each pair of indices $0 \leq j < \beta$ and $0 \leq k < \alpha$,   there is an index $1 \leq i(j,k)\leq \nu$ such that 
$$H([s_j,s_{j+1}] \times [t_k, t_{k+1}] ) \subset D_{\F}(H(s_j, t_k), \eFU) \subset U_{i(j,k)}$$
Then proceeding using the methods above and a standard induction argument, we obtain:
\begin{lemma}\label{lem-homotopy}
Let $\gamma, \gamma' \colon [0,1] \to \fM$ be   paths with $x = \gamma(0) = \gamma'(0)$ and $y = \gamma(1) = \gamma'(1)$, and suppose they are homotopic relative endpoints. Then the induced holonomy maps $h_{\gamma}$ and $h_{\gamma'}$ agree on an open neighborhood of $\xi_0 = \pi_{i_0}(x)$. 
\end{lemma}
The   strongly convex property of the plaques and the techniques above also yield   the following result.
 \begin{lemma} \label{lem-domainconst}
 Let  $\gamma, \gamma' \colon [0,1] \to \fM$ be   paths. Suppose that $x = \gamma(0), x' = \gamma'(0) \in U_i$ and
  $y = \gamma(1), y' = \gamma'(1) \in U_j$. If $d_{\fM}(\gamma(t) , \gamma'(t)) \leq \eU/4$ for all $0 \leq t \leq 1$, then the induced holonomy maps $h_{\gamma}, h_{\gamma'}$ agree on their common domain $\Dom(h_{\gamma}) \cap \Dom(h_{\gamma'}) \subset \fT_i$.
\end{lemma}

\subsection{Plaque-chains to leafwise paths}\label{subsec-domains}

 Recall another basic construction, which associates  to an admissible sequence $\cI = (i_0, i_1, \ldots, i_{\alpha})$ and $\xi_0 \in \Dom(h_{\cI_0})$  a leafwise path.
For each $0 < \ell \leq \alpha$,  choose  $z_{\ell} \in \cP_{\ell -1}(\xi_{\ell -1}) \cap \, \cP_{\ell}(\xi_{\ell})$.
Let
$\gamma_{\ell} \colon [(\ell -1)/\alpha , \ell / \alpha] \to L_{x_0}$ be the leafwise piecewise geodesic segment from $x_{\ell -1}$ to $z_{\ell}$ to $x_{\ell}$. Define the   path $\gamma_{\cI} \colon [0,1] \to L_{x_0}$ from $x_0$ to $x_{\alpha}$ to be the concatenation of these paths.
Note that $\| \gamma_{\cI} \| \leq 2 \alpha \dFU$ by Proposition~\ref{prop-regular}.\ref{item-uniform}.

Let $(\cI', \xi_0)$  be a plaque-chain covering of      $\gamma_{\cI}$, then for $U \subset \Dom(h_{\cI}) \cap \Dom(h_{\gamma_{\cI}})$ we have   $h_{\cI} |u  = h_{\gamma_{\cI}} | U$ by Lemma~\ref{lem-copc}. Of course, we may take $\cI' = \cI$ and then the conclusion is tautologous, or may take $\cI'$ to define 
a good plaque-chain covering as in Definition~\ref{def-goodPC}.

 The study of the dynamics of the pseudogroup $\cGF$ acting on $\fX$ differs from the case of a group $\G$ acting on $\fX$, in that for a group action each $\gamma \in \G$ defines a homeomorphism $h_{\gamma} \colon \fX \to \fX$.  For a pseudogroup action,  given $g  \in \cGF$ and $w \in \Dom(g)$, 
 there is some clopen neighborhood $w \in U \subset \Dom(g)$ for which $g | U = h_{\cI} | U$ where $\cI$ is   admissible sequence  with $w \in \Dom(h_{\cI})$. 
 By the definition of a pseudogroup, every $g \in \cGF$ is the ``union'' of such maps in $\cGF^*$. 
 For our applications to the structure of the dynamics of $\F$, it is very useful to estimate the maximal domain for the holonomy map defined by an admissible sequence.     
 Recall that  $\eTU$ was defined by \eqref{eq-transdiam}.

\begin{prop}\label{prop-domest} 
For each $0 < \epsilon \leq \eTU$ and integer $\alpha > 0$, there exists
$0 < \delta(\epsilon, \alpha) \leq \epsilon$
so that for each admissible sequence $\cI$ with length at most $\alpha$, $w \in \Dom(h_{\cI})$ and $w' = h_{\cI}(w)$,  there exists an admissible sequence $\cI'$ with $D_{\fX}(w_0 , \delta(\epsilon, \alpha)) \subset \Dom(h_{\cI'})$ such that $h_{\cI'}$ coincides with $h_\cI$ on the intersection of their domains, and  
\begin{equation}\label{eq-maxdomains}
h_{\cI'}(D_{\fX}(w_0, \delta(\epsilon, \alpha))) \subset D_{\fX}(w', \e)
\end{equation}
\end{prop}
\proof
Let $\gamma_{\cI}$ be the  piecewise geodesic from $x = \tau(w) \in \cT$ to $y=\tau(w')$ as constructed   above.
Then by the method of Section~\ref{subsec-pathstochains}, there is an admissible sequence $\cI'$ which defines a good plaque-chain covering $(\cI', w)$ of $\gamma_{\cI}$ satisfying the   condition \eqref{eq-domains} on domains. The result then follows using induction, as in the proof of \cite[Proposition~5.7]{CHL2013a}.
\endproof

\subsection{Dynamics of matchbox manifolds}\label{sec-mme}

We   recall two     definitions from topological dynamics, that of \emph{equicontinuous} and \emph{expansive} dynamics, as adapted to actions of pseudogroups. 
 
\begin{defn} \label{def-expansive}
Say that the action of the    pseudogroup $\cGF$ on $\fX$ is \emph{expansive}, or more properly that it is $\e$-expansive, if there exists $\e > 0$ such that for all $w, w' \in \fX$, there exists $g \in \cGF^*$ with $w, w' \in D(g)$ such that $d_{\fX}(g(w), g(w')) \geq \e$.
\end{defn}

\begin{defn} \label{def-equicontinuous}
Say that the action of the    pseudogroup $\cGF$ on $\fX$ is  \emph{equicontinuous} if for all $\epsilon > 0$, there exists $\delta > 0$ such that for all $g \in \cGF^*$, if $w, w' \in D(g)$ and $d_{\fX}(w,w') < \delta$, then $d_{\fX}(g(w), g(w')) < \epsilon$.
Thus, $\cGF^*$ is equicontinuous as a family of local group actions. 
\end{defn}

Equicontinuity for $\cGF$ gives \emph{uniform} control over the domains of arbitrary compositions of generators, so that  a much stronger conclusion than  Proposition~\ref{prop-domest} is true, as shown in \cite{ClarkHurder2013}:
\begin{prop} \label{prop-uniformdom}
Assume the holonomy pseudogroup $\cGF$ of $\F$ is equicontinuous. Then there exists $\dTU > 0$ such that for every leafwise path
$\gamma \colon [0,1] \to \fM$, there is a corresponding admissible sequence $\cI_{\gamma} = (i_0 , i_1 , \ldots , i_{\alpha})$ so that
$B_{\fX}(w_0, \dTU) \subset D(h_{\cI_{\gamma}})$, where $x = \gamma(0)$ and $w_0 = \pi_{i_0}(x)$.
Moreover, for all $0 < \e \leq \eTU$ there exists $0 < \delta  \leq \dTU$ independent of the path $\gamma$, such that
$h_{\cI_{\gamma}}(D_{\fX}(w_0, \delta)) \subset D_{\fX}(w', \e)$ where $w' = \pi_{i_{\alpha}}(\gamma(1))$.
\end{prop}

It is possible that these are the only two possibilities for a minimal action of a pseudogroup on a Cantor space. A proof of the following would generalize   results in    \cite{AY1980,AGW2007,DM2008}: 
\begin{conj}\label{conj-dichotomy}
Let $\cGF$ act minimally on a Cantor space $\fX$ as above. Then either the action is equicontinuous, or it is expansive for some $\e > 0$.
\end{conj}

  \subsection{Holonomy groupoids}
 We next consider the groupoid \cite{Haefliger1984} formed by germs of maps in $\cGF$.  
  
Let  $U, U', V, V' \subset \fX$ be open subsets with  $w \in U \cap U'$. Given homeomorphisms    $h \colon U \to V$  and $h' \colon U' \to V'$    with $h(w) = h'(w)$,   then   $h$ and $h'$ have the same \emph{germ at $w$}, and write    $h \sim_w h'$,   if there exists an open neighborhood $w \in W \subset U \cap U'$ such that $h | W= h' |W$. Note that $\sim_w$ defines an equivalence relation. 

\begin{defn}\label{def-germ}
The \emph{germ of $h$ at $w$} is the equivalence class $[h]_w$ under the relation ~$\sim_w$. The  map  $h \colon U \to V$ is called a \emph{representative} of  $[h]_w$.
The point $w$ is called the source of  $[h]_w$ and denoted $s([h]_w)$, while $w' = h(w)$ is called the range of  $[h]_w$ and denoted $r([h]_w)$.
\end{defn}

  The \emph{holonomy   groupoid} $\GF$ is the the collection of all   germs $[h]_w$ for   $h \in \cGF$ and $w \in \Dom(h)$,  equipped with  the sheaf topology  for maps over $\fX$. Let $\cRF \subset \fX \times \fX$ denote the equivalence relation on  $\fX$ induced by $\F$, where     $(w,w') \in \cRF$ if and only if $w,w'$ correspond to points on the same leaf of $\F$. The product map $s \times r \colon \GF \to \cRF$ is  \'etale; that is,  a local homeomorphism with discrete fibers.
  
 We introduce a convenient notation for elements of $\GF$.
 Let $(w,w') \in \cRF$, and let $\gamma$  denote a path from $x = \tau(w)$ to $y = \tau(w')$.
Proposition~\ref{prop-domest} implies that we can choose an admissible sequence $\cI_{\gamma}$ for which the domain of the induced map $h_{\cI}$ satisfies condition \eqref{eq-maxdomains}. 
By  Lemma~\ref{lem-homotopy}, the germ  $[h_{\cI_{\gamma}}]_w$   depends only on the endpoint-fixed homotopy class of $\gamma$, so that we may assume $\gamma$ is a geodesic between $x$ and $y$. Moreover, by Lemma~\ref{lem-domainconst} the germ $[h_{\cI_{\gamma}}]_w$ is independent of the admissible sequence covering $\gamma$. Thus, we introduce the notation $\gamma_w = [h_{\cI_{\gamma}}]_w$ for this germ, which depends only on the end-point fixed homotopy class of $\gamma$. Given $[h]_w \in \GF$ by Section~\ref{subsec-plaquechains} the germ is induced by a plaque-chain covering a piecewise geodesic between $x$ and $y$, so all elements of $\GF$ are represented by $\gamma_w$ for some path $\gamma$ from $x$ to $y$.

 These remarks imply there is a well-defined surjective homomorphism, the \emph{holonomy map},  
 \begin{equation}\label{eq-holodef}
h_{\F,x} \colon \pi_1(L_x , x) \to \G_w^w \equiv \left\{  [g]_w \in \GF  \mid    r([g]_w) =w \right\}  . 
\end{equation}
Moreover,  if $y,z \in L$ then the homomorphism
$h_{\F , y}$ is conjugate (by an element of $\cGF$) to the homomorphism $h_{\F , z}$.
A leaf $L$ is said to have \emph{non-trivial germinal holonomy} if for some $y \in L$, the homomorphism $h_{\F , y}$ is non-trivial. If the homomorphism $h_{\F , y}$ is trivial, then we say that $L_y$ is a \emph{leaf without holonomy}. This property depends only on $L$, and not the choice of  $y \in L$.

A matchbox manifold $\fM$ is said to be \emph{without holonomy}, if every leaf is without holonomy.

As an application of these  remarks, we obtain:
\begin{lemma}\label{lem-homotopymin}
Given a path $\gamma \colon [0,1] \to \fM$  with $x = \gamma(0)$ and $y = \gamma(1)$. Suppose that   $L_x$ is a leaf without holonomy. Then there exists a leafwise geodesic segment $\gamma'  \colon [0,1] \to \fM$  with $x = \gamma'(0)$ and $y = \gamma'(1)$,  such that $\|\gamma' \| = \dF(x,y)$,    and $h_{\gamma}$ and $h_{\gamma'}$ agree on an open neighborhood of $\xi_0$.
\end{lemma}
\proof
The leaf $L_x$ containing $x$ is a complete Riemannian manifold, so there exists  a geodesic segment $\gamma'$ which is length minimizing between $x$ and $y$.
Then the holonomy maps  $h_{\gamma}$ and $h_{\gamma'}$ agree on an open neighborhood of $\xi_0  = \pi_{i_0}(x)$ by the definition of germinal holonomy. 
\endproof

\subsection{Word and path length} \label{subsec-wordlength}
We recall the word length function on  $\cGF^*$  and $\GF$,   and the path length function on $\GF$, and give estimates comparing these notions of length.

For $\alpha \geq 1$, let  $\cGF^{(\alpha)}$ be the collection of holonomy homeomorphisms $h_{\cI} | U \in \cGF^*$ determined   by admissible paths $\cI = (i_0,\ldots,i_k)$ such that $k \leq \alpha$ and $U \subset \Dom(h_{\cI})$. For each $\alpha$, let $C(\alpha)$ denote the number of admissible sequences of length at most $\alpha$. As there are at most $\nu^2$ admissible pairs $(i,j)$, we have the basic estimate that $C(\alpha) \leq \nu^{2 \alpha}$. This upper bound estimate grows exponentially with $\alpha$, though the exact growth rate of $C(\alpha)$ may be much less.

For each $g \in \cGF^*$ there is some $\alpha$ such that $g \in \cGF^{(\alpha)}$. Let $\|g\|$ denote the least such $\alpha$, which is called the \emph{word length} of $g$.  Note that  $\cGF^{(1)}$ generates $\cGF^*$.

We use the word length on $\cGF^*$ to  define the word length on  $\GF$, where for  $\gamma_w \in \GF$, set
\begin{equation}
\| \gamma_w \| ~ = ~ \min ~ \left\{ \| g \| \mid [g]_w = \gamma_w ~ {\rm for}~ g \in \cGF^*  \right\} .
\end{equation}

  Introduce the \emph{path length} of $\gamma_w \in \GF$, by considering the infimum of the lengths $\| \gamma'\|$ for all   piecewise smooth curves $\gamma'$  for which $\gamma_w' = \gamma_w$. That is, 
\begin{equation}\label{eq-groupoidpathlength}
\ell(\gamma_w)  ~ = ~ \inf ~ \left\{ \| \gamma' \| \mid \gamma'_w = \gamma_w   \right\} .
\end{equation}
Note that if $L_w$ is a leaf without holonomy, set $x = \tau(w)$ and $y = \tau(w')$, then Lemma~\ref{lem-homotopymin} implies that 
$\ell(\gamma_w) = \dF(x,y)$. 
 This yields a   fundamental estimate:
\begin{lemma}\label{lem-comparisons}
Let  $[g]_w \in \GF$ where $w$ corresponds to a leaf without holonomy. Then
\begin{equation}\label{eq-comparisons}
\dF(x,y)/2\dFU ~ \leq ~ \| [g]_w \| ~ \leq ~  1 + \dF(x,y)/\eFU
\end{equation}
\end{lemma}
\proof
Let $[g]_w$ be represented by $h_{\cI}$ where $\cI = (i_0, i_1, \ldots, i_{\alpha})$   of length $\alpha = \| [g]_w \|$. Then $\cI$ defines a plaque-chain $\cP_{\cI}$ as in \eqref{eq-plaquechain} where each plaque is contained in a disk  $D_{\F}(x_{w_i} , \dFU)$, so the piecewise-geodesic $\gamma_{\cI}$ which $\cP_{\cI}$ defines as in Section~\ref{subsec-plaquechains} has length at most  $2  \alpha \, \dFU$. Thus, $\dF(x,y) \leq 2  \alpha \, \dFU$.

Conversely, as noted in Section~\ref{subsec-pathstochains}, given a path $\gamma$ from $x$ to $y$ there is a good plaque-chain $\cI$ which covers $\gamma$, with 
$\alpha \leq 1 + \| \gamma \|/\eFU$ which yields the right-hand-side estimate in \eqref{eq-comparisons}.
\endproof

 In general, the leafwise   distance function $\dF$ of distinct leaves cannot be compared for points which have a large separation, except when it is possible to define a ``shadowing'' of a path in one leaf by a path in a nearby leaf, as described in the following.
 
  Let $x \in \fM$ with $L_x$ the leaf containing it, and assume that $L_x$ is without holonomy. Given $y \in L_x$ let $\gamma_{x,y}$ be a geodesic from $x$ to $y$ with $\| \gamma_{x,y}\| = \dF(x,y)$.
    Choose a plaque-chain $\cP_{\cI}$ covering $\gamma_{x,y}$ as in Section~\ref{subsec-pathstochains} so that  Proposition~\ref{prop-domest}  holds for the holonomy map $h_{\cI}$ it induces. Let  $U_{i_0}$ be the first chart of the plaque-chain, so that   $x \in \cP_{i_0}(\xi_0)$, with notation as in \eqref{eq-cover}.
 Let $x' \in U_{i_0}$ have the same leaf coordinate as $x$, so $\lambda_{i_0}(x') = \lambda_{i_0}(x)$. 
 Assume that  $w' = \pi(x') \in \Dom(h_{\cI})$, then we   
   define $y' \in U_{i_{\alpha}}$ corresponding to $y$ in the same way, for the last plaque $\cP_{i_{\alpha}}(\xi_{\alpha})$ in the plaque-chain. 
       The  curve $\gamma_{x,y}$ is given in \eqref{eq-cover}  as a concatenation of piecewise-geodesic segments, each of which is contained in a plaque of the plaque-chain $\cP_{\cI}$.  Then by shadowing these geodesic segments within each respective foliation chart, we obtain a piecewise geodesic   $\wtgamma_{x',y'}$   from $x'$ to $y'$. 
     Then using the   estimates \eqref{eq-metricbounds} and Proposition~\ref{prop-domest},  we obtain:
\begin{lemma}\label{lem-pathlengths}
 Let $\e > 0$ and $R > 0$, then there exists $\delta(\e,R) > 0$ such that if  $x \in \fM$ is contained in a leaf $ L_x$ without holonomy, and $y \in L_x$ satisfies $\dF(x,y) \leq R$, then for   $x'$ as above with $w' = \pi(x') \in B_{\fX}(w, \delta(\e,R))$,  it follows that $w' \in \Dom(h_{\cI})$ and for the corresponding endpoint $y' \in L_{x'}$ we have
  \begin{equation}\label{eq-pathlengths}
\dF(x',y') ~ \leq ~ (1+\e) \, \dF(x,y)  \quad {\rm and} \quad \dF(x,y) ~ \leq ~ (1+\e) \, \dF(x',y')  ~ .
\end{equation}
 \end{lemma}
\proof
Let $\gamma_{x,y}$ be a leafwise geodesic from $x$ to $y$ with $\| \gamma_{x,y}\| = \dF(x,y)$. Construct a plaque-chain covering $\cP_{\cI}$ of $\gamma_{x,y}$ as above.  Let  $\delta_{\cU}(\e)$ be the constant appearing in \eqref{eq-metricbounds}, then by Lemma~\ref{lem-comparisons} and Proposition~\ref{prop-domest}, there exists $\delta(\delta_{\cU}(\e),R) > 0$ so that $ B_{\fX}(w, \delta(\delta_{\cU}(\e),R)) \subset \Dom(h_{\cI})$. Moreover,  for any $w' \in B_{\fX}(w, \delta(\delta_{\cU}(\e),R))$ the plaques defined by the shadowed points along the piecewise geodesic $\wtgamma_{x',y'}$ from $x'$ to $y'$   have separation at most  $\delta_{\cU}(\e)$. Then use
the estimates in \eqref{eq-metricbounds} in each plaque to deduce the first estimate in \eqref{eq-pathlengths}. The second estimate follows by considering the path $\gamma_{x,y}$ as shadowing the path $\wtgamma_{x',y'}$. 
\endproof

The technical result Lemma~\ref{lem-pathlengths} illustrates some of the difficulties working with the dynamics of pseudogroups, as   the notion of distance along paths in leaves requires care with the initial point.

\subsection{Induced groupoids} \label{subsec-inducedgpds}

For an open subset $W \subset \fX$,  consider the following subsets of $\GF$, 
\begin{eqnarray}
\G_W & = & \left\{ [g]_w \in \GF \mid w \in W \right\}  \nonumber\\
\G_W^W & = & \left\{ [g]_w \in \GF \mid w \in W ~, ~ r([g]_w) \in W\right\}  \nonumber\\
\G_w^W & = & \left\{  [g]_w \in \GF  \mid    r([g]_w) \in W \right\}  \nonumber
\end{eqnarray}
so we have the inclusions
$$ \G_w^w \subset \G_w^W \subset \G_W^W \subset \G_W \subset \GF ~ .$$

The set   $\G_W^W$ is called the induced  groupoid on $W$. 
On the other hand,  the sets $\G_W$ and $\G_w^W$ are not groupoids, as they need not contain the inverses of their elements.

 Given $w \in W$ and $R > 0$, define   the word and   path length filtration of the groupoid $\G_w^W$, 
\begin{equation}\label{eq-groupoidfiltraqtion}
\G_{w}^{W,R}   =   \left\{ \gamma_{w} \in \G_w^W  \mid  \ell(\gamma_w)   \leq R \right\} \quad , \quad \G_{w}^{W,\alpha}   =   \left\{ [g]_{w} \in \G_w^W  \mid  g \in \cGF^{(\alpha)} ~{\rm with }~ w \in \Dom(g) \right\}  .
\end{equation}
By Lemma~\ref{lem-comparisons}, we then have the estimates:
 \begin{lemma} Let $L_w$ be a leaf without holonomy.  For $\alpha > 0$, set $R^{(\alpha)} = 2\alpha \dFU$ and $R_{(\alpha)} = \alpha \eFU$ where we recall from Section~\ref{subsec-metricregular} that $\eFU < \dFU$ so that $R^{(\alpha)} > R_{(\alpha)}$, then 
\begin{equation}\label{eq-inclusions}
\G_{w}^{W, R_{(\alpha)}}   \subset \G_{w}^{W, \alpha}   \subset \G_{w}^{W, R^{(\alpha)}}  .
\end{equation}
\end{lemma}
 The inclusions in \eqref{eq-inclusions} are used in the definition of the coding of the orbits of $\cGF^*$ in Section~\ref{sec-partialcoding}.
\bigskip
 
\section{Delone sets and Voronoi tessellations} \label{sec-voronoi}

In this section, we consider   the    \emph{Delone sets}  and their associated \emph{Voronoi tessellations} on the leaves of $\F$ which are associated to a choice of a clopen transversal set $W \subset \fX$.  These ideas are discussed in extensive detail in \cite{CHL2013a}. It is important to note that,  as no geometric assumption is made on the leaves    of $\F$,   many standard techniques for Voronoi tessellations in the literature for Euclidean spaces  do not apply to this general case.     

Assume that  $\fM$ is a minimal matchbox manifold, so    for every $w \in \fX$,  its $\cGF$-orbit $\cO(w)$ is dense. 
This   implies   that the restriction of $\cGF$ to any clopen subset $U \subset \fX$ yields an induced dynamical system which is   ``equivalent'' to the action of $\cGF$ on $\fX$,  as discussed in \cite{CHL2013c}. This restriction property is fundamental    for the study  of minimal Cantor actions. 
Here is an essential application  of minimality: 
 \begin{lemma}\label{lem-finitegen}
Let $W \subset \fX$ be an open subset. Then there exists an integer  $\alpha_W$ such that $\fX$ is covered by the collection 
$\{h_{\cI} (W)  \mid   h_{\cI} \in \cGF^{(\alpha_W)} \}$. 
Moreover,    as   $\diamX(W) \to 0$, we have $\alpha_W \to \infty$. 

\end{lemma} 
\proof 
Consider the collection of open sets  $ \cQ = \{h_\cI (W) ~|~   h_{\cI}  \in \cGF~\}$ where we adopt the abuse of notation that 
$h_{\cI}(W) = h_{\cI}(W \cap \Dom(h_{\cI}))$. 
As the action of $\cGF$ on $\fX$ is minimal, for all $y \in \fX$ there exists $h_{\cI}$ with $h_{\cI}(y) \in W$. Thus, the collection $\cQ$ is an open covering of $\fX$. 
Since $\fX$ is compact, there exists a finite subcover $\{h_{\cI_1}(W), \ldots, h_{\cI_m}(W)\}$ of $\cQ$, where $\cI_i = (j_0,\ldots,j_{\alpha_i})$.  Let 
\begin{equation}\label{eq-alphaW}
\alpha_W = \max \{\alpha_i ~| ~ 1 \leq i \leq m\}.
\end{equation}
The conclusion that $\alpha_W \to \infty$  as  $\diamX(W) \to 0$ follows immediately, as assuming $\alpha_W$ is bounded implies that $\fX$ admits   coverings   by clopen sets with arbitrarily small diameters, yet bounded in number, which is impossible.  
\endproof

 \subsection{Delone sets} \label{subsec-delone}

We define  nets in the leaves of $\F$ related to the dynamics of $\cGF$.

\begin{defn} \label{def-net}
Let $(X, d_X)$ be a complete separable metric space. Given
$0 < \lambda_1 < \lambda_2$, a subset $\cN \subset X$ is a \emph{$(\lambda_1 , \lambda_2)$-net} (or \emph{Delone set}) if:
\begin{enumerate}
\item $\cN$ is $\lambda_1$-separated: for all $y \ne z \in \cN$, $d_{X}(y,z) \geq \lambda_1$;
\item $\cN$ is $\lambda_2$-dense: for all $x \in X$, there exists some $z \in \cN$ such that $d_{X}(x,z) \leq \lambda_2$.
\end{enumerate}
\end{defn}

Given    $w \in \fX$, let    $L_w$ denote the leaf containing $x = \tau(w)$.
Define $\cN_w = L_w \cap   \cT$.  Then \eqref{eq-dT} implies    that $\cN_w$ is an $\dT$-separated subset of $L_w$. On the other hand, $\cU$ is an open covering of $\fM$, and therefore by Proposition~\ref{prop-regular}.\ref{item-uniform}  each point of $L_w$ has distance at most   $\dFU$ from a point of $\cN_w$. 
That is, $\cN_w$ is always an $(\dT , \dFU)$-net in $L_w$.

 In the following, we    assume without loss of generality that the open set $W$ satisfies $W \subset \fT_1$, and set $\cT_W = \tau_1(W) \subset \fM$. For each  $w \in W$,   define $\cN_w^W = L_w \cap \cT_W \subset \cN_w$.   
 
 If $w $ is a fixed point for some non-trivial holonomy map $h_{\gamma}$ defined by a path $\gamma$, then the length of $\gamma$ gives an approximate upper bound on the constant $\lambda_1$ for the net $\cN_w^W$.  Define the function
\begin{equation}\label{eq-separationW}
\lambda_1(\delta) = \inf \left\{ \lambda \mid w \in W \subset \fX  , ~ \diamX (W) \leq \delta   ,  ~ \cN_w^W~ {\rm  is} ~  \lambda{\rm -separated}  \right\}.
\end{equation}
We always have $\lambda_1(\delta) \geq \dT$ as $\cN_w^W \subset \cN_w$ for all $w \in \fX$, and $\cN_w$ is   $\dT$-separated.  The following result   shows there is  no upper bound on the   function $\lambda_1$ for a foliation without holonomy.

\begin{lemma}\label{lem-separationW} 
Let  $\fM$ be a matchbox manifold  without holonomy, then 
 $\lambda_1(\delta)$ is increasing, and unbounded as $\delta \to 0$.
\end{lemma}
\proof
If $0 < \delta' < \delta$ then $\lambda_1(\delta') \geq \lambda_1(\delta)$ follows from the    definition, so the function $\lambda_1$ is increasing. 

We show that the function $\lambda_1$ is unbounded.  
  Assume there exists $\lambda_* > 0$, such that 
for each integer  $\ell > 0$, there exists  distinct points  $u_{\ell}, v_{\ell} \in \fX$ such that  $\dX(u_{\ell}, v_{\ell}) \leq 1/\ell$ and there is a path $\gamma_{\ell}$ joining $y_{\ell} = \tau(u_{\ell})$ and $z_{\ell} = \tau(v_{\ell})$   with $\|\gamma_{\ell} \| \leq \lambda_*$. 
As $\fX$ is compact, by passing to subsequences, we can assume that there exists $w \in \fX$ which is the limit of both   sequences $\{u_{\ell}\}$ and $\{v_{\ell}\}$.  

By Lemma~\ref{lem-comparisons}, we can choose a plaque-chain   $\cI_{\ell}$ which covers $\gamma_{\ell}$  with length at most   $ 1 + \lambda_*/\eFU$. Since there are only a finite number of admissible plaque-chains of   length at most   $ 1 + \lambda_*/\eFU$, we can pass to a subsequence, and assume without loss of generality that $\cI = \cI_{\ell}$ is independent of $\ell$. Then Proposition~\ref{prop-domest} implies that we can assume that there is some $\delta' >0$ and $\ell_0 > 0 $ with 
  $u_{\ell}, v_{\ell} \in B_{\fX}(w, \delta') \subset \Dom(h_{\cI})$ for all $\ell \geq \ell_0$. 
Then  the   holonomy map $h_{\cI}$ satisfies $h_{\cI}(u_{\ell}) = v_{\ell}'$ for all $\ell \geq \ell_0$, and thus we have $h_{\cI}(w) = w$. 
The assumption $u_{\ell} \ne v_{\ell}$ for all $\ell \geq \ell_0$ implies that   $h_{\cI}$ is not the identity map on any open neighborhood of $w$.
 Thus,  the germ $[h_{\cI}]_{w}$ is non-trivial, so that the leaf  $L_w$ has is   germinal holonomy, which contradicts our assumption that $\F$ is without holonomy.
\endproof

 \begin{remark}
 {\rm 
   The function $\lambda_1(\delta)$      is a fundamental property of the  dynamics of the pseudogroup $\cGF$ acting on $\fX$. From the above proof, it is clear that the function  is closely related to the ``return times'' of the action near to $w$, though little more seems to be known  about its properties in general. Note that   Forrest in \cite{Forrest2000} derives estimates for the constants $\{ \lambda_1, \lambda_2\}$ for the net $\cN_w^W$ in the case where $\cGF$ is induced from a minimal action of $\mZ^n$ on a Cantor set. It is a very interesting problem to derive estimates for   these constants  in the more general cases of group actions and pseudogroups. }
 \end{remark}

We next show that the density   of the net $\cN_w^W$ in the leaf $L_w$ has a uniform estimate in terms of the integer $\alpha_W$ introduced in Lemma~\ref{lem-finitegen}. 

 \begin{lemma}\label{lem-deloneW}
 The set  $\cN_w^W$ is $(2 \alpha_W \dFU)$-dense in the set $\cN_w$.
 \end{lemma} 
\proof 
Let $L_w$ be the leaf determined by $w \in W$. 
Let  $y \in \cN_w$ and $w' = \pi(y) \in \fX$.
The collection $\{h_{\cJ} (W)  \mid ~ h_{\cJ} \in \cGF^{(\alpha_W)} \}$ is an open covering of $\fX$ by  Lemma~\ref{lem-finitegen}, so there exists $\cJ$ with length at most $\alpha_W$ with $w' \in h_{\cJ}(W)$, which defines a plaque-chain from $x$ to $y$ with at most $\alpha_W$ plaque intersections.
Each  plaque is contained in a disk  $D_{\F}(x_{w_i} , \dFU)$,  thus there is a piecewise geodesic path in $L_w$ from $x$ to $y$ with length at most 
$2\alpha_W  \,  \dFU$.  
\endproof

  \begin{cor}\label{cor-deloneW}
For   $e_W = (2\alpha_W + 1) \dFU$ the set  $\cN_w^W$ is  $e_W$-dense in $L_w$.   That is, for every $x \in \fM$ there is a piecewise geodesic path with length at most $e_W$ joining $x$ to a point of $\cT_W$.
\end{cor} 
\proof
Given    $x \in L_w$  there exists $y \in   \cN_w$ with $\dF(x,y) \leq  \dFU$.
By Lemma~\ref{lem-deloneW},   there is $z \in \cN_w^W$ with $\dF(y,z) \leq 2 \alpha_W \dFU$ and thus $\dF(x, z) \leq (2\alpha_W + 1) \dFU$, 
as was to be shown.
\endproof

\subsection{Voronoi tessellations} \label{subsec-VT}

We   develop   some of the basic concepts of \emph{Voronoi tessellations}  for complete Riemannian manifolds, and in particular for the leaves of $\F$.    As the leaves  of $\fM$  are not necessarily Euclidean,   some of the usual results for tessellations  of $\mR^n$ with the Euclidean metric need not be true in this context (see \cite{LeibonLetscher2000} and  \cite[Introduction]{OBKC2000}).


 The  \emph{Voronoi tessellation} associated to $\cN_w^W$  is a partition of $L_w$ into compact   regions with disjoint interiors, called the ``cells'', where each cell is ``centered'' at a unique point of $\cN_w^W$.  We describe the construction of these cells and and some of their properties. 
 
Introduce the ``leafwise nearest--neighbor'' distance function, where for $y \in  L_w$,
\begin{equation}\label{eq-kF2}
\kwW(y) = \inf \left\{ \dF(y,z) \mid z \in \cN_w^W \right\}.
\end{equation}

Note that $\kwW(y) < \infty$ if and only if $y \in L_w$, and $\kwW(y) = 0$ if and only if $y \in \cN_w^W$.

\begin{defn}
For $y \in \cN_w^W$, define its   \emph{Voronoi cell}   in $L_w$ by
\begin{equation}\label{eq-voronoi2}
\cC_w^W(y) = \left\{ z \in L_w \mid \dF(y,z) = \kwW(z)\right\}.
\end{equation}
\end{defn}
That is, for $x \in \cN_w^W$ the Voronoi cell $\cC_w^W(y)$ consists of the points $z \in L_w$ which are at least as close to $y$ in the leafwise metric as to  any other   point of $\cN_w^W$. 

 Note that if $L_w$ is a convex metric space for the leaf metric, for example when the leaf $L_w$ is a Euclidean space, then  each cell $\cC_w^W(y)$  is a convex subset. However, in the general case,  the   cells $\cC_w^W(y)$ need not be convex, or even simply connected, especially if the leaf $L_w$ has non-trivial topology and the   constants $\lambda_1$ and $\lambda_2$ for the net $\cN_w^W$ are ``large''. 
In this work, we assume only that the leaves of $\fM$ are complete metric spaces, and so develop some of the   properties of the cells $\cC_w^W(y)$ in the decomposition \eqref{eq-Vdecomposition} for this generality.

First, note that by the definition,     for each $y \in L$,  the set $\cC_w^W(y) $ is an intersection of closed subsets of $L_w$, hence is closed. Also, from the definitions  every point of $L_w$ belongs to some $\cC_w^W(z)$ for $z \in \cN_w^W$, so we have:   
\begin{defn} \label{def-voronoitessel}
The \emph{Voronoi tessellation} of $L_w$ associated to  $\cN_w^W$ is the decomposition
\begin{equation}\label{eq-Vdecomposition}
L_w ~ = ~ \bigcup_{y \in \cN_w^W}  ~ \cC_w^W(y)
\end{equation}
\end{defn} 
In general, the geometry of a cell $\cC_w^W(z) \subset L_w$ may be difficult to describe, though we always have the following bounds, which are
 direct consequences of Lemma~\ref{lem-separationW}  and Corollary~\ref{cor-deloneW}.
\begin{lemma}\label{lem-celldiam}
Suppose that $W \subset B_{\fX}(w, \delta)$, then for each $y \in \cN_w^W$, 
\begin{equation}\label{eq-Vcellbounds}
D_{\F}(y, \lambda_1(\delta)/2)   \subset \cC_w^W(y) \subset B_{\F}(y, e_W)  
\end{equation}
 In particular, each cell  $\cC_w^W(y)$ is closed with diameter at most $2e_W$, and thus is compact.  
\end{lemma}

 If $\F$ is a foliation without holonomy, then Lemmas~\ref{lem-separationW}  and \ref{lem-celldiam} imply that  if  $W$ is chosen to have sufficiently small diameter, then each Voronoi cell in $L_w$ defined by the net $\ds \cN_w^W = L_w \cap \tau(W)$ contains a ball of   radius larger than any prescribed size.

Next, we introduce the \emph{star-neighborhoods} of Voronoi cells.  Given $y \in \cN_w^W$, introduce the \emph{vertex-set}
\begin{equation}\label{eq-vertexset}
\cV_w^W(y) = \{ z \in \cN_w^W \mid \cC_w^W(z) \cap \cC_w^W(y) \ne \emptyset \} .
\end{equation}

\begin{defn}\label{lem-star-nbhd}
For $y \in \cN_w^W$ the ``star-neighborhood'' of the Voronoi cell $\cC_w^W(y)$ is the set
\begin{equation}\label{eq-starset}
\cS_w^W(y) = \bigcup_{z \in \cV_w^W(y)} ~ \cC_w^W(z).
\end{equation}
\end{defn}

 \begin{lemma}(\cite{CHL2013a})\label{lem-star2}
 For each $y \in \cN_w^W$ and $z \in \cV_w^W(y)$, we have that $\dF(y,z) \leq 2e_W$.
 Hence, $\cS_w^W(y) \subset B_{\F}(y, 3e_W)$, and the collection $\cV_w^W(y) $ is   finite. 
\end{lemma}

For each $x \in \cN_w^W$ define the boundary of $\cC_w^W(x)$ to be the set  $\ds \partial \cC_w^W(x) \equiv \cC_w^W(x)- {\rm int}(\cC_w^W(x))$.
We  consider the relation between the star-neighborhood of a cell and its boundary $\ds \partial \cC_w^W(x)$.

\begin{lemma}  \label{lem-star3}
 For each $x \in \cN_w^W$, we have that $\ds \cC_w^W(x) \subset {\rm int}(\cS_w^W(x))$.
\end{lemma}
\proof
We always have ${\rm int}(\cC_w^W(x)) \subset {\rm int}(\cS_w^W(x))$, so we must show $\partial \cC_w^W(x) \subset {\rm int}(\cS_w^W(x))$.

Let $y \in \partial \cC_w^W(x)$ with $y \not\in {\rm int}(\cS_w^W(x))$. Then for every $0 < \e < e_W$, the open ball $B_{\F}(y,\e)$ must intercept the complement of $\cS_w^W(x)$ in $L_w$ in a non-empty set; let $y'_{\e}$ be a point in this set.  Thus, there is some $z_{\e} \in \cN_w^W$ for which $y'_{\e} \in \cC_w^W(z_{\e})$ so that 
$\dF(x,z_{\e}) \leq 2e_W + \e < 3 e_W$. There are at most a finite number of such net points $z_{\e}$, so we can choose a subsequence of the points $y'_{1/\ell}$ with $\dF(y, y'_{1/\ell}) < 1/\ell$ and $y'_{1/\ell} \in \cC_w^W(z')$ for all $\ell$, where $z' \in \cN_w^W$ is fixed. As  $\cC_w^W(z')$ is compact, this implies that $\cC_w^W(z') \cap \cC_w^W(x) \ne \emptyset$ and so $y'_{1/\ell} \in \cS_w^W(x)$ contrary to choice.
\endproof

 For $z \in \cN_w^W$ with $z \ne y$, introduce the closed subsets of the leaf $L_w$ defined by 
\begin{eqnarray*}
H_w^W(y,z) & = & \{ \xi \in L_w \mid \dF(y,\xi) \leq \dF(z,\xi)\}\\
L_w^W(y,z)  & = & \{ \xi \in L_w \mid \dF(y,\xi) = \dF(z,\xi)\} 
\end{eqnarray*}

If $L_w$ is isometric to Euclidean space $\mR^n$, then    $H_w^W(y,z)$ is a closed half-space, and   $L_w^W(y,z)$ is the boundary hyperplane.
In   general,   the closed subspace $L_w^W(y,z)$ need not   be a manifold.

The next result is a ``local'' version of  the observation that if $L_w$ is isometric to Euclidean space $\mR^n$, then    each Voronoi cell   $\cC_w^W(y) \subset \mR^n$ is the intersection of half-spaces defined by its boundary planes.

\begin{lemma}\label{lem-cellshyperplanes}
For each $y \in \cN_w^W$ we have \quad 
$\ds \cC_w^W(y) ~ = ~ \cS_w^W(y) ~ \cap \bigcap_{z \in \cV_w^W(y) } ~ H_w^W(y,z)$.  
\end{lemma}
 \proof
 The   identity below and the subsequent inclusion follow from the definitions:
 \begin{equation}\label{eq-halfplanes}
\cC_w^W(y) = \bigcap_{z \in \cN_w^W} ~ H_w^W(y,z) ~ \subset ~ \cS_w^W(y) \cap  \bigcap_{z \in \cV_w^W(y) } ~ H_w^W(y,z) ~ .
\end{equation}

We must show the reverse inclusion in \eqref{eq-halfplanes} holds as well.
Observe  that for $z \ne y$, we have $\cC_w^W(y) \cap \cC_w^W(z) \subset L_w^W(y,z)$ by  the definition.
We claim that  also,   $\cC_w^W(y) \cap L_w^W(y,z) \subset \cC_w^W(z)$. 
Let $\xi \in \cC_w^W(y) \cap L_w^W(y,z)$, so that $\xi$ is closer to $y$ than any other point in $ \cN_w^W$ and $\dF(y,\xi) = \dF(z,\xi)$. Suppose that $\xi \not\in \cC_w^W(z)$, then there exists $z' \in \cN_w^W$ with $\dF(z', \xi) < \dF(z,\xi) = \dF(y,\xi)$, so $z' \not\in \cC_w^W(y)$, which contradicts the choice of $\xi$.
Thus, $\xi \in \cC_w^W(z)$.
 
Now consider $z \in \cN_w^W$ and $z \not\in \cV_w^W(y)$. Then $\cC_w^W(y) \cap L_w^W(y,z) = \emptyset$.
The cell $\cC_w^W(y)$ is compact by Lemma~\ref{lem-celldiam}, so there exists $\epsilon > 0$ such that the $\epsilon$-neighborhood of $\cC_w^W(y)$ is also disjoint from  $L_w^W(y,z)$. Since $\cC_w^W(y) \subset H_w^W(y,z)$,  this $\epsilon$-neighborhood is also contained in $H_w^W(y,z)$.
Note that for all $z' \in \cN_w^W$, either $\cC_w^W(z') \subset H_w^W(y,z)$ or $\cC_w^W(z') \subset H_w^W(z,y)$ so by Lemma~\ref{lem-star3}, we conclude that  $\cC_w^W(z') \subset H_w^W(y,z)$ for all $z' \in \cV_w^W(y)$. That is, $\cS_w^W(y) \subset H_w^W(y,z)$, which implies the reverse inclusion, as was to be shown.  
 \endproof

\section{Induced pseudogroups}\label{sec-induced}

Let $W\subset \fX$ be an open subset, and let  $\cGF^*(W)$  denote  the restriction of $\cGF^*$ to $W$.  That is, 
\begin{equation}\label{eq-restrpseudogroup} 
\cGF^*(W)  = \left\{ h := (h_\cI,U_h) \in \cGF^* ~|~ \emptyset \ne U_h\subseteq \Dom(h_\cI) \cap W ~{\rm and} ~ h_\cI(U_h) \subseteq W \right\}. 
\end{equation}
 
 The goal of this section is to show  that if   the action of $\cGF$ on $\fX$ is minimal,  then 
 $\cGF^*(W)$ admits a \emph{finite} symmetric generating set  $\cGF^*(W)^{(1)}$, whose elements have ``maximal'' domains in $W$.   
 We begin with some preliminary   considerations. 
Introduce the constants:  
\begin{equation*}\label{eq-transdiamwohol0}
\e_{H,i} < \max \left\{\e \mid  ~\forall ~ x \in \oU_i  ~ {\rm such ~ that} ~   D_{\fM}(x, \eU/4) \subset \oU_i   ,  {\rm then} ~ D_{\fX}(\pi_i(x),\e) \subset \pi_i\left( D_{\fM}(x, \eU/4)\right)\right\}
\end{equation*}
 \begin{equation}\label{eq-transdiamwohol}
\e_H = \min \left\{\e_{H,i} \mid \forall ~ 1 \leq i \leq \nu \right\}.
\end{equation}
Note that $\e_H \leq \eTU$ follows by comparing  the definition \eqref{eq-transdiamwohol} of $\e_H$ with  the definition  \eqref{eq-transdiam}  of $\eTU$.

Recall that $\alpha_W$ is an upper bound on the length of plaque-chains $\cI$ required for a covering of $\fX$ by open sets of the form $h_{\cI}(W)$.
  and that Lemma~\ref{lem-deloneW} shows that  $h_{\cI}$ can be realized by the holonomy along a piecewise geodesic path of length at most  $2\alpha_W \,\dFU$.

Note  that while $\alpha_W$ bounds the lengths of plaque-chains required to obtain a covering of $\fX$ by images of $W$, there is no control over the domains of the holonomy maps defined by the overlaps between the plaques in these chains.   Introduce the integer
\begin{equation}\label{eq-betaW}
\beta_W = \ceil{2 \alpha_W \, \dFU/\eFU} \geq 2 \alpha_W
\end{equation}
where the inequality follows from the inequality $\dFU/\eFU > 1$ of Section~\ref{sec-concepts}. 
By allowing   plaque-chains of length $\beta_W$   we can obtain uniform estimates on the domains of the holonomy maps used, as shown by the following remarks. 
First,  introduce the constant     defined by Proposition~\ref{prop-domest} for $\e = \e_H$, 
 \begin{equation}
\delta_W =  \delta(\e_H, \beta_W) \leq \e_H .
\end{equation}
    Let $h_{\cI} |U_h \in \cGF^*(W) $ where $\cI$ has length at most $\alpha_W$.
 Let $w \in U_h \subset W$ define a leaf $L_w$ without holonomy.  Set   $w' = h(w) \in W$ and so   $w' \in \cN_w^{W}$.  
 Recall that the leaves of $\F$ are complete Riemannian manifolds, so   there exists a geodesic segment $\gamma_{x,y}$ from $x = \tau(w)$ to $y = \tau(w')$ with length $\dF(x,y) \leq   2\alpha_W \,\dFU$. 
By the method of   Section~\ref{subsec-pathstochains}  used in the proof of Proposition~\ref{prop-domest}, we can choose a plaque-chain $\cI'$ covering $\gamma_{x,y}$ with length at most 
$\alpha \leq  \dF(x,y)/\eFU \leq \beta_W$, and such that the holonomy map $h_{\cI'}$   satisfies
  $D_{\fX}(w, \delta_W) \subset \Dom(h_{\cI'})$, and  
\begin{equation}\label{eq-maxdomains44}
h_{\cI'}(D_{\fX}(w, \delta_W)) \subset D_{\fX}(w', \e_H) .
\end{equation}
 As $L_w$ is a leaf without holonomy,   by Lemma~\ref{lem-homotopymin} we have   $[h]_w = [h_{\cI'}]_w$. We thus obtain the main result of this section:

\begin{prop} \label{prop-domaincontrol}
Let $V \subset W$ be an open subset with $\diamX(V) \leq \delta_W$, and let $w \in V$ be such that $L_w$ is without holonomy.  
Then for all $h \in \cGF^*(W)$ and $w \in V \cap U_h$ with $w' = h(w)$ represented by a path $\gamma_{x,y}$  
with length $\dF(x,y) \leq 2\alpha_W \dFU$, then there exists an admissible sequence $\cI'$ with length at most $\beta_W$ such that 
  $V \subset \Dom(h_{\cI'})$ and $[h]_w = [h_{\cI'}]_w$.
\end{prop}
 Introduce the following subset of $\cGF^*(W)$,  
\begin{equation}
\cGF^*(W)^{(1)} = \left\{h_{\cI} \mid \cI = (i_0, i_1, \ldots , i_k) , ~   k \leq \beta_W , ~ 
\Dom(h_{\cI}) \cap W \ne \emptyset , ~ h_{\cI}(W) \cap W \ne \emptyset   \right\}.
\end{equation}
The number of admissible sequences $\cI$ with length at most $\beta_W$ is bounded above by $\nu^{(\beta_W+1)}$, so $\cGF^*(W)^{(1)}$ is a finite collection of maps.

\section{Dynamical partitions of the transversal} \label{sec-partialcoding}

   In this section, we construct  coding functions for a minimal action of the    pseudogroup $\cGF^*$ on $\fX$, and   show some basic  properties for such codings.  
 We note that   Gromov  \cite{Gromov1981,Gromov1987} and   Fried \cite{Fried1987} described applications of the coding technique in dynamics to the study of  actions by finitely-generated groups, as described in 
 the text  by Coornaert and Papadopoulos  \cite{CP1993} which gives an excellent overview.     Our development of  codings for the orbits for a   pseudogroup action  is an extension of this method.

    The coding technique was developed for equicontinuous actions of   pseudogroups  on a Cantor space $\fX$ in the work  \cite{ClarkHurder2013}, where it was also shown that the domains of the coding maps have uniform estimates, which fails in the general case considered  here. Thus,   an important aspect of our development of the coding of the action of the pseudogroup $\cGF^*$ is the careful treatment of the relation between the dynamics of the holonomy maps,   and the domains on which the coding function is defined.
  
 After establishing the basic properties of the coding functions, we give    an inductive method to     
  construct a   nested   sequence of  partitions  of $\fX$ into clopen sets, defined using the  level sets of the coding function.  
  The resulting clopen covers of $\fX$ so defined   yield   ``multi-dimensional'' Kakutani-Rokhlin towers for the action of $\cGF^*$, analogous to those constructed by Forrest in \cite{Forrest2000} for minimal $\mZ^n$-actions.  Our construction here is also related to the   ``zooming'' procedure  in Bellisard, Benedetti and Gambaudo  \cite[Proposition~2.42]{BBG2006} whose proof   uses the   Euclidean geometry of the leaves. The work of this section and Section~\ref{sec-reebslabs}  is also related to the construction    of a ``tower system'' for a tiling space $\Omega_{\bf T}$ as sketched  in the work of   Benedetti and  Gambaudo  \cite[Theorem~3.1]{BG2003},  where the leaves of $\Omega_{\bf T}$ are quotients of a connected Lie group. The method of ``inflation'' described  in Section~3.1 of  \cite{ALM2011} can also be seen as analogous, though the lack of explicit details in these previous works makes  the connections between the approaches to be mostly on an intuitive level.

In this section, we allow $\F$ to have non-trivial holonomy, but work with a leaf without holonomy, which always exists. Fix a point  $w_0 \in \fX$     without holonomy,  and without loss of generality assume that $w_0 \in  \fT_1$.
 Let $L_0$ denote the leaf determined by $w_0$ and let  $\cN_0 = L_0 \cap \cT$ be the induced net.

\subsection{Setting constants}\label{sec-constants}
The analysis of the domains of compositions of elements in $\cGF^*$ requires careful attention to the metric properties of the action. We first fix some basic constants.

Let $0 < \e_1 \leq \e_H$ where $\e_H$  is defined by \eqref{eq-transdiamwohol}. 
Choose  $V_1 \subset \fT_1 \subset \fX$   a clopen neighborhood of $w_0$ with   $\diamX(V_1) < \e_1$. 
Set  $\cN_1 = L_0 \cap \tau(V_1)$, which is a subnet of $\cN_0$.

Let $\alpha_1$ be the word length of the elements needed to cover $\fX$ by translates of  $V_1$, as defined by Lemma~\ref{lem-finitegen}. 
 Then   the generators of the induced pseudogroup    $\cGF(V_1)^{(1)}$ are represented by elements of $\cGF^*$  defined by
 the holonomy along paths of lengths at most $2 \alpha_1 \dFU$, and thus by Proposition~\ref{prop-domaincontrol} are represented by maps   in $\cGF^{(\beta_1)}$,  
 where $\beta_1 = \ceil{2 \alpha_1 \dFU/\eFU}$.

Now let  ${\theta_1} =  (2\alpha_1 +1) \, \dFU $,  then for all $w \in V_1$,  the net  $\cN_w^{V_1} = L_w \cap \tau(V_1)$ is  ${\theta_1}$-dense in $L_w$ by Corollary~\ref{cor-deloneW}.
Moreover, by Lemma~\ref{lem-celldiam} for  each $y \in \cN_w^{V_1}$ the Voronoi cell $\cC_w^{V_1}(y) \subset B_{\F}(y, {\theta_1})$, and by Lemma~\ref{lem-star2}     the star-neighborhood  $\cS_w^{V_1}(y) \subset B_{\F}(y, 3{\theta_1})$.  
Introduce the constants 
\begin{equation}\label{eq-RV1}
  R_1' =  2 \theta_1 + \lF \quad  ,  \quad R_1 = 2 R_1' =  4 \theta_1 + 2\lF ~ .
\end{equation}
Finally, observe that a leafwise path of length $R_1$ in $\fM$ can be covered by a good plaque-chain of length   at most $\whbeta_1 +1$,  where
$\ds \whbeta_1 = \ceil{R_1/\eFU}$. 

  Proposition~\ref{prop-domest}, for the case  $\e = \e_H$ and $\alpha = \whbeta_1$,  yields  a uniform  estimate on the radius of balls contained in the domains of holonomy maps associated to plaque-chains of length at most $R_1$, and such that the images of these balls have diameter at most $\e_H$. As $L_0$ is a leaf without holonomy, the holonomy along a path depends only on the endpoints. As there are only a finite number of such plaque-chains, there are at most a finite number of homotopies required between them as in the proof of Lemma~\ref{lem-homotopy} which will have length at most $\whbeta_1'$.

In     the proof of  Theorem~\ref{thm-specialpartition} below, we construct  a transverse Cantor foliation for $\fM$, based on the   techniques    used in the proof of Theorem~1.3 in  \cite{CHL2013a}. Applying these techniques   requires a (possibly) strong  restriction on  the diameter of the   disks transverse to $L_0$. In particular, the proof of  \cite[Theorem~1.3]{CHL2013a}   shows that for the given bound  $R_1$ on path length distances,  there exists $\delta^*_1 > 0$   such that a Reeb neighborhood  with base given by compact set  of diameter less than $R_1$ and transversal of diameter less than $\delta^*_1$ admits a transverse Cantor foliation.  
Let  $\delta^*_1$ denote this constant.

Note that    $\{V_1, \fX-V_1\}$ defines a clopen partition of $\fX$, and we introduce a constant which bounds the distance between the sets in it. 
Let $\ve_1 = \min \{ \e_H, \dF(V_1, \fX-V_1), \delta^*_1 \}$ if   $V_1 \subset \fX$ is a proper inclusion; otherwise set $\ve_1 = \min\{ \e_H , \delta^*_1\}$. 
Define the constant $\whdelta_1$  using Proposition~\ref{prop-domest},  where 
 \begin{equation}\label{eq-whdelta1}
\whdelta_1 ~ = ~  \delta(\ve_1, \whbeta_1' ) ~ \leq ~  \ve_1 ~. 
\end{equation}
  This constant is used to conclude that the partial coding functions introduced next are well-defined, and  so  assures that the  Reeb neighborhoods introduced in Section~\ref{sec-reebslabs} are well-defined as well.

 \subsection{Coding functions}\label{sec-coding}
  
We begin with a basic notion.
\begin{defn}\label{def-codepartV}
Let   $0 < \e_1' < \diamX(V_1)/2$.
An  \emph{$\e_1'$-coding partition}  of $V_1$  is a collection of disjoint \emph{clopen} sets $\cW_1 =\{W(1;1) ,\ldots,W(1;\kappa_1)\}$  with $w_0 \in W(1;1)$ whose union is $V_1$, satisfying:
\begin{equation}\label{eq-codepartV}
\diamX(W(1;i)) < \min~   \{\e_1' ,  \whdelta_1    \} \quad {\rm for} ~ 1 \leq i \leq \kappa_1.
\end{equation}
\end{defn}

We then have that $\diamX(W(1;i)) <  \diamX (V_1)/2$ so that  $\kappa_1 \geq 2$, thus the partition is not trivial. The assumption   $\diamX(W(1;i)) < \whdelta_1$ implies that if  $w \in W(1;i)$ and $w \in \Dom(\gamma)$ with $||\gamma|| \leq R_1'$, then $W(1;i) \subset \Dom(\gamma)$.

The index set  $\cA_1 = \{0,1,2, \ldots, \kappa_1\}$ of  $\cW_1$ is the ``alphabet''  used to code the action of the induced pseudogroup $\cGF^*(V_1)$ on $V_1$.  This coding  is used to   subdivide the sets of $\cW_1$ into   smaller sets with constant codes, in a way compatible with the dynamics of the induced pseudogroup $\cGF^*(V_1)$.  Let 
\begin{equation}\label{eq-spacingV}
\eta_{1} = \min\left\{~\ve_1,~ d_{\fX}(W(1;i) ,W(1;j))\,~|~\,i \neq j\,~\right\} \quad {\rm and ~set}  \quad \zeta_1=\delta(\eta_1,\whbeta_1')  ~ .
\end{equation}
Note that $\eta_1 > 0$,  as $V_1$ is clopen in $\fX$, so  the constant   $\delta(\eta_1, \whbeta_1')$  is     well-defined by Proposition~\ref{prop-domest}. 

\medskip
 
 With the above preparations, we now give a key construction. Given points $w,w' \in \cN_1$,    let $\gamma_{w,w'}$ be a length minimizing geodesic from $y = \tau(w)$ to $y'=\tau(w')$.    
 Assume that $\dF(y,y') \leq R_1$  then   there is an admissible sequence  $\cI_{w,w'}$   of length $\alpha \leq \whbeta_1'$ which defines  a good plaque-chain covering for $\gamma_{w,w'}$ as in  Definition~\ref{def-goodPC}.   
 Then $\cI_{w,w'}$ yields a holonomy map $h_{\cI_{w,w'}}$ satisfying   $B_{\fX}(w, \whdelta_1) \subset \Dom(h_{\cI_{w,w'}})$. Note that as $L_0$ has no holonomy, the map $h_{\cI_{w,w'}}$  is independent of the choice of $\gamma_{w,w'}$ on the subset   $B_{\fX}(w, \whdelta_1)$ of its domain.

For $w \in \fX$, recall that  $\G_{w}^{V_1,R}$ was introduced in \eqref{eq-groupoidfiltraqtion}.
For simplicity on notation in the following, given $\gamma \in \G_{w}^{V_1,R}$ we  let  $\gamma \equiv h_{\cI_{w,w'}}$ denote the holonomy map defined by the path.

The  \emph{partial coding function} for   $\G_{w}^{V_1,R}$ is defined   as follows. 
\begin{defn}\label{def-coding-1} 
Given   $0 < R \leq R_1$ and  $w \in W(1;i)$ for  $1 \leq i \leq \kappa_1$,  the $\cC_{w,i}^{R}$-code of $u \in W(1;i)$ is the function $\ds   C^{i,R}_{w,u} \colon \G_{w}^{V_1,R} \to \cA_1$ 
defined as
  \begin{equation}
  C^{i,R}_{w,u}(\gamma) =    k  ~ {\rm where} ~  \gamma(u) \in W(1;k).
    \end{equation}
\end{defn}

The value $C^{i,R}_{w,u}(\gamma)$   encodes the   index  $k $ for the path starting at $\tau(u)$, 
shadowing the path $\gamma_{w,w'}$ and terminating at $\tau(u')$.  
By the choice of $\whdelta_1$ and $\e_1'$ above,  $w \in W(1;i)$ implies  $W(1;i) \subset \Dom(\gamma)$ and thus the function $C^{i,R}_{w,u}$ is well-defined.

The next   result implies that the coding partition to be defined next will consist of clopen subsets. Recall that $\zeta_1$ is defined above in \eqref{eq-spacingV}.
\begin{lemma} \label{lem-locconstant}
Let $0 < R \leq R_1$,  $w \in V_1$ and $u,v \in W(1;i)$ with $d_{\fX}(u,v) < \zeta_1$.  Then $C^{i,R}_{w,u}(\gamma) = C^{i,R}_{w,v}(\gamma)$ for all $\gamma \in \G_w^{V_1 ,R}$. Hence, the function $C^{i,R}_w$ defined by $C^{i,R}_w(u) = C^{i,R}_{w,u}$ is locally constant. 
\end{lemma}
\proof
Let  $\gamma \in \G_w^{V_1 ,R}$, and suppose that $u,v \in W(1;i)$ with $d_{\fX}(u,v) < \zeta_1$.
Set $u' = \gamma(u)$ and $v' = \gamma(v)$. Then $d_{\fX}(u,v) < \zeta_1$ implies $d_{\fX}(u',v') < \eta_1$.
If  $u'   \in W(1;i)$ then
 $d_{\fX}(u',v') < \eta_1 \leq d_{\fX}(W(1;i),W(1;j))$ for all $j \ne i$ implies that $v' \in W(1;j)$. 
Thus, $C^{i,R}_{w,u}(\gamma_w) = C^{i,R}_{w,v}(\gamma_w)$.
\endproof

  We next introduce clopen partitions of the sets $W(1;i)$ chosen in Definition~\ref{def-codepartV}, so that    the coding function along paths  of length at most  $R_1'$ as defined by \eqref{eq-RV1} are constant on each set of the partition. 
 For each $1 \leq i \leq \kappa_1$ define the partition $\ds \cV(1;i) = \left\{ V(1;i,1) \cup \cdots \cup  V(1;i,{\kappa_{1,i}})\right\}$ with 
\begin{equation}\label{eq-clopenpart1}
  W(1;i) = V(1;i,1) \cup \cdots \cup  V(1;i,{\kappa_{1,i}})   ~ .
\end{equation}
where the coding function $C^{i,R_1'}_w$ is constant on each set $V(1;i,j)$, and if $j \ne j'$ then $V(1;i,j)$ and $V(1;i,j')$ have distinct codes. 
Assume that $w_0 \in V(1;1,1)$.

\begin{lemma} \label{lem-clopenpart1i}
For $1 \leq i \leq \kappa_1$ and each $1 \leq j \leq {\kappa_{1,i}}$ the set  $V(1;i,j)$ is clopen. Thus, $\cV(1;i) $ is a clopen partition of $W(1;i)$ and the collection 
$\cV(1;1)  \cup \cdots \cup \cV(1;\kappa_1)$ is a clopen partition of $V_1$. 
\end{lemma}
\proof
Lemma~\ref{lem-locconstant} implies that each set  $V(1;i,\kappa_{1,i})$ is open, and as they form a finite cover of the clopen set $W(1;i)$,  they are also clopen sets. \endproof

If the action of $\cGF$ on $\fX$ is equicontinuous, then it may happen that     points  $u,u' \in V(1;i,j)$ have the same code for all paths starting at $w \in V_1$, not just those of length at most $R_1'$. In fact, this fact is used in  \cite{ClarkHurder2013}   to define   clopen partitions which are invariant under the action of $\cGF$. However, if the action of $\cGF^*$ is $\e$-expansive, its action  separates points by at least $\e$, so will have distinct codes if the diameters of the sets in the clopen partition $\cW_1$ are less than $\e$.

 We   introduce a further partition of the clopen sets in $\cV(1;i)$ according to their coding for paths of length at most $R_1 = 2R_1'$.  This additional decomposition    is used  for the patching arguments  appearing in later sections to construct the transverse Cantor foliation for $\F$.

 For $1 \leq i \leq \kappa_1$ and each $1 \leq j \leq \kappa_{1,i}$ let $\ds \cV(1;i,j) = \left\{ V(1;i,j,1) \cup \cdots \cup  V(1;i,j,{\kappa_{1,i,j}})\right\}$ with 
 \begin{equation}\label{eq-clopenpart1ij}
V(1;i,j) = V(1;i,j,1) \cup \cdots \cup  V(1;i,j,{\kappa_{1,i,j}}) 
\end{equation}
where the coding function $C^{i,R_1}_w$ is constant on each set $V(1;i,j,k)$, and if $k \ne k'$ then $V(1;i,j,k)$ and $V(1;i,j,k')$ have distinct codes. 
Assume that $w_0 \in V(1;1,1,1)$.  As before, each set in the partition $\cV(1;i,j) $ is clopen and the sets are disjoint.
   
 Note that  $V(1;i,j)  \subset W(1;i)$ and so satisfies    $\diamX(V(1;i,j)) \leq  \min~   \{\e_1' ,  \whdelta_1    \}$ by \eqref{eq-codepartV}.


 \subsection{Inductive step} \label{subsec-induct}

We next   extend the above  process, so that  the diameters of the transverse partition sets decrease in a systematic manner. 
The results are  used in Sections~\ref{sec-nicestabletransversals}  and \ref{sec-compatiblecantor}    for the construction of  nested   Reeb neighborhoods.
 We first show how to proceed to step 2, which refines the construction above, and the general inductive step follows analogously.
 
   Let $0 < \e_2 < \e_1/2$ be chosen sufficiently small so that 
  the function $\lambda_1({\delta})$ defined in Lemma \ref{lem-separationW}  satisfies   $\lambda_{\e_2} = \lambda_1({\e_2}) > R_1$. 
  
  By choice, we have  $w_0 \in V(1;1,1,1)$, and let   $w_0 \in V_2 \subset V(1;1,1,1)$ be a clopen neighborhood with $\diamX(V_2) < \e_2$. Set $\cN_2 = L_0 \cap \tau(V_2)$,  which is a   subnet of $\cN_1$.
 
 Let $\alpha_2 \geq \alpha_1$ be an integer such that   $\fX$ is covered by translates of  $V_2$ using    words of length at most $\alpha_2$, 
 as given  by Lemma~\ref{lem-finitegen}. 
Define  $\beta_2 = \ceil{2 \alpha_2 \dFU/\eFU}$ and  ${\theta_2} =  (2\alpha_2 +1) \, \dFU$, and as before, set 
 \begin{equation}\label{eq-RV2}
R_2' = 2 {\theta_2} +   \lF   \quad, \quad R_2 = 2 R_2'   \quad, \quad\whbeta_2 = \ceil{R_2/\eFU} ~. 
\end{equation}
 Choose $\whbeta_2'$ so that homotopies of paths  of length at most $\whbeta_2$ contain chains of length at most  $\whbeta_2'$. Let
  \begin{align}\label{eq-nestedreeb} \ve_2 & < \min\{~\delta^*_2,~  d_\fX(V(1;i),V(1,j)),~ d_\fX(V_2,\fX-V_2))~|~i\ne j~\}, \end{align}
where $\delta_2^*$ is given by Corollary \ref{cor-disjoint}, and then set
  \begin{equation}\label{eq-whdelta2}
\whdelta_2 =  \delta(\ve_2, \whbeta_2' ).
\end{equation}

For $\e_2' = \diamX(V_2)/2 < \e_2/2 < \e_1/4$, choose  a clopen partition of $V_2$ as in Definition~\ref{def-codepartV}, given by 
 $\cW_2 =\{W(2;1) ,\ldots, W(2;\kappa_2)\}$  with $w_0 \in W(2;1)$  and which satisfies
$$\diamX(W(2;i)) < \min~   \{\e_2' ,  \whdelta_2    \} \quad {\rm for} ~ 1 \leq i \leq \kappa_2.$$
 This yields a new coding alphabet $\cA_2 = \{1,2, \ldots , \kappa_2\}$ corresponding to $\cW_2$. 
For each $1 \leq i \leq \kappa_2$, use the coding function for the sets in $\cW_2$ to refine $\cW_2$ into its coding partition  for the holonomy along paths of lengths at most $R_2'$, where  $ \cV(2;i) = \left\{ V(2;i,1) \cup \cdots \cup  V(2;i,{\kappa_{2,i}})\right\}$ with 
\begin{equation}\label{eq-clopenpart2}
W(2;i) = V(2; i,1) \cup \cdots \cup  V(2;i,{\kappa_{2,i}}) .
\end{equation}
We assume that $w_0 \in V(2;1,1)$, and observe that  the clopen sets satisfy 
 $$V(2; i,j) \subset W(2;i) \subset V(1;1,1,1) \subset W_1 \subset V_1 ~ .$$  
In particular, this implies that     $\diamX(V(2;i,j)) \leq  \min~   \{\e_2' ,  \whdelta_2    \} \leq \e_H/4$.

We also introduce the subdivision of these partitions   for the holonomy of paths with length at most $R_2 = 2 R_2'$.
 For $1 \leq i \leq \kappa_2$ and  $1 \leq j \leq \kappa_{2,i}$,  let $\ds \cV(2;i,j) = \left\{ V(2;i,j,1) \cup \cdots \cup  V(2;i,j,\kappa_{2,i,j})\right\}$ with 
$$V(2;i,j) = V(2;i,j,1) \cup \cdots \cup  V(2;i,j,{\kappa_{2,i,j}}) ~ , ~ w_0 \in V(2;1,1,1)$$
where the coding function   is constant along paths of length at most $R_2$ on each set $V(2;i,j,k)$, and if $k \ne k'$ then $V(2;i,j,k)$ and $V(2;i,j,k')$ have distinct codes.   As before, each set in the partition $\cV(2;i,j) $ is clopen and disjoint.

 The above process can now be  repeated recursively to obtain:
 
 \begin{prop}\label{prop-dynamicpart}
Let  $w_0 \in \fT_1 \subset \fX$   such that the leaf $L_0$ it determines is without holonomy.    
 Then each $\ell \geq 1$, there exists:
 \begin{enumerate}
\item  \label{list-partitem1}  clopen neighborhoods  $w_0 \in   V_{\ell} \subset V_{\ell -1} \subset \cdots \subset V_1 \subset \fT_1$ ~ with ~  $\diamX(V_i) < \e_H/2^i$
\item \label{list-partitem2} constants $\e_{\ell}$,   $\alpha_{\ell}$, $\theta_{\ell}$, $R_{\ell}$, $\whbeta_{\ell}$, $\whdelta_{\ell}$ such that $\lambda_{\e_{\ell}} = \lambda_1(\e_\ell) > R_{\ell -1}$
\item \label{list-partitem3} nets $\cN_{\ell} = L_0 \cap \tau(V_{\ell})$ which are $\lambda_{\e_\ell}$-separated and $\theta_{\ell}$ dense
\item \label{list-partitem4} a clopen partition $\cW_{\ell}$ of $V_{\ell}$ with elements $W(\ell;  i)$ labeled by $\cA_{\ell} = \{1,2, \ldots , \kappa_{\ell}\}$
\item  \label{list-partitem5} a clopen partition $\cV(\ell ; i)$ of each $W(\ell, i)$ with elements $V(\ell; i,j)$ whose elements are clopen sets with  $\diamX(V(\ell; i,j)) < \min\{\whdelta_{\ell}, \e_H/4^{\ell}\}$ for $1 \leq j \leq \kappa_{\ell, i}$, and the coding function is constant on each $V(\ell; i,j)$ along paths of length at most $R_{\ell}'$
\item \label{list-partitem6} a clopen partition $\cV( \ell ; i,j)$ of each clopen set $V(\ell; i,j)$ such that the coding function is constant on each $V(\ell; i,j,k)$ along paths of length at most $R_{\ell} = 2 R_{\ell}'$ for $1 \leq k \leq \kappa_{\ell,i,j}$.
\end{enumerate}
\end{prop}

The clopen sets $V(\ell; i,j)$ of the partition $\cV(\ell;i)$ will be used to create Reeb neighborhoods in the following Sections~\ref{sec-reebslabs} and \ref{sec-nicestabletransversals}, by translating them along the holonomy of paths of length at most $R_{\ell}'$.

Note that this technique of decomposing  a descending chain of clopen transversals   $V_{\ell}$ for $\ell \geq 1$ into code blocks $\cV(\ell;i)$ is \emph{implicit} in the studies of $\mR^n$-actions by Bellisard, Benedetti and Gambaudo  \cite[Proposition~2.39]{BBG2006} and  by Forrest in \cite{Forrest2000}, and  in the sketch of the proof of Theorem~3.1 in  \cite{BG2003} for case when the leaves of $\F$ are defined by a transitive action of a connected Lie group.

\section{Transverse Cantor foliations}\label{sec-reebslabs}

A \emph{Cantor foliation} $\cH$ on a continuum $\Omega$ is a ``continuous decomposition'' of it into Cantor sets, which are the ``leaves'' of $\cH$.   
Cantor foliations arise naturally in many dynamical and geometric contexts.  In this work, we require Cantor foliations  which  are ``transverse'' to the foliation $\F$ of a  matchbox manifold $\fM$, as made precise in Definition~\ref{def-cantorfol} below.  
All constructions in the literature  which build up inverse limit representations for a particular  class of matchbox manifolds, such as the Williams solenoids, tiling spaces, and suspensions of minimal actions on Cantor sets  \cite{APC2011,AP1998,BDHS2010,BBG2006,BG2003,ClarkHurder2013,CHL2013a,LR2013,Sadun2003,EThomas1973}  use such a transverse Cantor foliation either explicitly or implicitly. The  branched manifold   quotients of a matchbox manifold are obtained via transverse projections along subsets of the leaves of $\cH$ that are the ``axes of collapsing'' in foliated compact subsets of  $\fM$. As such,  the existence of such $\cH$ is fundamental to the proof of Theorem~\ref{thm-main1}.

The existence of a transverse Cantor foliation $\cH$ can be part of the given data, such as when a matchbox manifold is a fiber bundle over a closed manifold whose fibers are Cantor sets. In this case, the leaves of $\cH$ are given by Cantor fibers of the bundle.  
For a tiling space $\Omega_{\bT}$ formed from a   tiling of $\mR^n$ with finite local complexity, 
the foliation by tilings is defined by a free action of $\mR^n$ on $\Omega_{\bT}$, 
and the Cantor foliation is defined by the choice of punctures in each of the proto-tiles,   then taking their translates under the $\mR^n$-action.   Finally,  a Cantor foliation may be  determined by dynamical properties of the ambient space, such as in the case   of Williams solenoids \cite{Williams1967,Williams1974} where the Cantor foliation is  given by intersections   of the (unstable) attractor with leaves of the stable foliation. This dynamical construction has been generalized  to the notion of  a ``Smale space'' as introduced by Ruelle \cite{Ruelle1988}, and studied in the  works of Putnam \cite{APSG2013,Putnam1996,PS1999,Putnam2013} and    Wieler \cite{Wieler2012b}.

For the general  matchbox manifold,  the existence of  a Cantor foliation on the ``Big Boxes'' in $\fM$ is shown in  the authors' work \cite{CHL2013a},
where this is proven using       local constructions based entirely on intrinsic  techniques.  
The results of that paper are used in the rest of this section, and in Sections~\ref{sec-nicestabletransversals} and \ref{sec-compatiblecantor}, where   we extend the results   to show that a Cantor foliation can always be defined on a minimal matchbox manifold $\fM$ without holonomy.
Moreover, the Cantor foliation on $\fM$ can be constructed to be  compatible   with the constructions in Section~\ref{sec-partialcoding}.  

In Section~\ref{sec-defCantor}, we give a rigorous definition of a Cantor foliation, then in Section~\ref{sec-defBigBox} we  recall a main result proved in \cite{CHL2013a} which gives the existence of a Cantor foliation on a   neighborhood of a compact \emph{path-connected} subset $K_x$ in $\fM$,  called a \emph{Reeb neighborhood}. Section~\ref{sec-cantorreeb} shows the existence of a covering of $\fM$ by ``Big Boxes'', which are Reeb neighborhoods with a transverse Cantor foliation.

\subsection{Definition of transverse Cantor foliations} \label{sec-defCantor}

We assume there is given a   regular covering $\{U_{i} \mid 1 \leq i \leq \nu\}$ of $\fM$ by foliation charts, as in Proposition~\ref{prop-regular}, with charts  
  $\vp_i \colon \oU_i \to [-1,1]^n \times \fT_i$ where $\fT_i \subset \fX$ is a clopen subset. Moreover,  by construction, each chart admits a foliated extension $\whvarp_i \colon \whU_i \to (-2,2)^n \times \fT_i$ where $\oU_i \subset \whU_i \subset \fM$ is an open neighborhood of the closure $\oU_i$ and $\whvarp_i | \oU_i = \vp_i$.
 For a clopen set $V \subset \fT_i $ set
 \begin{equation}\label{eq-saturatedset}
\fU^{V}_{i} = \pi_{i}^{-1}(V) \subset \oU_{i} \quad ; \quad   \whfU^{V}_{i} = \whvarp_i^{-1}(V) \subset \whU_{i} ~ .
\end{equation}

\medskip

  \begin{defn}\label{def-cantorfol}
Let $\fM$ be a matchbox manifold, and  $\fB \subset \fM$ a closed subset. An equivalence relation $\approx$ on $\fB$ is said to define a \emph{transverse Cantor foliation} $\cH$ of $\fB$ if for each 
 $x \in \fB$, the class $\cH_x = \{y \in \fB \mid y \approx x\}$ is a Cantor set. Moreover, assume that for each $x \in \fB$,    there exists:
\begin{enumerate}
\item   $1 \leq i_x \leq \nu$ with $x \in U_{i_x}$, 

\item  a clopen subset $V_x \subset \fT_{i_x}$ with $w_x = \pi_{i_x}(x) \in V_x$ and $\fU^{V_x}_{i_x} \subset \fB$;

\item    a homeomorphism into $\Phi_{x} \colon [-1,1]^n \times V_{x} \to \whU_{i_x}$ such that for the point $w_x$ we have
$$\Phi_{x}(\xi , w_x) = \whvarp^{-1}(\xi , w_x)  ~ {\rm for} ~   \xi \in (-1,1)^n ,$$
\item   for   $\xi \in (-1,1)^n$ and  $z = \whvarp^{-1}(\xi , w_x)$, the image $\Phi_x(\{ \xi \}  \times V_{x}) = \cH_z \cap \whfU^{V_x}_{i_x}$.
\end{enumerate}
The leaves of the ``foliation'' $\cH$ are defined to be the equivalence classes $\cH_x$ of ~ $\approx$ in $\fB$.
\end{defn}

We give some additional remarks on  the conditions in the definition above. First, the index $i_x$ in Condition~\ref{def-cantorfol}.1 selects a foliation chart containing $x$ in its interior. Then      Condition~\ref{def-cantorfol}.2 specifies a clopen subset $V_x$  which defines a ``plaque''  in the leaf   $\cH_x$, so that $\fB$ is a union of such plaques. 

Conditions~\ref{def-cantorfol}.3  and  \ref{def-cantorfol}.4 specify the leaf $\cH_x$ containing $x$ as the graph of the function $\Phi_x \colon \{\xi\} \times V_x \to  \whU_{i_x}$ where $\xi = \lambda_{i_x}(x) \in (-1,1)^n$ is the horizontal coordinate of $x$ for the coordinates $\vp_{i_x}$.
In particular, as $V_x$ is clopen, it is   a Cantor set, so the leaf $\cH_x$ is a Cantor set.

Note that  the function $\Phi_{x}$ depends on  the index $i_x$ as the leaves of $\cH$ are not assumed to coincide with the local ``vertical'' foliation of 
$U_{i_x}$ defined by the transverse coordinate.  The function $\Phi_x$ gives the necessary ``adjustment'' to the local vertical foliation.  Also note, the   image of the graph of $\Phi_x$ through $z \in  \oU_{i_x}$ is contained in the open neighborhood $\whU_{i_x}$ of $\oU_{i_x}$ but  need not be contained in $\oU_{i_x}$. 

 Recall that the constants $0< \dFU < \lF/5$ were introduced in Proposition~\ref{prop-regular}, and there is a   constant $\lF^* \leq \lF/5$ introduced in 
Section~15.4 of \cite{CHL2013a}, where $\lF^*$ is chosen so that the leafwise balls of radius $\lF^*$ are ``approximately Euclidean'' as required in the constructions there.  We can also  assume without loss of generality  that  $\lF^* < \eFU < \dFU$.

\subsection{Reeb neighborhoods} \label{sec-defBigBox}

Cantor foliations as in Definition~\ref{def-cantorfol} are constructed using the ``Big Box'' Theorem   \cite[Theorem~1.3]{CHL2013a} as follows.  
First we define    \emph{Reeb neighborhoods} as in   \cite[Section~6]{CHL2013a}. 

For $x \in \fM$,  let $L_x$ be the leaf containing $x$, and set  $w = \pi_{i_x}(x) \in \fX_{i_x}$ for some $1 \leq i_x \leq \nu$. We assume that $L_x$ is without  holonomy.   

 Let $Q_x \subset L_x$ be a compact connected subset with $x \in int(Q_x)$ and assume that $Q_x$ is a union of plaques. Then    $Q_x \cap \cT$ is a net in $Q_x$, and for each   $z \in Q_x \cap \cT$,  there is a path $\gamma_z$ in $Q_x$ between $x$ and $z$. Let   $h_{x,z}$ be  
the  holonomy homeomorphism  along $\gamma_{x,z}$ as constructed  in Section~\ref{subsec-pathstochains}. 

For     a clopen neighborhood $V_x \subset \fX_{i_x}$ of $w$, we assume that    $V_x \subset \Dom(h_{x,z})$ for each  $z \in Q_x \cap \cT$.
  Then set  $V_x^z = h_{x,z}(V_x)$, which is well-defined, as $L_x$ has no holonomy. 
 Let   $1 \leq i_z \leq \nu$ denote  the index of the foliation chart containing $z$, which is the last index for the plaque-chain between $x$ and $z$  used to define $h_{x,z}$. 
The \emph{Reeb neighborhood}   of $Q_x$ associated to   $V_x$ is   a union of saturated neighborhoods, defined  as in \eqref{eq-saturatedset}:
\begin{equation}\label{eq-Reeb}
\fN_{Q_x}^{V_x} = \bigcup_{z \in Q_x \cap \cT} \fU^{V_x^z}_{i_z}.
\end{equation}
Note that each path-connected component in $\fN_{Q_x}^{V_x}$ is a union of plaques. Here is one of the main results of  \cite{CHL2013a}. 

\begin{thm}[Big Box]\label{thm-tessel}
Let $\fM$ be a matchbox manifold, $x \in \fM$ and suppose that $L_x$ is a leaf without holonomy. 
Let  $Q_x$ be a compact connected subset of  $L_x$ with $x \in int(Q_x)$  and assume that $Q_x$ is a union of plaques. Then there exists   a 
 clopen set $V_x \subset \fX$, a  closed connected subset $\widehat{Q}_x$ which is a union of plaques with     $Q_x \subset int(\widehat{Q}_x)$, and    a foliated homeomorphic inclusion $\Phi \colon  \widehat{Q}_x \times V_x \to \fM$   such that the images $\ds \Phi   \left\{ \{y\} \times V_x \mid y \in \widehat{Q}_x \right\}$ form a continuous family of Cantor transversals on a neighborhood of  $\F | \fN_{Q_x}^{V_x}$.
\end{thm}
\proof 
We give a sketch of the proof from \cite{CHL2013a}, so that we can refer to the steps involved in extending this result in the following. Given $Q_x$, we want to choose a neighborhood $V_x$ small enough so that Theorem \ref{thm-tessel} follows from Theorem~1.3 of \cite{CHL2013a}.  

We require some notations and constants from   \cite{CHL2013a}. The reader can just take these as given, as the reasons for and details of their choices is part of the most technical aspects of the work \cite{CHL2013a}.
The constant $\ve_0>0$ defined in \cite{CHL2013a} is chosen so that there is a prescribed bound on the metric distortions in charts. This is used in  \cite[Section~15.5]{CHL2013a} to define a   ``transverse'' constant $r^*>0$, such that  if $x$ and $y$ are points in the same plaque, and $x'$ and $y'$ are points in a plaque in the same chart which is at transverse distance at most $r^*$ from the given one, and $x'$ and $y'$ have the same coordinates as $x$ and $y$ respectively, then the metric distances between $x$ and $y$, and between $x'$ and $y'$,  differ by at most $\ve_0 \lF^*$, where $\ve_0$ is very small and depends on the geometry of leaves. 
The definition of $\ve_0 > 0$ as given in \cite[Section~15.3]{CHL2013a}  is the most delicate part of the estimates in that work. In particular, we have that 
$\ds   \ve_0 < 1/2000$ and $\ds  \lF^*< \eFU$, so  that $\ds \ve_0 \lF^* < \eFU/2000$.

We require the following   consequence of the choice of $r^*$ and $\ve_0$.  Suppose that $z \in  U_i \cap  U_j$ for $i \ne j$.
Then the  divergence $div(z,i,j,r) \leq \ve_0\lF^*$ as defined in \cite[Section~14.1]{CHL2013a}. This technical statement has a simple geometric interpretation. 
The set of points in $U_i$ with the same horizontal coordinate as $z$ defines a ``standard section''   $\fZ_{z,i} \subset U_i$ (see  \cite[(17.2)]{CHL2013a}). Likewise, $z$ defines a standard section $\fZ_{z,j} \subset U_j$. Let $x' \in \fZ_{z,i}$ and $y' \in \fZ_{z,j}$ lie on the same plaque in $U_i \cap  U_j$ and suppose this plaque  is within $r^*$ from the plaque containing $z$. Then the leafwise distance   $\ds \dF(x',y') \leq \ve_0\lF^* < \eFU/2000$.

 For $Q_x$ given, let  $R$ denote its    diameter so that $Q_x \subset D_{\F}(x,R)$.    By   Proposition~\ref{prop-regular}.\ref{item-uniform}, the radii   of the plaques in the atlas $\cU$ are uniformly bounded by $\dFU$, so the   saturation of $Q_x$ is contained in the ball $D_\F(x,R + 4\dFU)$. Form the closed plaque-saturation of $D_\F(x,R + 4\dFU)$  and denote the set obtained   by $\widehat{Q}_x'$. Then the diameter of $\widehat{Q}_x'$ is bounded by $R + 8\dFU < R + 2\lF$.

Apply Proposition \ref{prop-domest} for $\epsilon = r_*/2$ to conclude that there exists $\delta_* = \delta(r_*/2, R + 2\lF)$ such that if $x \in \widehat{V}_x$ is a clopen neighborhood with $\widehat{V}_x \subset B_\fX(w_x,\delta_*)$, then for any path $\gamma$ with initial point $x$ and length at most $R + 2\lF$, the holonomy translate $h_\gamma(\widehat{V}_x)$ of the set $\widehat{V}_x$ has diameter less than $r_*/2$.

For $\widehat{V}_x \subset B_\fX(w_x,\delta_*)$ then 
by the ``Big Box'' Theorem \cite{CHL2013a} there exists a Cantor foliation on a subset $\widehat{Q}_x \subset \widehat{Q}_x'$ with  $ Q_x \subset int(\widehat{Q}_x)$  and   there is a continuous map $\Phi_x: \widehat{Q}_x \times \widehat{V}_x \to \fM$. Since $\widehat{Q}_x \subset L_x$ is a subset of a leaf without holonomy, the map $\Phi_x$ is injective on $\widehat{Q}_x$, and so there exists $V_x \subset \widehat{V}_x$ such that the restriction $\Phi_x|_{\widehat{Q}_x \times V_x}$ is a homeomorphism onto its image. As $Q_x \subset \widehat{Q}_x$, this completes the proof of Theorem~\ref{thm-tessel}.
\endproof

\subsection{``Big Box'' Coverings} \label{sec-cantorreeb}
We next use  Theorem~\ref{thm-tessel} to construct      coverings of $\fM$ by ``Big Boxes'' which are adapted to the coding partitions of $\fX$ defined in Section~\ref{sec-partialcoding}. This will use the   constant  $\dFU$   introduced in Proposition~\ref{prop-regular},  the constant   $\e_H$   as defined by \eqref{eq-transdiamwohol}, and $\lF^*$  introduced in Section~\ref{sec-defCantor}.

\begin{thm}\label{thm-specialpartition}
Let $\fM$ be a minimal matchbox manifold without holonomy, and there is given a clopen subset    $V_1 \subset \fX$ with $diam(V_1) < e_H/2$.
Choose $w_0 \in V_1$, and let $L_0$ denote the leaf containing $x_0 = \tau(w_0)$. Then there 
  exists a constant $\delta_1^*>0$ such that for any clopen partition $\cW_{1} = \{W(1;1), \ldots, W(1;\kappa_1)\}$ of $V_1$ with $diam(W(1;  i))< \delta_1^*$ for $1 \leq i \leq \kappa_1$, there exists: 
  \begin{enumerate}
\item     compact connected  subsets $\{ K_1,K_2,\ldots,K_{\kappa_1}\}$ of  $L_0$;
\item   for each $1 \leq i \leq  \kappa_1$, a compact connected subset   $\whK_i \subset L_0$ which is the closure of a union of plaques, such that $K_i \subset int(\whK_i)$;
\item  for each $1 \leq i \leq  \kappa_1$, a foliated homeomorphic inclusion $\ds  \Phi_i \colon  \whK_i \times W(1; i)  \to \fM$  which form a continuous family of Cantor transversals for $\F | \fN_{K_i}^{W(1;i)}$; 
\end{enumerate}
such that the interiors $\ds U(1,i) \equiv   \Phi_i\left(int(K_i) \times W(1; i)\right)$ form an open  covering of $\fM$.
\end{thm}
\proof

The idea of the proof is use the Voronoi cells defined by the net $\cN_{w_0}^{V_1} = L_0 \cap \tau(V_1)$ to construct Reeb neighborhoods as in \eqref{eq-Reeb}     which cover $\fM$, and then choose an appropriate subcover.   

We assume without loss of generality that $V_1 \subset \fT_1 \subset \fX$, so that $\tau(V_1) = \tau_1(V_1)  \subset \cT_1$.

Let $\alpha_1$ be the integer defined by  Lemma~\ref{lem-finitegen} for the clopen set $V_1$, and let ${\theta_1} =  (2\alpha_1 +1) \, \dFU $ be the density constant for the net $\cN_{w_0}^{V_1}$ as established in Corollary~\ref{cor-deloneW}. Set     $R_1' =  2 \theta_1 + \lF$ and $R_1 = 2R_1'$.

  For $w \in V_1$, let $\cO(V_1, w) = \pi(\cN_{w}^{V_1}) \subset V_1$ denote the   orbit  of $w$ in $V_1$ which is dense by minimality.   
Then for  $u \in \cO(V_1, w)$, set $z = \tau(u) \in L_w$,   and  recall that $\cC^{V_1}_{w}(z)$ denotes the Voronoi cell in $L_w$ as defined by \eqref{eq-voronoi2} which contains $z$ in its interior. 
   Lemma \ref{lem-star2} implies for all  $w \in V_1$ and $z \in \cN_w^{V_1}$, 
\begin{equation}
\cC^{V_1}_{w}(z) \subset \cS^{V_1}_{w}(z) \subset B_\F(z, 3   \theta_1) ~ .
\end{equation} 
 By assumption, for each $w \in V_1$    the leaf $L_w$ has no holonomy. Thus,  given  $y,z\in \cN_{w}^{V_1}$ with $u =\pi(y) \in V_1$ and $v = \pi(z)$, there is a geodesic path $\gamma_{y,z}$ in $L_x$ joining them which defines a holonomy transformation $h_{u,v}$ defined on a sufficiently small open neighborhood of $u$.  
More precisely,   Lemma~\ref{lem-pathlengths} implies that  for    $\e' = \dFU/4 R_1$,  there exists $\delta_1' > 0$ so that if $\dX(u, u') < \delta_1'$ then for a length-minimizing geodesic $\gamma$ contained in $B_\F(x, R_1)$,    there is a shadowing piecewise geodesic $\gamma'$ as in the proof of Lemma~\ref{lem-pathlengths} from $y'$ to $y'$ with the length estimate  
\begin{equation}\label{eq-lengthbound}
 \| \gamma\| - \dFU/4 ~ \leq ~  \| \gamma'\| ~ \leq ~ \|  \gamma + \dFU/4\| ~. 
\end{equation}
For the constant $\e = \min\{\e_H/2, \delta_1'/2\}$,    Proposition \ref{prop-domest} implies 
  there exists $\delta_1''>0$ such that  for $y,z \in \cN_{w}^{V_1}$,
$\dF(y,z) \leq R_1$ implies  $B_{\fX}(u,\delta_1'') \subset \Dom(h_{y,z})$, 
   and  $h_{y,z}(B_{\fX}(u,\delta_1'')) \subset B_{\fX}(v, \e)$. 
  
 Next, observe that   
 $\ds \pi \left( B_\F(x, R_1') \cap \cT \right) \subset   \cN_w^{V_1} \subset \fX$ is a finite subset for all   $w \in V_1$.
As  $V_1$ is compact,  its    cardinality has a uniform upper bound    for $w \in V_1$. It follows that there exists $\delta_1''' \leq \delta_1''$ so that for all $w \in V_1$, then for a clopen neighborhood  of $w$, $V \subset V_1$ with $\diamX(V) \leq \delta_1'''$,    
  \begin{equation}
\left\{ h_{w,u}(V)   \mid u \in \pi \left( B_\F(x, R_1') \cap \cT \right)  \right\} 
\end{equation}
forms a disjoint collection of clopen sets in $\fX$. Thus,  there is a   Reeb neighborhood of the plaque-saturation of $B_\F(x, \theta_1 + \lF)$, defined  as in    \eqref{eq-Reeb}.

\begin{defn}\label{def-defbase}
For each $x \in \cN_w^{V_1}$, let    $K_x \subset L_w$ denote the closure of the plaque-saturation of $B_\F(x, \theta_1 + \dFU)$, and let 
 $\whK_x$ denote the closure of the plaque-saturation of $B_\F(x, \theta_1 + \lF)$, so that  $ K_x\subset int(\whK_x)$ as $4\dFU < \lF$. 
\end{defn}

We next define a constant $\delta_1 > 0$ which sets the scale for the coding partition of $V_1$. 
  Theorem \ref{thm-tessel} for $Q_x = \whK_x$ implies that there exists $0 < \delta_1 \leq \delta_1'''$ so that if $W \subset V_1$ is a clopen neighborhood $w \in V_1$   with $\diamX (W) < \delta_1$, and $\whK_x$ is defined as previously for $x = \tau(w)$,  then the corresponding Reeb neighborhood $\fN_{\whK_x}^{W}$ admits a Cantor foliation. 
 Moreover, from the proof of Theorem \ref{thm-tessel} in \cite{CHL2013a},  the choice of $\delta_1$ is uniform for all $w \in V_1$. We use this uniformity as follows.

For  $w_0 \in V_1$ given, consider  the collection  $\ds \{ B_{\fX}(w, \delta_1/4) \cap V_1 \mid w \in \cO(V_1, w_0)  \}$.
This forms   a covering of $V_1$ by open sets, as $\cO(V_1, w_0)$ is dense in $V_1$. We modify this covering to obtain a covering of $V_1$ by clopen subsets.
Note that the  closure $D_{\fX}(w, \delta_1/4) = \overline{B_{\fX}(w, \delta_1/4)}$ is compact with   $D_{\fX}(w, \delta_1/4) \subset  B_{\fX}(w, \delta_1/2)$, so there exists a finite covering of $D_{\fX}(w, \delta_1/4)$ by clopen subsets of $B_{\fX}(w, \delta_1/2)$. 
Apply this remark to each set in the covering $\ds \{ B_{\fX}(w, \delta_1/4) \cap V_1 \mid w \in \cO(V_1, w_0)  \}$, and we 
  obtain a countable covering of $V_1$ by \emph{clopen} sets with diameters less than $\delta_1$. 
   Since $V_1$ is compact, we can choose  a finite subcover $\{ W_{1},\ldots, W_{m} \}$ of this clopen cover.     
   This finite subcovering need not be disjoint, but  we can use it to form a   \emph{disjoint clopen} covering as follows:
   \begin{align}\label{eq-hatcover} 
  \widehat{W}_1& = W_{1}, ~ \widehat{W}_2 = W_{2} - W_{1}, \ldots , \widehat{W}_i = W_{i} - \bigcup_{k=1}^{i-1} W_{k},  \ldots
  \end{align}
If some of $\widehat{W}_i$ are empty, we discard them and renumber the remaining sets, thus obtaining a clopen partition 
$\ds \widehat{W}_1 = \{\widehat{W}_1,\ldots,\widehat{W}_s\}$ of $V_1$ where $s \leq m$.  Then define: 
  \begin{align*} 
  \delta_1^* < \min\{~\delta_1, ~d_\fX(\widehat{W}_i, \widehat{W}_j)~|~i,j=1,\ldots, s~\}.
  \end{align*}
Next, choose a clopen partition   $\cW_1 = \{W(1;1), \ldots, W(1;\kappa_1)\}$ of $V_1$   with $\ds \diamX(W(1;i)) \leq  \delta_1^*$.
 Note that  every subset of $V_1$ of diameter less than $\delta_1^*$ must be contained in one of the sets $\widehat{W}_i$, for some $1 \leq i \leq s$, 
 and so   for each   $1 \leq i \leq \kappa_1$,  there is a unique $1 \leq m_i \leq s$ with  $W(1;i) \subset \widehat{W}_{m_i}$.

For  $1 \leq i \leq \kappa_1$, choose $x_i \in \tau(W(1;i)) \cap L_0$ and set $w_i = \pi(x_i) \in W(1;i)$.

Then set $K_i$ to be the closed plaque-saturation of  $B_{\F}(x_i, \theta_1 + \dFU)$, and let $\whK_i = \whK_{x_i}$ be   defined  as in Definition~\ref{def-defbase}.   
By the choice of $\delta_1^*$ above we conclude that there exists $\whK_i'$ with  $\whK_i \subset int(\whK_i')$ and a foliated homeomorphic inclusion
$\ds  \Phi_i \colon  \whK_i' \times W(1; i)  \to \fM$  which restricts to  a continuous family of Cantor transversals for $\F | \fN_{\whK_i}^{W(1;i)}$

  To conclude the    proof of Theorem~\ref{thm-specialpartition}, it 
 remains to show that the sets $\ds  int\left(\Phi_i \left\{ K_i \times W(1; i)\right\}\right)$ yield an open of covering $\fM$.
 Recall that for any $w \in V_1$ the leaf $L_w$ admits a Voronoi decomposition as in \eqref{eq-Vdecomposition}, so it suffices to show that for each $y \in \cN_w^{V_1}$ the cell $\cC_w^{V_1}(y) \subset  int\left(\Phi_i(K_i \times W(1; i))\right)$ for some $1 \leq i \leq \kappa_1$. 
 
 Now let  $y \in \cN_w^{V_1}$ with $w = \pi(y) \in V_1$ then there exists   $1 \leq i \leq \kappa_1$ for which $w \in W(1;i)$. 
 Then $\cC_w^{V_1}(y) \subset B_{\F}(y, \theta_1)$ by the choice of $\theta_1$, so it suffices to show that 
 $\ds B_{\F}(y, \theta_1) \subset \Phi_i(K_i \times W(1; i))$. 
 
 By construction, we have that $\ds B_{\F}(x_i, \theta_1 + \dFU) \subset \Phi_i(K_i \times W(1; i))$. By the choice of $\delta_1'$ above, given a   geodesic path starting at     $y \in \tau(W(1,i))$ with length at most $\theta_1$, there is a piecewise-geodesic path    starting at $x_i$ with length at most $\theta_1 + \dFU/4$ which shadows the given geodesic. The endpoint of the shadowing curve  is contained in  $\ds B_{\F}(x_i, \theta_1 + \dFU)$, which yields the inclusion  $\ds B_{\F}(y, \theta_1) \subset \Phi_i(K_i \times W(1; i))$ as was to be shown.  
\endproof

A similar argument may of course be repeated for any clopen subset $V_\ell \subset V_1$.  We assume the clopen set $V_{\ell} \subset V_1$ is given with $w_0 \in V_{\ell}$. Let ${\theta_{\ell}} =  (2\alpha_{\ell} +1) \, \dFU $ be the density constant for the net $\cN_{\ell} = L_0 \cap \tau(V_{\ell})$ as established in Corollary~\ref{cor-deloneW}. Set 
   $R_{\ell} = 4 \theta_{\ell} +  \lF$ and  $R_{\ell}' =  2 \theta_{\ell} + \lF$.

\begin{cor}\label{cor-disjoint}
 There exists a constant $\delta_\ell^*>0$ such that for any clopen subset $W\subset V_\ell$ of $\diamX(W)< \delta_\ell^*$ and $w \in W$,  the collection  
$\ds \left\{ h_{w,u}(V)   \mid u \in \pi(B_\F(x, R_{\ell}') \cap \cT)  \right\}$
forms a disjoint collection of clopen sets in $\fX$.
\end{cor}
 
 \begin{remark}
 {\rm 
 Observe that in the proof above, we cover the leaves of $\F$ with plaque-saturations of the balls $B_{\F}(y, \theta_1 +\dFU)$ centered at points  $y \in V_1$. 
 The fact that these sets provide a covering of $\fM$ is a consequence of Corollary~\ref{cor-deloneW} and the definition of $R_1'$, which implies that, for each $y \in \cN_w^{V_1}$,  the ball  $B_{\F}(y, \theta_1)$  contains the Voronoi cell $\cC_w^{V_1}(y)$,   and these Voronoi cells form a covering of the leaf $L_w$. 
  Note that the balls  $B_{\F}(y, \theta_1 )$ do not necessarily provide an efficient covering, as for $y,z \in  \cN_w^{V_1}$ the intersection   $B_{\F}(y, \theta_1 ) \cap  B_{\F}(z, \theta_1 )$ may be ``large''. It would be more efficient to use   small neighborhoods of the Voronoi cells $\cC_w^{V_1}(y)$ themselves to form the covering, except that in the construction of the Reeb neighborhoods, there is no guarantee that the Cantor foliation of the neighborhoods of the form 
 $\ds \Phi_i \left\{K_i \times W(1; i)\right\}$ has bounded distortion when used to translate  the Voronoi cells between leaves.   }
 \end{remark}

\section{Existence of  nice stable transversals} \label{sec-nicestabletransversals}

Theorem~\ref{thm-specialpartition} shows that a minimal matchbox manifold $\fM$ without holonomy has a finite covering by Big Box neighborhoods, each of which is endowed with     a transverse Cantor foliation.  However, on the overlaps of the sets in this covering, the Cantor foliations need not agree (see Figure \ref{fig:overlaps3}), so the union of the leaves of these foliations does not define a Cantor foliation on all of $\fM$. In this section, we give an inductive procedure for modifying the construction of the transverse Cantor foliations on the Big Boxes in a covering, to obtain a Cantor foliation defined on all of $\fM$, which is also adapted to a given initial coding partition on the transversal.

 The approach we take here, is to recall the construction of the transverse Cantor foliations on the Reeb neighborhoods from the authors' work  \cite{CHL2013a}, which  reduces the construction to the existence of a \emph{nice stable transversal} on each, which is a finite collection 
 $\cX = \{ \fZ(\xi,i_z,W(1;i)^z) \mid 1 \leq i \leq p\}$ 
 of transversals  which intersect each leaf in a net whose  separation and spanning constants $(d_1 , d_2)$ are sufficiently small. The  leaves of the  Cantor foliation on each Big Box neighborhood are   then defined as collections of points with the same barycentric coordinates in each simplex of the foliated simplicial decomposition defined by the leafwise triangulations.

  The constructions of this section   require patience to proceed through the multiple steps required. However, the idea is simple in principle, in that we start with a nice stable transversal defined on a subset of an initial Reeb neighborhood, defined using a coding partition of the transversal $V_1$, then extend this to its intersections with neighboring     Reeb neighborhoods also defined by the coding partition, and observe  that the extension so obtained is well-defined. We then recursively repeat this argument with intersections of the previous results of this procedure, an so obtain a recursive procedure for obtaining a nice stable transversal on all of $\fM$.

\subsection{Reeb neighborhoods from partitions}\label{subsec-reebneighborhoods}

Assume that $\fM$ is a minimal matchbox manifold without holonomy.  Given a clopen subset    $V_1 \subset \fX$ with $\diamX(V_1) < e_H/2$ and 
   $w_0 \in V_1$, let   $L_0$ denote the leaf containing $x_0 = \tau(w_0)$. Let  $\delta_1^*>0$ be the constant defined in the proof of   Theorem~\ref{thm-specialpartition}.
   
 Let    $\cW_{1} = \{W(1;1), \ldots, W(1;\kappa_1)\}$ be a   clopen partition of $V_1$ with $\diamX(W(1;  i))< \delta_1^*$ for all $1 \leq i \leq \kappa_1$ as constructed in   the proof of   Theorem~\ref{thm-specialpartition}, with ``basepoints''   $x_i \in \tau(W(1;i)) \cap L_0$ for which $w_i = \pi(x_i) \in W(1;i)$.
    In addition,    for each   $1 \leq i \leq \kappa_1$,  there is a unique $1 \leq m_i \leq s$ with  $W(1;i) \subset \widehat{W}_{m_i}$, and     
  $K_i$ is  the closed plaque-saturation of  $B_{\F}(x_i, \theta_1 + \dFU)$, and  $\whK_i = \whK_{x_i}$  is defined as in Definition~\ref{def-defbase}.  
 Recall the constant $R_1' = R_1' =  2 \theta_1 + \lF$ was defined by \eqref{eq-RV1}.
   
   Then there are  
 coding partitions derived  from the partition $\cW_{1}$ as in Section~\ref{sec-coding}.
  For each $1 \leq i \leq \kappa_1$,   we have $\ds \cV(1 ; i) = \{ V(1; i,j)  \mid 1 \leq j \leq \kappa_{1, i}\}$ which is a clopen partition of $W(1; i)$ such that the coding function is constant on each $V(1; i,j)$ along paths of length at most $R_1'$.  
Also recall the refinements of these partitions, where the collection  $\cV( 1 ; i,j)$   is a clopen partition of the clopen set $V(1; i,j)$, where for each $1 \leq k \leq \kappa_{1,i,j}$,   the coding function is constant on   $V(1; i,j,k)$ along paths of length at most $R_{1} = 2R_1'$.

Next choose ``basepoints'' for these refined partitions.  
For each $1 \leq i \leq \kappa_{1}$, $1 \leq j \leq \kappa_{1;i}$, let $v(1;i,j) \in V(1;i,j)$ be a point such that $x(1;i,j) = \tau(v(1;i,j)) \in L_0$, 
and let $\whK_{1;i,j}$ be the path-connected subset of  $\fN^{W(1;i)}_{\whK_i}$  containing this point. 
Denote by $\fN_{1;i,j} \subset \fN^{W(1;i)}_{\whK_i}$ the Reeb neighborhood formed from the base $\whK_{1;i,j}$ and transversal $V(1;i,j)$.

Repeat this procedure   for the partitions $\cV( 1 ; i,j)$. 
For each $1 \leq i \leq \kappa_{1}$ and  $1 \leq j \leq \kappa_{1;i}$, let $v(1;i,j,k) \in V(1;i,j,k)$ be a point such that $x(1;i,j,k) = \tau(v(1;i,j,k)) \in L_0$.
Let $\whK_{1;i,j,k}$ be a path-connected subset of $\fN_{1;i,j}$  which contains $x(1;i,j,k)$ in its interior. Denote by $\fN_{1;i,j,k}$,   the   subsets of $\fN_{1;i,j}$  formed from the base $\whK_{1;i,j,k}$ and transversal $V(1;i,j,k)$.
Observe that we then have three  collections of Reeb neighborhoods, each of which forms a covering of $\fM$:
\begin{eqnarray}
\fM ~ & = & ~  \bigcup ~\left\{   \fN_{1;i} \mid  1 \leq i \leq \kappa_{1} \right\} \label{eq-decompositions1} \\
  ~ & = & ~ \bigcup ~\left\{   \fN_{1;i,j} \mid  1 \leq i \leq \kappa_{1} ; 1 \leq j \leq \kappa_{1;i} \right\}  \label{eq-decompositions2}    \\
~ & = & ~   \bigcup ~\left\{  \fN_{1;i,j,k} \mid 1 \leq i \leq \kappa_{1} ; 1 \leq j \leq \kappa_{1;i} ; 1 \leq k \leq \kappa_{1;i,j}\right\} ~ .\label{eq-decompositions3}
\end{eqnarray}
The decomposition \eqref{eq-decompositions1} is that given in the proof of Theorem~\ref{thm-specialpartition}. 
The decomposition \eqref{eq-decompositions2}   arises naturally when considering the intersections of Reeb neighborhoods in \eqref{eq-decompositions1}, while the   decomposition \eqref{eq-decompositions3} is used to show that we obtain a well-defined Cantor foliation. 

There is no natural notion of ``forcing the border'' for the decompositions of $\fM$ that we consider, as is the case for tiling spaces as discussed  in \cite{AP1998,BDHS2010,FHK2002,Sadun2003,Sadun2008}, and as a result we use the additional structure given by the extended coding partitions used to define the spaces in  \eqref{eq-decompositions2} and \eqref{eq-decompositions3} to show that our decompositions are well-defined.

\subsection{Nice stable transversals}\label{subsec-stabletransversals}
 
For $1 \leq i \leq \kappa_1$    and    $z \in \whK_{1;i}  \cap \cT$, there exists  $1 \leq i_z \leq \nu$ such that   $z \in \cT_{i_z} \subset  \oU_{i_z}$. 
Choose a path  in $\whK_{1;i}$ joining $x(1;i) = x_i$ and $z$, which determines a holonomy homeomorphism  $h_{x(1;i),z}$  which will be denoted by $h_z$ for simplicity of notation. Note that $h_z(w_i) \in \fX_{i_z}$ and  define
\begin{equation}\label{eq-transversal1}
W(1;i)^z = h_z(W(1;i)) \subset \fX_{i_z} ~ .
\end{equation}
Then for     $\xi \in \cP_{i_z}(h_z(w_i))$,  a \emph{standard transversal} through $\xi$ is given by 
\begin{equation}\label{eq-transversal1a}
\fZ(\xi,i_z,W(1;i)^z)   \equiv \vp_{i_z}^{-1}\left(\lambda_{i_z}(\xi) , W(1;i)^z \right) ~ 
\end{equation}
which is   a vertical coordinate slice for the coordinates on $U_{i_z}$ which passes through the point $\xi$, and projects to the set $W(1;i)^z$.
Note that if $\xi \in U_{i_{z'}} \cap U_{i_z}$ for some $i_{z'} \ne i_z$ then the section  $\fZ(\xi,i_z,W(1;i)^z)$ need not be a vertical section in  the coordinate chart $U_{i_{z'}}$, and so the choice of the    transversal depends on the index $i_z$ and not just the point $\xi$.

A   \emph{nice stable transversal} $\cX(1;i)$ on the Reeb  neighborhood $\fN_{1;i}$  of $\whK_{1;i} $   is a disjoint union
\begin{equation}\label{eq-transversal2}
\cX(1;i) =\cX_{1}^{1;i} \cup \cdots \cup \cX_{p_i}^{1;i}
\end{equation}
for some ordering of the individual transversals, where each $\cX_{j}^{1;i}$  has the form \eqref{eq-transversal1a}. The choices of the data defining each $\cX_{j}^{1;i}$,  as given in the proof of \cite[Theorem~1.3]{CHL2013a}, are made so that the intersection $\cX(1;i) \cap \whK_{1;i}$ is  a $(d_1,d_2)$-net on $\whK_{1;i}$ where the separation and spanning constants  $0 < d_1 < d_2$ are chosen sufficiently small so that the leafwise nets they induce yield stable leafwise triangulations, as discussed in \cite{CHL2013a}. 
  In particular, we note that the constant $d_2 = \lF^*/5$ satisfies $ d_2 \leq \lF/25$, although    $d_2$ may     be much smaller, as it depends on the  sectional curvatures of the leaves of $\F$, as explained in \cite[Sections~15, 17]{CHL2013a}

The   simplicial decomposition on $\whK_{1;i}$ corresponding to $\cX(1;i)$,  and therefore the transverse Cantor foliation constructed on $\fN_{1;i}$, depends on the ordering of the local sections comprising  a transversal.  In the proof of Proposition \ref{prop-Cantorintersections} below, we make an initial  choice of ordering for the section in $\cX(1;1)$, and then extend this ordering to the other  transversals $\cX(1;i)$ for $1 < i \leq \kappa_1$, using a recursive    procedure. This will yield an ordered  nice stable transversal for all of $\fM$  as a result, which is   used to construct  the transverse Cantor foliation for $\fM$.

We will consider the restriction of the   transversals in  $\cX(1;i)$ to subsets $\fN_{1;i,j}$, and so obtain a partial transversal $\cX(1;i,j)$ for $\fN_{1;i,j}$. To distinguish between these various sections, we will denote by $\cX^{1;i}_{s}$ a section in $\fN_{1;i}$ and by $\cX^{1;i,j}_s$ the corresponding section   for  $\fN_{1;i,j}$ obtained as a restriction of $\cX^{1;i}_s$. We note that the section $\cX^{1;i,j}_s$ is assigned the same number $s$ in the ordering of transverse sections. 

\begin{figure}[H]
\centering
\includegraphics [width=7cm] {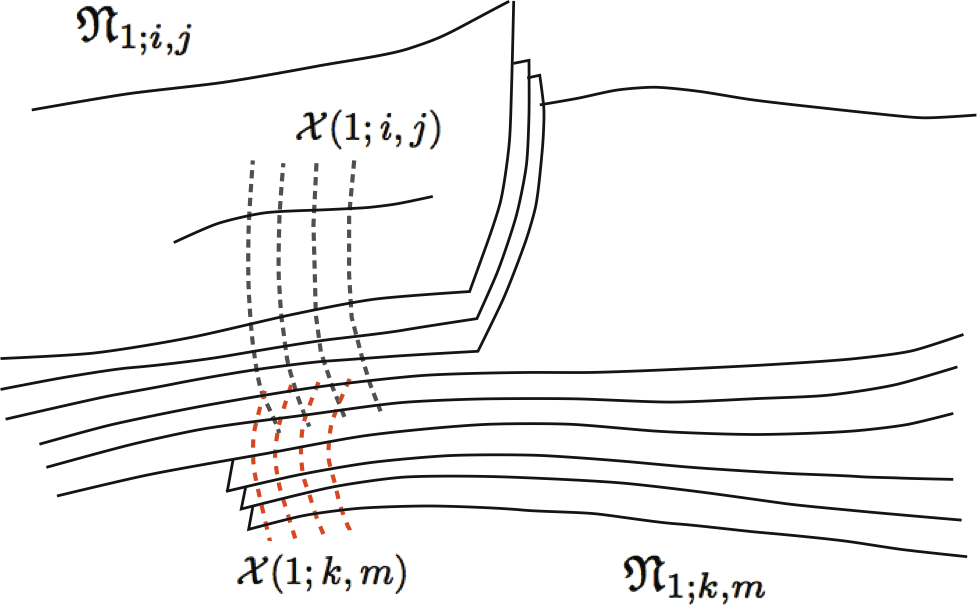}
\caption{Transversals on the overlap of Reeb neighborhoods $\fN_{1;i,j}$ and $\fN_{1;k,m}$ need not match.}
 \label{fig:overlaps3}
\vspace{-20pt}
\end{figure}

 If  the neighborhoods $\fN_{1;i,j}$ and $\fN_{1;k,m}$ intersect, in general the transversals $\cX(1;i,j)$ and $\cX(1;k,m)$ need not match, as illustrated in Figure~\ref{fig:overlaps3}.  The recursive procedure developed in the proof of   Proposition \ref{prop-Cantorintersections} below   yields a collection of transversals   $\cX(1;i,j)$ which match on the intersections of the corresponding Reeb neighborhoods $\fN_{1;i,j}$.
\vspace{-20pt}

\subsection{Constructing nice stable transversals}\label{subsec-existencestabletransversals}
We give a recursive procedure for modifying the   nice stable transversals on the Reeb covering of $\fM$ as in  \eqref{eq-decompositions2} to obtain a global transversal. 

\begin{prop}\label{prop-Cantorintersections}
Nice stable transversals $\cX(1;i,j)$ on Reeb neighborhoods $\fN_{1;i,j}$ associated to the partitions $\cV(1;i)$ can be chosen so that the following property is satisfied: if $x, y \in \fN_{1;i,j} \cap \fN_{1;k,m}$, and $x, y \in \cX^{1;i,j}_s$ and $x \in \cX^{1;k,m}_s$, then $y \in \cX^{1;k,m}_s$. 
In other words, sections of $\cX(1;i,j)$ and $\cX(1;k,m)$ coincide on the intersection $\fN_{1;i,j} \cap \fN_{1;k,m}$.
\end{prop}
\proof 
 The recursive construction of the transversals $\cX(1;i,j)$ follows a multiple step process.   Recall that    $R_1' =  2 \theta_1 + \lF$ was defined by \eqref{eq-RV1} and  $R_1 = 2 R_1'$.

\emph{Step 1.} Choose     nice stable transversals $\ds \cX(1;i) = \{ \cX_{1}^{1;i} ,  \ldots , \cX_{p_i}^{1;i} \}$ on the Reeb  neighborhoods $\fN_{1;i}$ for $1 \leq i \leq   \kappa_{1}$ as in \eqref{eq-transversal2}.
 By assumption, $w_0 \in V(1;1,1)$. Let $\cX(1;1,1)$ be a nice stable transversal for $\fN_{1;1,1}$ obtained by restriction of $\cX(1;1)$.

\emph{Step 2.} The sets $V(1;i,j)$ for $1 \leq i \leq \kappa_1$ and $1 \leq j \leq \kappa_{1;i}$ are defined using the coding function for paths of length at most $R_1'$ as above, and their indices are given the  
 lexicographic order. We next extend the transversals in $\cX(1;1,1)$ to those Reeb neighborhoods which intersect $\fN_{1;1,1}$.

\begin{figure}[H]
\centering
\includegraphics [width=9cm] {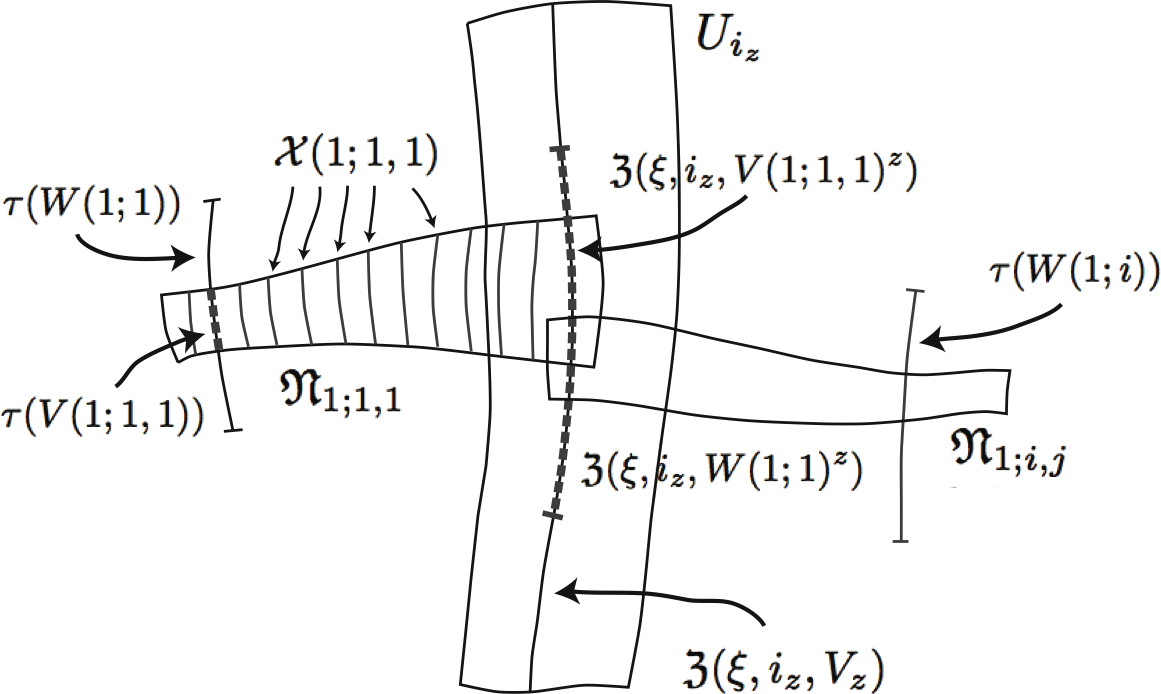}
\caption{A  transversal $\cX(1;1,1)$   constructed on a Reeb neighborhood $\fN_{1;1,1}$. A section $\cX_\xi^{1;1,1} = \fZ(\xi,i_z,V(1;1,1)^z)$ intersects the Reeb neighborhood $\fN_{1;i,j}$ which can be extended to a section $\cX_\xi^{1;i,j}$ since the diameter of $V(1;i,j)$ is so small that $\tau(V(1;i,j)^z) \subset \tau(W(1;1,1)^z) \subset \fZ(\xi, i_z, W(1;1)^z)$.}
 \label{fig:pic1}
\vspace{-10pt}
\end{figure}

 Let $(1;i,j)$ be the least index not equal to $(1;1,1)$ such that  $\fN_{1;1,1} \cap \fN_{1;i,j} \ne \emptyset$, if such exists. 
 For any $\cX_s ^{1;1,1}\in \cX(1;1,1)$ such that $\cX_s^{1;1,1} \cap \fN_{1;i,j} \ne \emptyset$, we define $\cX_s^{1;i,j}$ by
\begin{equation}\label{eq-extendCantor} 
\cX_s^{1;i,j}   = \cX_s^{1;1} \cap \fN_{1;i,j} = \fZ(\xi_s,i_z,W(1;1)^z) \cap \fN_{1;i,j} ~ .
\end{equation}
That is, we extend the transversal to $\fN_{1;i,j}$ while preserving the ordering of the sections.   Figure~\ref{fig:pic1} illustrates the intersections of  the Reeb Neighborhood $\fN_{1;1,1}$ with $\fN_{1;i,j}$. 

\vspace{-10pt} 

We repeat this procedure for every $\fN_{1;i,j}$ which intersects $\fN_{1;1,1}$, following the lexicographical ordering on indices.
Note that  the   transversals constructed in the sets $\cX(1;i,j)$ do not depend on the order in which the neighborhoods $\fN_{1;i,j}$ are considered, since they are all obtained  by restriction of the   transversal $\cX(1;1)$.  
We thus obtain   a   transversal $\cX(1;1,1)$ on $\fN_{1;1,1}$, and  (partial)  transversals $\cX(1;i,j)$ defined on   the neighborhoods $\fN_{1;i,j}$, $1 \leq i \leq \kappa_1$, $1 \leq j \leq \kappa_{1;i}$ which intersect $\fN_{1;1,1}$.

\begin{figure}[H]
\centering
\includegraphics [width=10cm] {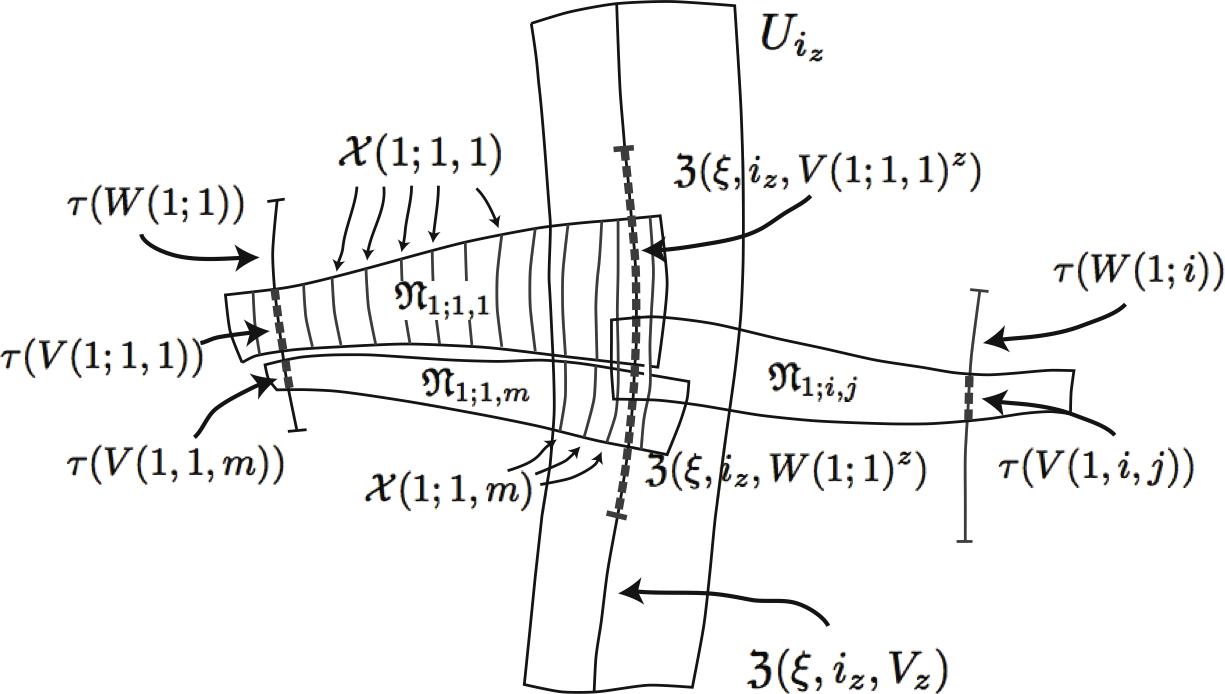}
\caption{A nice stable transversal $\cX(1;1,1)$ is continued to define partial transversals $\cX(1;i,j)$ and $\cX(1;1,m)$.}
 \label{fig:pic2}
\vspace{-10pt}
\end{figure}

Next, we consider the case where there is a Reeb neighborhood $\fN_{1;1,m}$ which does not intersect $\fN_{1;1,1}$, but does intersect a subset of $\fN_{1;i,j}$ on which  a partial transversal has been defined in the previous step. This  situation is illustrated in Figure \ref{fig:pic2}, and involves   a  fundamental point for the construction of the global Cantor foliation on $\fM$. 

Let  $z \in \cX_s^{1;i,j}  = \cX_s^{1;1} \cap \fN_{1;i,j}$.
  Then by the construction of the bases for the Reeb neighborhoods in Definition~\ref{def-defbase}, and by definitions  \eqref{eq-transversal1},  \eqref{eq-transversal1a} and \eqref{eq-transversal2},    there is $x \in \tau(V(1;1))$ and a path  $\gamma$ of length at most $\theta_1 + 2\dFU$ starting at $x$ and ending at  $z$. 
  Let $h_{x,z} \colon V(1;1) \to \cX_s^{1;1}$ be the holonomy map defined by $\gamma$. 
   
Suppose there exists  $v \in   \cX_s^{1;i,j} \cap \fN_{1;1,m}$, then   there is $u \in \tau(V(1;i,j))$ and a path  $\gamma'$ of length at most $\theta_1 + 2\dFU$ starting at $u$ and ending at  $v$.   Let $h_{u,v} \colon V(1;i) \to \cX_s^{1;1}$ be the holonomy map defined by $\gamma'$. 
  Then using the shadowing method in the proof of Lemma~\ref{lem-pathlengths}, we can define a path $(\gamma')^{-1} * \gamma$ with length at most $R' = 2\theta_1 + \lF$, which defines the holonomy map      $h_{u,v}^{-1} \circ h_{x,y}$.
  
  Then $h_{u,v}^{-1}(\pi(v)) = u \in V(1;i,j)$ and $z \in   \fN_{1;i,j}$ implies that $h_{u,v}^{-1}(\pi(z))   \in V(1;i,j)$. 
  This implies that $h_{u,v}^{-1} \circ (\pi(x)) \in  V(1;i,j)$. Since the inverse path $\gamma^{-1} * \gamma'$ has length at most $R_1'$ the definition of the coding partition implies that  $h_{x,y}^{-1} \circ h_{u,v}(V(1;i,j)) \subset W(1;1)$.

  It follows that         $\cX_s^{1;1,m} = \cX_s^{1;1} \cap \fN_{1;1,m}$ defines a transversal for $\fN_{1;1,m}$.

We repeat this procedure for every $\fN_{1;1,m}$ which intersects one or more of Reeb neighborhoods with partial transversals along the subset where the transversal is defined. Since the sections $\cX_s^{1;i,j}$ are obtained by restriction of the same section $\cX^{1;1}_s$, this procedure is well-defined in the case when $\fN_{1;1,m}$ intersects two or more distinct Reeb neighborhoods $\fN_{1;i,j}$ with partial transversal. The partial transversals $\cX(1;1,m)$ do not depend on the order in which the Reeb neighborhoods are considered, and we choose to follow the lexicographic order for definiteness (see Figure \ref{fig:pic2}).

We now have a nice stable transversal $\cX(1;1,1)$ on $\fN_{1;1,1}$, and   transversals $\cX(1;i,j)$ on some of the neighborhoods $\fN_{1;i,j}$, $1 \leq i \leq \kappa_1$, $1 \leq j \leq \kappa_{1;i}$, some of them for $i=1$ and $j \ne 1$.

We must show that the procedure of defining partial transversals eventually stops. For that, let $\gamma, \gamma'$ be homeomorphisms associated to paths of length at most $\theta_1 + 2\dFU$. We notice that if for $u \in V(1;1,1)$ we have $V(1;1,1) \subset \Dom(\gamma)$ and $\gamma(u) \in V(1;i,j)$, then for every $\gamma' $ with $V(1; i,j) \in \Dom(\gamma')$ the code $C^{1,R_1}_{w_0}(\gamma' \circ \gamma)$ is known. It follows that $\gamma^{-1}(\gamma(V(1;1,1)) \cap V(1;i,j))$ is a union of a finite number of sets in the partition $\cV(1;1,1)$ on which the coding function for paths of length at most $R_1$ is constant. By a similar argument, $\gamma(\gamma^{-1}(V(1;i,j)) \cap V(1;1,1))$ is the union of a finite number of clopen sets in the partition $\cV(1;i,j)$. Therefore, an intersection of two slabs $\fN_{1;1,1}$ and $\fN_{1;i,j}$ is always along a subset corresponding to a clopen subset in $V(1;1)$ which is a union of sets in $\cV(1;1,1)$. Similarly, $\fN_{1;1,m}$ intersects $\fN_{1;i,j}$ along a subset corresponding to a clopen subset in $V(1;i,j)$ which is a union of subsets in $\cV(1;i,j)$.
Since the partitions $\cV(1;1)$ and $\cV(1;1,1)$ are finite, and intersections always happen over Reeb neighborhoods within the same range of $1 \leq i \leq \kappa_1$, the procedure of defining a partial transversal stops after a finite number of steps.

 \medskip

\emph{Step 3 (Inductive step)} Suppose for all $(1,1) \leq (i',j') < (i,j)$ a stable transversal has been defined  on $\fN_{1;i',j'}$, and numbers $1, \ldots, p$ were assigned to sections in the nice stable transversals so that on each neighborhood, there is at most one section which is assigned a given number $s$, and two sections are numbered by the same number $s$ in two different neighborhoods if and only if these sections intersect, i.e. $\cX^{1;i,j}_s \cap \cX^{1;k,m}_c \ne \emptyset$ if and only if $s = c$. Then two cases are possible.

\emph{Case 1.} There is a partial transversal $\cX(1;i,j)$ on $\fN_{1;i,j}$. Then use \cite[Proposition 17.4]{CHL2013a} to complete $\cX(1;i,j)$ starting to number the newly defined sections in $\cX(1;i,j)$ from $p+1$. 

We now have to repeat the procedure of defining a partial transversal of \emph{Step 2}. We note that it is enough to consider $(k,m)> (i,j)$, since for all $(k,m) < (i,j)$ such that $\fN_{1;k,m} \cap \fN_{1;i,j} \ne \emptyset$, the transversal $\cX(1;k,m)$ has already been defined, and, moreover, possibly gave rise to a subset of sections in a partial transversal $\cX(1;i,j)$. This ensures that a section $\cX_s^{1;i,j}$ in the transversals $\cX(1;i,j)$ which intersect the Reeb neighborhood $\fN_{1;k,m}$, $(k,m) < (i,j)$, coincides on this intersection with a section $\cX^{1;k,m}_s$, as required.

So suppose there exists $(k,m)>(i,j)$ such that $\fN_{1;k,m} \cap \fN_{1;i,j} \ne \emptyset$. If $\cX_s^{1;i,j} \cap \fN_{1;k,m} \ne \emptyset$, then there exists $z \in K_{1;i,j}$ and a chart $U_{i_z}$ such that 
 \begin{align}\label{eq-conditionextend} 
 \cX_s^{1;i,j} \subset \fZ(\xi_s, i_z, W(1;i)^z) \quad {\rm and}  \quad \gamma(V(1;k,m)) \subset W(1,i) 
 \end{align}
for an appropriate homeomorphism $\gamma$ associated to a path with the same name. Condition \eqref{eq-conditionextend} makes sure that a section $\cX^{1;k,m}_s$ is extended to a maximal subset of $\fN_{1;k,m}$, and in particular to the corresponding subset of $\fN_{1;k,m}$.  We define
  \begin{align*} \cX_s^{1;k,m} & = \fZ(\xi_s, i_z, W(1;i)^z) \cap \fN_{1;k,m}, \end{align*}
and repeat the procedure for every Reeb neighborhood intersecting $\fN_{1;i,j}$.
We may need to repeat the procedure of defining a partial nice stable transversal for $\fN_{1;k',m'}$ such that $\fN_{1;k',m'} \cap \fN_{1;k,m} \ne \emptyset$, and so on. By an argument similar to the one in \emph{Step 2}, the procedure stops after a finite number of repetitions.

\emph{Case 2.} There is no partial transversal on $\fN_{1;i,j}$. Then implement \emph{Step 1} and \emph{Step 2} with $(1;i,j)$ instead of $(1;1,1)$, starting to number sections in $\cX(1;i,j)$ from $p+1$.

Since $\cV(1;i)$ is a finite partition, and $1 \leq i \leq \kappa_\ell$, the procedure stops after a finite number of steps. 
\endproof

\section{Construction  of the transverse Cantor foliation} \label{sec-compatiblecantor}

In this section, we discuss how to obtain the transverse Cantor foliations $\cH_{\ell}$ from the nice stable transversal constructed in Section~\ref{sec-nicestabletransversals}.
The foliation $\cH_1$ is constructed in Section~\ref{subsec-restrictReeb}. Then in  Section~\ref{subsec-Cantorinductive}, we show that $\cH_1$ induces Cantor foliations $\cH_{\ell}$ on families of Reeb neighborhoods associated with smaller coding partitions, which is the key to obtaining the presentation \eqref{eq-mainthmbonding} in Theorem~\ref{thm-main1}.

\subsection{Stable triangulations}\label{subsec-restrictReeb}

The constructions in Section~\ref{subsec-stabletransversals} yield a    nice stable transversal 
\begin{equation}\label{eq-NST}
\cX = \bigcup ~ \left\{\cX(1;i,j)  \mid 1 \leq i \leq \kappa_1 ~, ~ 1 \leq j \leq \kappa_{1,i}  \right\}
\end{equation}
on $\fM$,  whose restriction to each   Reeb  neighborhood $\fN_{1;i,j}$ yields 
  a $(d_1,d_2)$-net for each path-connected component   of $\fN_{1;i,j}$. 
The     net  $\cX \cap \whK_{1;i,j}$  determines a triangulation of $\whK_{1;i,j}$, using the circumscribed sphere method as given in \cite[Sections~12,17]{CHL2013a}. 
The vertices of the triangulation are contained in the net, and the transverse stability property for $\cX$ implies that the triangulations so obtained are isomorphic for each path-connected component   of $\fN_{1;i,j}$. As remarked earlier,   the   leaves of a Cantor foliation on $\fN_{1;i,j}$ are defined by identifying points with the same barycentric coordinates in this simplicial decomposition. This procedure is described in detail in   \cite[Section~11]{CHL2013a}. 

In this section, we take care of a technical issue that arises in this construction of the Cantor foliation. Since taking barycentric coordinates in a simplex only makes sense if all vertices of a simplex are contained in $\fN_{1;i,j}$ we have to restrict to a subset of $\fN_{1;i,j}$. In addition, we have to make sure that after shrinking, the new Reeb neighborhoods still cover $\fM$. This is achieved if for each $w \in V(1;i,j)$ the restricted neighborhood contains the leafwise ball $B_{\F}(\tau(w), \theta_1)$.

\begin{prop}\label{prop-honeycombReeb}
A Reeb neighborhood $\fN_{1;i,j}$ associated to a clopen set $V(1;i,j)$ contains a foliated neighborhood $\widehat{\fN}_{1;i,j}$ with the following properties.
\begin{enumerate}
\item for every $w \in V(1;i,j)$ we have $B_{\F}(\tau(w), \theta_1) \subset \widehat{\fN}_{1;i,j}$.
\item $\widehat{\fN}_{1;i,j} \cap K_{1;i,j}$ is the union of simplices of maximal dimension $n$ for the simplicial decomposition of $K_{1;i,j}$ associated to the nice stable transversal $\cX(1;i,j)$.
\item $\widehat{\fN}_{1;i,j}$ is transversely stable; that is, for every path-connected component $K \subset \fN_{1;i,j}$ with $\tau(w) \in K$,  the points $\cX^{1;i,j}_{i_0}\cap K, \ldots, \cX^{1;i,j}_{i_p} \cap K$ are vertices of a $p$-simplex for $1 \leq p \leq n$, if and only if $\cX^{1;i,j}_{i_0}\cap K_{1;i,j}, \ldots, \cX^{1;i,j}_{i_p} \cap K_{1;i,j}$ are vertices of a $p$-simplex.
\end{enumerate}
\end{prop}

\proof 
The intersection $\cX \cap \whK_{1;i,j}$   is  a $(d_1,d_2)$-net for $\whK_{1;i,j}$, with $d_2  \leq \lF/25$.
Thus, for each $y \in \cX(1;i,j)$ with $w = \pi(y)$, there is a Voronoi cell $\cC^{\cX}_{w}(y) \subset D_{\F}(y,d_2)$. 
A simplex in the simplicial decomposition associated to the ordered transversal $\cX$ which has $y$ as one of the vertices, is contained in the star-neighborhood $\cS^{\cX}_{w}(y)$. Then  by Lemma~\ref{lem-star2} for this case,  we have  $\ds \cS^{\cX}_{w}(y) \subset D_\F(y,3d_2)$.

Therefore, for each simplex such that   one of its vertices is contained in a ball $B_{\F}(y,\theta_1 +\dFU)$, the entire simplex  is contained in the ball  $B_{\F}(y,\theta_1 +\dFU + 4d_2)$. From the estimates $\dFU < \lF/5$ and $d_2   \leq \lF/25$ we have $\dFU + 4d_2 < \lF/2$ so that
$B_{\F}(y,\theta_1 +\dFU + 4d_2) \subset B_{\F}(y,\theta_1 +\lF/2)$.

For  $v(1;i,j) \in V(1;i,j)$ with $y = x(1;i,j) = \tau(v(1;i,j))$,  the ball $B_{\F}(y, \theta_1 +\dFU) \subset \whK_{1;i,j}$ where
 $ \whK_{1;i,j}$ is the closure of the plaque-saturation of $B_\F(y, \theta_1 + \lF)$.

Let $K^t_{1;i,j}$ denote the closed triangulated set obtained from the union of all simplices of maximal dimension $n$ associated to the net $\cX \cap \whK_{1;i,j}$ which intersect $K_{1;i,j}$. It then follows that $\ds K^t_{1;i,j} \subset \whK_{1;i,j}$ so the map $\Phi_i$ given in Theorem~\ref{thm-specialpartition} is defined on  
$\ds K^t_{1;i,j} \times V(1;i,j)$. 

Set $\ds \widehat{\fN}_{1;i,j} = \Phi_i \left\{ K^t_{1;i,j} \times V(1;i,j) \right\}$ for $1 \leq i \leq \kappa_1$ and $1 \leq j \leq \kappa_{1;i}$. 

The inclusions $B_{\F}(x(1;i,j) , \theta_1 +\dFU) \subset K_{1;i,j}$ implies that these sets form a cover of $\fM$, exactly as in the proof of Theorem~\ref{thm-specialpartition}.
 The transverse stability of the sets $\widehat{\fN}_{1;i,j}$ follows from the transverse stability of the transversals $\cX(1;i,j)$.   
\endproof

We have a clopen set $V_1$, and a clopen partition $V(1;i,j)$ of $V_1$ with associated system of Reeb neighborhoods $\whfN_{1;i,j}$, $1 \leq i \leq \kappa_{1}$, $1 \leq j \leq \kappa_{1;i}$ defined in Proposition \ref{prop-honeycombReeb}. Since each path-connected component of $\whfN_{1;i,j}$ is a union of simplices, and the simplicial decomposition is stable in the transverse direction, the following is well-defined.

\begin{defn}\label{cantor-onreeb}
A Cantor foliation $\cH_{1;i,j}$ on $\whfN_{1;i,j}$ is defined as follows: two points $x,y \in \whfN_{1;i,j}$ are equivalent, $x \approx_{1;i,j} y$ if and only if the following conditions are satisfied:
\begin{enumerate}
\item if $x \in \Delta(x_0,x_1, \ldots, x_p)$ where $\Delta(x_0,\ldots,x_p)$ is a $p$-simplex with vertices $x_k \in \cX^{1;i,j}_{i_k}$, then $y \in \Delta(y_0,\ldots,y_p)$ with $y_k \in \cX^{1;i,j}_{i_k}$, the same section of $\cX(1;i,j)$,
\item $x$ and $y$ have the same barycentric coordinates in $\Delta(x_0, \ldots,x_p)$ and $\Delta(y_0,\ldots,y_p)$ respectively.
\end{enumerate}
For a point $x \in \fM$, $\cH_{1;i,j}(x)$ denotes a leaf of $\cH_{1;i,j}$ through $x$.
\end{defn}

If two neighborhoods $\whfN_{1;i,j}$ and $\whfN_{1;k,m}$ intersect, then the foliations $\cH_{1;i,j}$ and $\cH_{1;k,m}$ match on the intersection, since the sections of transversals $\cX(1;i,j)$ and $\cX(1;k,m)$ are continuations of each other, and the ordering of sections in nice stable transversals matches as well. We extend ${\cH_{1;i,j}}$ on $\fM$ in the following way.

\begin{defn}\label{extended-equiv}
Let $\fM$ be a matchbox manifold, and $\whfN_{1;i,j}$ be a Reeb neighborhood with a Cantor foliation $\approx_{{1;i,j}}$. We define $\approx_{1;i,j}$ on $\fM$ in the following way: let $x,y \in \fM$. Then $x \approx_{1;i,j} y$ if either $x,y \in \whfN_{1;i,j}$ and $x \approx_{{1;i,j}} y$, or $x,y \notin \whfN_{1;i,j}$ and $x = y$.

\end{defn}

Thus the equivalence classes of $\approx_{1;i,j}$ coincide with leaves of ${\cH_{1;i,j}}$ on $\whfN_{1;i,j}$, and points outside of $\whfN_{1;i,j}$ have only themselves in their equivalence class.

\begin{defn}\label{cantor-fol-comp}
Let $\fM$ be a matchbox manifold, and let $\approx_{1;i,j}$, $1 \leq i \leq \kappa_1$, $1 \leq j \leq \kappa_{1;i}$ be a family of equivalence relations on $\fM$. For any $x,y \in \fM$ define $x \approx_1 y$ if and only if there exists a finite collection of points $z_1,z_2,\ldots, z_k$ such that
  \begin{align*} x & \approx_{1;i,j} z_1  \approx_{1;i_1,j_1} z_2  \approx  \cdots  \approx z_k \approx_{1;i_k,j_k} y \end{align*}
A Cantor foliation $\cH_1$ has equivalence classes of $\approx_1$ as leaves. For a point $x \in \fM$, $\cH_1(x)$ denotes a leaf of $\cH_1$ through $x$.
\end{defn}

It is straightforward that $\approx_{1}$ is reflexive, symmetric and transitive. Geometrically, this relation identifies $x \in \whfN_{1;i,j}$ and $y \in \whfN_{1;k,m}$ if and only if the leaves $\cH_{1;i,j}(x), \cH_{1;i_1,j_1}(z_1), \ldots, \cH_{1;i_k,j_k}(z_k)$ form an intersecting chain, and $\cH_{1;i_k,j_k}(z_k) \owns y$. The equivalence relation $\approx_1$ determines the Cantor foliation  $\cH_1$.


Thus, by Propositions~\ref{prop-Cantorintersections} and \ref{prop-honeycombReeb},  there is a system of Reeb neighborhoods $\whfN_{1;i,j}$, $1 \leq i \leq \kappa_1$, $1 \leq j \leq \kappa_{1;i}$, and a Cantor foliation $\cH_{1;i,j}$ on each $\whfN_{1;i,j}$ such that the Cantor foliations match on the intersections of neighborhoods, and their union equals $\fM$. 
These  Cantor foliations combine to yield the Cantor foliation $\cH_1$  on $\fM$.

 \medskip
 
\subsection{Families of nested Cantor foliations}\label{subsec-Cantorinductive}

  In this section, we give an inductive procedure for restricting the foliation $\cH_1$ to the Reeb neighborhoods $\fN_{\ell;i,j}$ defined using the partitions defined in Proposition \ref{prop-dynamicpart}.
   
 \eject
   
 \begin{prop}\label{prop-nextedfoliations}
 There exists a   nested sequence of Cantor foliations $\cH_\ell$, for all $\ell \geq 1$, with $\cH_1$ given as above, such that each Reeb neighborhood $\fN_{\ell;i,j}$ associated to a clopen set $V(\ell;i,j)$ contains a foliated neighborhood $\widehat{\fN}_{\ell;i,j}$ with the following properties.
\begin{enumerate}
\item for every $w \in V(\ell;i,j)$ we have $B_{\F}(\tau(w), \theta_{\ell}) \subset \widehat{\fN}_{\ell;i,j}$.
\item $\widehat{\fN}_{\ell;i,j} \cap K_{\ell;i,j}$ is the union of simplices of maximal dimension $n$ for the simplicial decomposition of $K_{\ell;i,j}$ associated to the nice stable transversal $\cX(\ell;i,j)$.
\item $\widehat{\fN}_{\ell;i,j}$ is transversely stable; that is, for every path-connected component $K \subset \fN_{\ell;i,j}$ with $\tau(w) \in K$,  the points $\cX^{\ell;i,j}_{i_0}\cap K, \ldots, \cX^{\ell;i,j}_{i_p} \cap K$ are vertices of a $p$-simplex for $1 \leq p \leq n$, if and only if $\cX^{\ell;i,j}_{i_0}\cap K_{\ell;i,j}, \ldots, \cX^{\ell;i,j}_{i_p} \cap K_{\ell;i,j}$ are vertices of a $p$-simplex.
\item Each $\widehat{\fN}_{\ell;i,j}$ has a  Cantor foliation $\cH_{\ell;i,j}$ such that they defined a Cantor foliation $\cH_{\ell}$ on $\fM$.
\item The identity map on $\fM$ induces an inclusion of leaves of $\cH_{\ell}$ into the leaves of $\cH_{\ell -1}$.
\item For each $x \in \fM$ and $\ell \geq 1$, the leaves of $\cH_{\ell}$   containing $x$  are nested    Cantor sets.
  Moreover,   $\diamM(\cH_{\ell}(x)) \leq \whdelta_{\ell}$.
\end{enumerate}
 \end{prop}

\proof
 The idea of the inductive construction of $\cH_{\ell}$ is to restrict the given foliation $\cH_{\ell-1}$ on the Reeb neighborhoods $\whfN_{\ell -1;i,j}$ to the Reeb neighborhoods $\fN_{\ell;i,j}$, and then modify the neighborhoods $\fN_{\ell;i,j}$ to obtain the neighborhoods $\whfN_{\ell;i,j}$, 
 as in the proof of Proposition~\ref{prop-honeycombReeb}. The procedure is straightforward, though care is needed, to  show that the    restricted Cantor foliations are nested.

Let $\ell>1$, and assume there is given a   Cantor foliation  $\cH_{\ell-1}$ on $\fM$ satisfying the conditions of Proposition~\ref{prop-nextedfoliations} for $\ell -1$.

We also assume there is given a collection of nested clopen sets as given by Proposition \ref{prop-dynamicpart}.  
Recall that $\tau(V_\ell) \cap L$ is a $(\lambda_1(\e_\ell), \theta_\ell)$-net for each leaf $L \subset \fM$, and in particular for the leaf $L_0$. 
Define the constants $R_{\ell} = 4 \theta_{\ell} +   2\lF$ and $R_{\ell}' =  2 \theta_{\ell} + \lF$.
Then we have Reeb neighborhoods which cover $\fM$, defined as in Section~\ref{subsec-stabletransversals} using this data, so that 
\begin{eqnarray}
\fM ~  & = & ~ \bigcup ~\left\{   \fN_{\ell;i,j} \mid  1 \leq i \leq \kappa_{\ell} ; 1 \leq j \leq \kappa_{\ell;i} \right\} \label{eq-decompositions4} \\
~ & = & ~   \bigcup ~\left\{  \fN_{\ell;i,j,k} \mid 1 \leq i \leq \kappa_{\ell} ; 1 \leq j \leq \kappa_{\ell;i} ; 1 \leq k \leq \kappa_{\ell;i,j}\right\} ~ .\label{eq-decompositions5}
\end{eqnarray}

 For each $1 \leq i \leq \kappa_\ell$ and $1 \leq j \leq \kappa_{\ell;i}$, the inclusion    $\fN_{\ell;i,j} \subset  \fM$ defines  the foliation $\cH_{\ell}$ by restriction. However, it is necessary that this restriction is of the form in Definition~\ref{def-cantorfol}. This is guaranteed by the following technical result.

\begin{lemma}\label{lem-nestedslabs}
For each $1 \leq i \leq \kappa_\ell$ and $1 \leq j \leq \kappa_{\ell;i}$, with
\begin{equation}
\fN_{\ell;i,j} \subset  \bigcup \fN_{\ell-1;k,m}
\end{equation}
where $k$ and $m$ run over a subset of the indexing sets of $\cW(\ell-1)$ and $\cV(\ell-1;k)$, the following condition is satisfied:

{\bf (A)}  For each $z \in K_{\ell;i,j} \cap \cT$, if  
   $\ds \fZ(0, i_z, V(\ell;i,j)^z) \cap \fN_{\ell-1;k,m} \ne \emptyset$,  then $\fZ(0, i_z, V(\ell;i,j)^z) \subset \fN_{\ell-1; k,m}$.

 \end{lemma}

\proof 
Since the sets of $\cV(\ell;i)$, $1 \leq i \leq \kappa_\ell$, partition $V_\ell = V(\ell-1,1,1) $, the condition $(A)$ is satisfied for each $V(\ell;i,j)$. Also this implies that for every $z \in K_{\ell-1;1,1}$ we have $V(\ell;i,j) \subset V_\ell \subset \Dom(h_z)$, which means that $(A)$ is satisfied for all $z' \in K_{\ell;i,j} \cap \fN_{\ell-1;1,1} \cap \cT$.

Suppose $\fN_{\ell;i,j} \cap \fN_{\ell-1;k,m} \ne \emptyset$. We assume for now that $\fN_{\ell;i,j}$ has a single intersection with $\fN_{\ell-1;k,m}$, i.e. for each path-connected component $K$ of $\fN_{\ell;i,j}$ the intersection $K \cap \fN_{\ell-1;k,m}$ is path-connected. 

Let $z \in K \cap \fN_{\ell-1;k,m} \subset K'$, where $K'$ is a path-connected component of $\fN_{\ell-1;k,m}$.
Recall from conditions \eqref{eq-nestedreeb} and  \eqref{eq-whdelta2} that the restriction $\diamX(V(\ell;i,j))< \whdelta_\ell$ ensures that for every translation of $V(\ell;i,j)$ along a path $\gamma$ of length at most $R_\ell$, points in $V(\ell;i,j)$ do not move apart further than the constant $\ve_\ell$, which was chosen to be smaller than the distance between any two sets in the partition $\cV(\ell-1;k)$. This means that if $\gamma(V(\ell;i,j))$ hits $V_{\ell-1}$, then the image $\gamma(V(\ell;i,j))$ is contained in one of the sets of the clopen partition $\cV(\ell-1;k)$ of $V_{\ell-1}$. 

Let $\gamma'$ be a path between the point $K' \cap \tau(V(\ell-1;k,m))$ and $z$. Then 
$$\gamma(V(\ell;i,j)) \subset  V(\ell-1;k,m) \subset \Dom(\gamma')$$ 
and the condition $(A)$ is satisfied for the section $\tau(V(\ell;i,j)^z) = \tau(\gamma'(\gamma(V(\ell;i,j))))$. 

\begin{figure}[H]
\centering
\includegraphics [width=8cm] {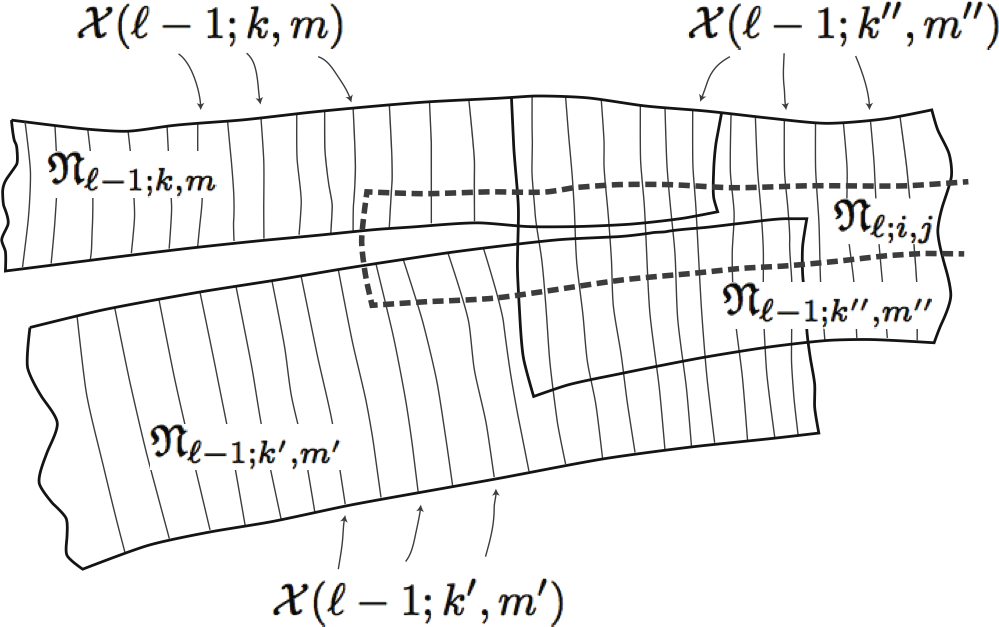}
\caption{$\fN_{\ell;i,j}$ cannot intersect the neighborhoods $\fN_{\ell-1;k,m}$ and $\fN_{\ell-1;k',m'}$ as above.}
 \label{fig:Reebnested}
\vspace{-10pt}
\end{figure}

Next, suppose that  $\fN_{\ell;i,j}$ has multiple intersections with $\fN_{\ell-1;k,m}$, as illustrated in Figure~\ref{fig:Reebnested}, where $\fN_{\ell;i,j}$ intersects 
Reeb neighborhoods $\fN_{\ell-1;k,m}$,$\fN_{\ell-1;k',m'}$ and $\fN_{\ell-1;k'',m''}$.

As the bases of the Reeb neighborhoods are plaque saturated,   by compactness the number of possible intersections is finite. That is, for each path-connected component $K \subset \fN_{\ell;i,j}$ the set $K \cap \fN_{\ell-1;k,m}$ has a finite number of path-connected components. Then the proof is completed by repeating  the above argument for each intersection.
\endproof
 
Lemma~\ref{lem-nestedslabs} allows us to define a nice stable transversal on each Reeb neighborhood $\fN_{\ell;i,j}$ as follows.

\begin{defn}\label{defn-CantorHell}
Let $\cX(\ell-1;k,m)$ be a family of nice stable transversals defined on Reeb neighborhoods $\fN_{\ell-1;k,m}$ associated to the coding partition $\cV(\ell-1;k)$, and let $\fN_{\ell;i,j}$ be a collection of Reeb neighborhoods associated to the coding partition $\cV(\ell;i)$. 

Define $\cX(\ell;i,j)$ as follows: for each $\cX_\xi^{\ell-1;k,m}$ such that $\cX_\xi^{\ell-1;k,m} \cap \fN_{\ell;i,j} \ne \emptyset$, define
  \begin{align*} \cX_\xi^{\ell;i,j} & = \cX_\xi^{\ell-1;k,m} \cap \fN_{\ell;i,j}. \end{align*}
\end{defn}

Now use argument similar to the proof of Proposition~\ref{prop-nextedfoliations} to restrict each Reeb neighborhood $\fN_{\ell;i,j}$ to a Reeb neighborhood $\whfN_{\ell;i,j}$ where each connected component is the union of simplices in the simplicial decomposition associated to the transversal $\cX_{\ell}$.
Then use an argument similar to those  in Section~\ref{subsec-restrictReeb} to define a family of equivalence relations $\approx_{\ell;i,j}$ on $\fM$, and to obtain  an equivalence relation $\approx_\ell$ and a Cantor foliation $\cH_\ell$ on $\fM$.

The claims (5) and (6) in Proposition~\ref{prop-nextedfoliations}  follow from the construction of $\cH_{\ell}$.
\endproof

\section{Inverse limit approximations of matchbox manifolds}\label{sec-Cantorprojections}

In this section we complete the proof of Theorem~\ref{thm-main1}. Namely, having established a sequence of nested equivalence relations $\approx_\ell$ on a matchbox manifold $\fM$, we prove that each quotient $M_{\ell} = \fM/\approx_\ell$ is a \emph{topological branched manifold}, and 
$\ds \fM \cong \varprojlim \, \{q_{\ell' \ell }: M_{\ell'} \to M_\ell\}$, where $q_{\ell',\ell}$ is map induced by inclusions of equivalence classes of $\approx_{\ell'}$ into equivalence classes of $\approx_{\ell}$.

\begin{thm}\label{branched}
Let $\fM$ be a minimal matchbox manifold with a Cantor foliation $\cH_\ell$ and associated equivalence relation $\approx_\ell$. Then the quotient $M_\ell = \fM/\approx_\ell$ is a compact connected metrizable space. For each $x \in \fM$, an open neighborhood of $x$ is homeomorphic to a finite union of open disks in $\mathbb{R}^n$ modulo identifications. Thus,  $M_\ell$ is a \emph{topological branched manifold}.
\end{thm}

\proof We first show that the quotient $M_\ell = \fM/\approx_\ell$ is a compact connected metrizable space.

Since equivalence classes of $\approx_\ell$ are closed, for every closed set $U \subset \fM$ its saturation $S(U)$ by equivalence classes of $\approx_\ell$ is closed. Therefore,  $\approx_{\ell}$ is a closed equivalence relation \cite{Bourbaki1989}.

Since $\fM$ is compact and Hausdorff, it is locally compact and so it is completely regular \cite[p.42]{Kur1968}. Since $\approx_\ell$ is a closed equivalence relation, its graph is closed in $\fM \times \fM$ \cite[p.82]{Bourbaki1989}. It follows that $M_\ell$ is Hausdorff. The quotient $M_\ell$ is compact and connected since $\fM$ is compact and connected. It is a continuous image of a compact metric space $\fM$, and so is metrizable.

Since each $x \in \fM$ is contained in at most a finite number of Reeb neighborhoods $\whfN_{\ell;i,j}$, a point $[x]_\ell \in M_\ell$ has a neighborhood which is a finite union of open disks modulo identifications. Although leaves of $\fM$ are Riemannian manifolds, there is no canonical way to put a differentiable structure on $M_\ell$, and therefore, it is a topological branched manifold (compare \cite{Williams1974}).
\endproof

Denote by $q_\ell:\fM \to \fM/\approx_\ell$ the projection map.
Let $\ell' > \ell$, and $M_{\ell'} = \fM/\approx_{\ell'}$ and $M_\ell = \fM/\approx_\ell$ be quotient branched manifolds. We are going to define bonding maps between them.

\begin{lemma}\label{bondingmaps}
For $\ell' > \ell$, there is a natural continuous map $q_{\ell' \ell} \colon M_{\ell'} \to M_{\ell}$ such that 
$q_{\ell'} = q_{\ell} \circ q_{\ell', \ell}$.
\end{lemma}

\proof Note that for $\ell' > \ell $ we have the following diagram
  \begin{equation*} \xymatrix{ \fM \ar[rr]^{\rm id} \ar[d]_{q_{\ell'}} & & \fM \ar[d]^{q_\ell} \\ M_{\ell'} \ar@{-->}[rr]_{q_{\ell'\ell}} && M_\ell} 
\end{equation*}
so we set $q_{\ell, \ell'} \colon \fM /\approx_{\ell'} \to \fM /\approx_{\ell}: \bar{x}' = [x]_{\ell'} \to [x]_\ell$. By Lemma~\ref{lem-nestedslabs} the equivalence classes of $\approx_{\ell'}$ are subsets of equivalence classes of $\approx_\ell$, and so this map is well-defined. It is straightforward that $q_{\ell' \ell}$ is continuous.
\endproof

Recall that the inverse limit of the directed systems of maps
$\ds
\left\{ q_{\ell,\ell+1} \colon M_{\ell+1} \to M_{\ell} \mid \ell \geq \ell_0 \right\}
$
is the topological space
\begin{equation}
\cS\{q_{\ell,\ell'} \colon M_{\ell'} \to M_{\ell}\} ~ \equiv ~  \{ \omega = (\omega_{\ell_0} , \omega_{\ell_0 +1} , \ldots) \in \prod_{\ell \geq \ell_0}^{\infty} ~ M_{\ell} ~ | ~ q_{\ell, \ell+1}(\omega_{\ell +1}) = \omega_{\ell}  \}
\end{equation}

\medskip

\begin{thm}\label{prop-invlim}
There is a homeomorphism $q \colon \fM \to \cS\{q_{\ell,\ell'} \colon M_{\ell'} \to M_{\ell}\}$ of foliated spaces.
\end{thm}
\proof
For $x \in \fM$ define
$$q(x) = ([x]_{1} , [x]_{2}, \ldots ) \in \cS\{q_{\ell,\ell'} \colon M_{\ell'} \to M_{\ell}\}$$
which is well-defined by Lemma~\ref{bondingmaps}. The map to each factor, $x \mapsto [x]_{\ell} = q_{\ell}(x)$ is continuous, hence $q$ is continuous. 

Let $x, y \in \fM$ such that $q(x) = q(y)$.
Then $ q_{\ell}(x) = q_{\ell}(y)$ for all $\ell \geq 1$.
That is, $x \approx_{\ell} y$ for all $\ell \geq 1$. 

This means that $x$ and $y$ are always in the same leaf of $\cH_{\ell}$.  Proposition~\ref{prop-nextedfoliations}.6 shows  that the maximal diameter of $\cH_{\ell}$-equivalence classes    tend to zero as $\ell \to \infty$. It follows that $x$ and $y$ cannot be separated by an open set, and, therefore, $x=y$.

Surjectivity of $q$ follows again from the fact that equivalence classes of $\approx_\ell$ are nested, and their maximal length tends to zero.
\endproof

This completes the proof of Theorem~\ref{thm-main1}.


\end{document}